\documentclass[12pt]{amsart}
\usepackage{amsfonts,amssymb}
\RequirePackage{ifpdf}
\usepackage{amsmath}
\usepackage{color}
\usepackage{pstricks,pst-node}
\usepackage{pst-plot}
\def\appendix#1{
\addtocounter{section}{1} \setcounter{equation}{0}
\renewcommand{\thesection}{\Alph{section}}
\section*{Appendix \thesection\protect\indent\quad
#1}
}
\renewcommand{\theequation}{\thesection.\arabic{equation}}

\catcode`\@=11
\def\marginnote#1{}

\newcount\hour
\newcount\minute
\newtoks\amorpm
\hour=\time\divide\hour by60 \minute=\time{\multiply\hour by60
\global\advance\minute by-\hour}
\edef\standardtime{{\ifnum\hour<12 \global\amorpm={am}%
        \else\global\amorpm={pm}\advance\hour by-12 \fi
        \ifnum\hour=0 \hour=12 \fi
        \number\hour:\ifnum\minute<10 0\fi\number\minute\the\amorpm}}
\edef\militarytime{\number\hour:\ifnum\minute<100\fi\number\minute}

\newcommand{\tcr}{\textcolor{red}}
\newcommand{\tcb}{\textcolor{blue}}
\newcommand{\tcg}{\textcolor{green}}
\newcommand{\tcw}{\textcolor{white}}

\def\draftlabel#1{{\@bsphack\if@filesw {\let\thepage\relax
      \xdef\@gtempa{\write\@auxout{\string
          \newlabel{#1}{{\@currentlabel}{\thepage}}}}}\@gtempa \if@nobreak
    \ifvmode\nobreak\fi\fi\fi\@esphack} \gdef\@eqnlabel{#1}}
    \def\@eqnlabel{}
\def\@vacuum{}
\def\draftmarginnote#1{\marginpar{\raggedright\scriptsize\tt#1}}

\def\draft{
  \oddsidemargin -.5truein
  \def\@oddfoot{\footnotesize \sl preliminary draft \hfil
    \rm\thepage\hfil\sl\today\quad\militarytime}
  \let\@evenfoot\@oddfoot \overfullrule 3pt
    \let\label=\draftlabel
    \let\marginnote=\draftmarginnote
  \def\@eqnnum{(\theequation)\rlap{\kern\marginparsep\tt\@eqnlabel}%
    \global\let\@eqnlabel\@vacuum}

  }
\voffset= - 0.5in \hoffset= - 1.0in

\newcommand{\tr}{\,{\rm Tr}\,}

\def\be{\begin{equation}}
\def\ee{\end{equation}}
\def\bea{\begin{eqnarray}}
\def\eea{\end{eqnarray}}
\def\<{\langle}
\def\>{\rangle}

\def\tr{\mathop{\rm{tr}}}

\def\ocomma{{\phantom{\Bigm|}^{\phantom {X}}_{\raise-1.5pt\hbox{,}}\!\!\!\!\!\!\otimes}}

\newtheorem{theorem}{Theorem}[section]
\newtheorem{lemma}[theorem]{Lemma}
\newtheorem{proposition}[theorem]{Proposition}
\theoremstyle{definition}
\newtheorem{definition}[theorem]{Definition}

\newtheorem{remark}[theorem]{Remark}

\allowdisplaybreaks

\long\def\rem#1{}

\begin{document}

\title[Fenchel--Nielsen coordinates and Goldman brackets]
{Fenchel--Nielsen coordinates and Goldman brackets}
\author{Leonid O. Chekhov$^{\ast}$}\thanks{$^{\ast}$Steklov Mathematical
Institute of Russian Academy of Sciences, Moscow, Russia. Email: chekhov@mi-ras.ru. This work was performed at the Steklov International Mathematical Center and supported by the Ministry of Science and Higher Education of the Russian Federation (agreement no. 075-15-2019-1614).}

\begin{abstract}
We explicitly show that the Poisson bracket on the set of shear coordinates introduced by V.V. Fock in 1997 induces the Fenchel--Nielsen bracket on the set of gluing parameters (length and twist parameters) for pairs of pants decomposition for Riemann surfaces with holes $\Sigma_{g,s}$. We generalize these structures to the case of Riemann surfaces $\Sigma_{g,s,n}$ with holes and bordered cusps.
\end{abstract}

\maketitle

\section{Introduction}\label{s:intro}
\setcounter{equation}{0}

Constructing Darboux coordinates  for moduli spaces of Riemann surfaces $\Sigma_{g,s}$ of genus $g\ge 0$ with $(s>0)$, or without (s=0), holes has a long and successful history. It is intrinsically related to the Poincar\'e uniformization of Riemann surfaces in which they are represented as quotients of the hyperbolic upper half-plane $\mathbb H_2^+$ under the action of a finitely generated discretely acting Fuchsian subgroup $\Delta_{g,s}$ of $PSL(2,\mathbb R)$: $\Sigma_{g,s}=\mathbb H_2^+/ \Delta_{g,s}$. Parameterizations of this action for a given $g$ and $s$ are called Teichm\"uller spaces ${\mathfrak T}_{g,s}$, and the action of $\Delta_{g,s}$ can be naturally lifted to mapping-class group transformations of ${\mathfrak T}_{g,s}$, so we expect any natural structure on a Teichm\"uller space to be consistent with the Fuchsian group action.

Historically, the first set of symplectic coordinates were Fenchel--Nielsen length--twist coordinates $\{ \ell_{A_i}, \tau_{B_i}\}_{i=1}^{3g-3+s}$ based on pair-of-pant decompositions of $\Sigma_{g,s}$; the corresponding symplectic form was merely $\sum_{i=1}^{3g-3+s}  d\ell_{A_i}\wedge d\tau_{B_i}$. Scott Wolpert  then used\cite{Wlp82}, \cite{Wlp83} the Fenchel--Nielsen form to derive a bracket between geodesic functions $G_\gamma=e^{\ell_\gamma/2}+e^{-\ell_\gamma/2}$ corresponding to two intersecting geodesic lines $\gamma_1$ and $\gamma_2$ on $\Sigma_{g,s}$: it was given locally by the sum over points of intersections of these two geodesics,
\be\label{Wolpert}
\{ G_{\gamma_1},G_{\gamma_2}\}=\sum_{P\in\gamma_1\# \gamma_2} \cos\varphi_{P},
\ee
where $\varphi_{P}$ is the signed angle between the corresponding geodesic lines at their crossing at the point $P$.
Almost simultaneously, William Goldman introduced \cite{Gold} his celebrated bracket on the set of geodesic functions,
\be\label{Goldman}
\{ G_{\gamma_1},G_{\gamma_2}\}=\sum_{P\in\gamma_1\# \gamma_2}\frac12 \bigl(G_{\gamma_1 \circ_P \gamma_2}-G_{\gamma_1 \circ_P \gamma_2^{-1}} \bigr),
\ee
where two geodesic functions in the right-hand side correspond to geodesic lines $\gamma_1 \circ_P \gamma_2$ and $\gamma_1 \circ_P\gamma_2^{-1}$ (possibly with self-intersections) obtained by resolving the crossing at the point $P$ in two possible ways (see Fig.~\ref{fi:skein}); throughout this text we often use the fact that every homotopy class of closed curves in $\Sigma_{g,s}$ contains a unique closed geodesics having the minimum length in the hyperbolic geometry. A clear advantage of the Goldman bracket is that it is manifestly mapping-class group invariant since both geodesic lengths and homotopy relations are invariant under the action of the mapping class group.

The proof that the brackets  (\ref{Wolpert})  and (\ref{Goldman})  give the same answer is a relatively easy exercise in hyperbolic geometry, see, e.g., \cite{McShane}, which we reproduce in the Appendix to this paper.

The inverse statement, that is, that the Goldman Poisson structure on the set of geodesic functions implies the Fenchel--Nielsen bracket is more difficult  technically because it requires finding proper normalizations of twist coordinates $\tau_{B_i}$ to ensure the vanishing Poisson bracket between them (note that brackets between $\ell_{A_i}$ vanish automatically because geodesics constituting a pair-of-pant decomposition do not intersect, and their geodesic functions therefore commute). We demonstrate in Theorem~\ref{th:twist} that the canonical normalization of twist coordinates ensures their commutativity. The corresponding formulas coincide (modulo fixing some typos) with those obtained by Nekrasov, Rosly, and Shatashvili in \cite{NRS} for the basic constructing blocks of the pair-of-pants decomposition: $\Sigma_{0,4}$ and $\Sigma_{1,1}$.

Another powerful approach to the description of Teichm\"uller spaces ${\mathfrak T}_{g,s}$ is due to ideal-triangle decomposition of $\Sigma_{g,s}$ and related Thurston's shear coordinates \cite{ThSh} and Penner's lambda lengths \cite{Penn1} for surfaces with punctures generalized by V.V.Fock \cite{Fock1} to surfaces with holes. A log-canonical mapping-class group invariant Poisson structure was introduced by Fock on the set of shear coordinates in \cite{F97}.

An explicit combinatorial construction of the corresponding classical geodesic  functions in terms of shear coordinates of decorated Teichm\"uller spaces for Riemann surfaces with $s>0$ holes was proposed in \cite{ChF2}: it was shown there that geodesic  functions of all closed geodesics are Laurent polynomials of exponentiated shear coordinates with positive integer coefficients;  these results were extended to orbifold Riemann surfaces in  \cite{ChSh} where generalized cluster transformations (cluster algebras with coefficients) were introduced. Note that shear coordinates can be identified with the $Y$-type cluster variables \cite{FZ},~\cite{FZ2}, and mapping-class morphisms can be identified with cluster mutations.

The main result of  \cite{ChF2} is that the constant Poisson brackets on the set of shear coordinates induce Goldman brackets on the set of geodesic functions. Combined with the result of Theorem~\ref{th:twist} this immediately implies that {\em the Fock Poisson structure on the set of shear coordinates for $\Sigma_{g,s}$ with $s>0$ is the Fenchel--Nielsen Poisson structure on the set of length--twist coordinates}.

The next step was to generalize results obtained for Riemann surfaces with holes to Riemann surfaces $\Sigma_{g,s,n}$ with holes and with $n$ marked points on the hole boundaries (geometrically, these marked points are bordered cusps decorated with horocycles). A quantitative description of surfaces with marked points on  boundary components was given from different perspectives by Fock and Goncharov \cite{FG1}, Musiker, Schiffler and Williams \cite{MSW1}, \cite{MSW2}, \cite{MW}, and S. Fomin, M. Shapiro, and D. Thurston \cite{FST}, \cite{FT}.

In \cite{ChMaz2}, the quantum bordered Riemann surfaces $\Sigma_{g,s,n}$ and the corresponding quantum Teichm\"uller spaces $\mathcal T^\hbar_{g,s,n}$ were constructed. Note that having at least one bordered cusp on one of boundary components enables constructing an ideal-triangle decomposition of such a surface $\Sigma_{g,s,n}$ in which all arcs start and terminate at bordered cusps; holes without cusps and orbifold points then have to be confined in monogons (loops), and one has to restrict mapping-class morphisms to generalized cluster mutations preserving this condition. For such ideal-triangle decompositions of $\Sigma_{g,s,n}$ we have a bijection between the extended shear coordinates and lambda lengths, which enables one to determine Poisson and quantum relations for lambda lengths; they appear to be \cite{ChMaz2} correspondingly Poisson and quantum cluster algebras by Berenstein and Zelevinsky \cite{BerZel}. On the other hand, a mapping-class group invariant 2-form \cite{Penn92} on the set of lambda-lengths generates the invariant 2-form on the set of extended shear coordinates \cite{BK}, \cite{Ch20}.

It was proved in  \cite{ChMaz2} that the Poisson and quantum structure on the set of extended shear coordinates, which is also clearly mapping-class group invariant, induces the extended Goldman bracket on the set of geodesic functions and $\lambda$-lengths of {\em arcs}---geodesics starting and terminating at bordered cusps. Hence, the second objective of this paper carried out in Sec.~\ref{s:arcs} is to construct a canonical extension of the Fenchel--Nielsen coordinates to surfaces with bordered cusps on the subset of decoration-independent variables.

For the integrity of presentation, we leave aside two very interesting topics: the first is quantization of Poisson structures under consideration. Recall that algebras of observables are representations of quantum geodesic  functions constructed out of exponentiated quantum shear coordinates \cite{ChF1}; the pivotal observation was  that quantum flips enjoy quantum pentagon identity \cite{ChF1}, \cite{Kashaev} based  on the quantum dilogarithm function \cite{Faddeev}. The second, novel topic, is a generalization of Fenchel--Nielsen coordinates to the case of monodromies of $SL(n,\mathbb R)$ Fuchsian systems using Fock--Goncharov higher Teichm\"uller spaces \cite{FG1}. The corresponding classical and quantum monodromies were constructed in \cite{ShSh} and \cite{ChSh2} for $\Sigma_{g,s,n}$ with $n>0$, where it was shown that arc elements of these monodromies satisfy the Goldman brackets \cite{Gold} for $SL(n,\mathbb R)$ and are subject to Fock--Rosly algebras \cite{FockRosly} developed for $SL(n,\mathbb R)$-monodromies in \cite{CMR3}. Approaches to corresponding Fenchel--Nielsen coordinates include the construction of {\sl spectral networks} (see, e.g.,~\cite{HK}) and higher Labourie--McShane identities \cite{LM09} used by Huang and Sun \cite{HuangSun} for constructing special potentials in higher Teichm\"uller theory. However, a complete construction of generalizations of Fenchel--Nielsen coordinates for higher-rank algebras is still lacking and deserves further studies.

\section{Classical and Poisson algebras of geodesic functions and  $\lambda$-lengths (cluster variables) and combinatorial models of Teichm\"uller spaces}\label{s:preliminaries}
\setcounter{equation}{0}

In this section, we first recall the Goldman Poisson bracket on the set of geodesic functions and $\lambda$-lengths of arcs and make a brief excursion into the combinatorial description of Teichm\"uller space ${\mathfrak T}_{g,s,n}$ of Riemann surfaces of genus $g$ with $s>0$ holes/orbifold points, and with $n>0$ decorated {\it bordered cusps} situated on the hole boundaries.

\subsection{Goldman brackets and skein relations}\label{ss:Goldman}
In Fig.~\ref{fi:skein}, we present three basic relations valid for all geodesic functions $G_\gamma$ corresponding to closed geodesics $\gamma$ and for $\lambda$-lengths $\lambda_{\mathfrak a}$ of arcs $\mathfrak a$ on any Riemann surface $\Sigma_{g,s,n}$ (we allow only arcs starting and terminating at decorated bordered cusps). We can replace any (or both) $G_{\gamma_i}$ in these relations by $\lambda_{\mathfrak a_i}$ with corresponding natural adjustments of the right-hand sides.

\begin{figure}[tb]
\begin{pspicture}(-5,-1.5)(5,1.5){\psset{unit=1}
\rput(-3,0){
\psclip{\pscircle[linewidth=1.5pt, linestyle=dashed](0,0){1}}
\rput(0,0){\psline[linewidth=1.5pt,linecolor=red, linestyle=dashed]{<-}(0.7,-0.7)(-0.7,0.7)}
\rput(0,0){\psline[linewidth=3pt,linecolor=white](-1,-1)(1,1)}
\rput(0,0){\psline[linewidth=1.5pt,linecolor=blue,  linestyle=dashed]{->}(-0.7,-0.7)(0.7,0.7)}
\endpsclip
\rput(0,-0.2){\makebox(0,0)[ct]{$p$}}
\rput(-0.8,-.8){\makebox(0,0)[rt]{$G_{\gamma_1}$}}
\rput(0.8,-.8){\makebox(0,0)[lt]{$G_{\gamma_2}$}}
}
\rput(0,0){
\psclip{\pscircle[linewidth=1.5pt, linestyle=dashed](0,0){1}}
\rput(0,-1.4){\psarc[linewidth=1.5pt,linecolor=green, linestyle=dashed](0,0){1}{45}{135}}
\rput(0,1.4){\psarc[linewidth=1.5pt,linecolor=green, linestyle=dashed](0,0){1}{225}{315}}
\endpsclip
\rput(0,-1.2){\makebox(0,0)[ct]{$G_{\gamma_1\circ_P\gamma_2}$}}
\rput(-1.2,0){\makebox(0,0)[rc]{$1\cdot$}}
}
\rput(3,0){
\psclip{\pscircle[linewidth=1.5pt, linestyle=dashed](0,0){1}}
\rput(-1.4,0){\psarc[linewidth=1.5pt,linecolor=green, linestyle=dashed](0,0){1}{-45}{45}}
\rput(1.4,0){\psarc[linewidth=1.5pt,linecolor=green, linestyle=dashed](0,0){1}{135}{225}}
\endpsclip
\rput(0,-1.2){\makebox(0,0)[ct]{$G_{\gamma_1\circ_P\gamma_2^{-1}}$}}
\rput(-1.2,0){\makebox(0,0)[rc]{$1\cdot$}}
}
\rput(-1.7,0){
\rput(0,0){\makebox(0,0){$=$}}}
\rput(1.3,0){
\rput(0,0){\makebox(0,0){$+$}}}
}
\end{pspicture}
%$$
%\tcb{Goldman brackets for geodesic functions}
%$$
\begin{pspicture}(-4,-1.5)(4,1.5){\psset{unit=1}
\rput(-3,0){
\psclip{\pscircle[linewidth=1.5pt, linestyle=dashed](0,0){1}}
\rput(0,0){\psline[linewidth=1.5pt,linecolor=red, linestyle=dashed]{<-}(0.7,-0.7)(-0.7,0.7)}
\rput(0,0){\psline[linewidth=3pt,linecolor=white](-1,-1)(1,1)}
\rput(0,0){\psline[linewidth=1.5pt,linecolor=blue,  linestyle=dashed]{->}(-0.7,-0.7)(0.7,0.7)}
\endpsclip
\rput(0,-1.2){\makebox(0,0)[ct]{$\{G_{\gamma_1},G_{\gamma_2}\}_P$}}
\rput(0,-0.2){\makebox(0,0)[ct]{$p$}}
\rput(0.8,0.8){\makebox(0,0)[lb]{$\gamma_1$}}
\rput(0.8,-0.8){\makebox(0,0)[lt]{$\gamma_2$}}
}
\rput(0,0){
\psclip{\pscircle[linewidth=1.5pt, linestyle=dashed](0,0){1}}
\rput(0,-1.4){\psarc[linewidth=1.5pt,linecolor=green, linestyle=dashed](0,0){1}{45}{135}}
\rput(0,1.4){\psarc[linewidth=1.5pt,linecolor=green, linestyle=dashed](0,0){1}{225}{315}}
\endpsclip
\rput(0,-1.2){\makebox(0,0)[ct]{$G_{\gamma_1\circ_P\gamma_2}$}}
\rput(-1.2,0){\makebox(0,0)[rc]{$\dfrac12$}}
}
\rput(3,0){
\psclip{\pscircle[linewidth=1.5pt, linestyle=dashed](0,0){1}}
\rput(-1.4,0){\psarc[linewidth=1.5pt,linecolor=green, linestyle=dashed](0,0){1}{-45}{45}}
\rput(1.4,0){\psarc[linewidth=1.5pt,linecolor=green, linestyle=dashed](0,0){1}{135}{225}}
\endpsclip
\rput(0,-1.2){\makebox(0,0)[ct]{$G_{\gamma_1\circ_P\gamma_2^{-1}}$}}
\rput(-1.2,0){\makebox(0,0)[rc]{$\dfrac12$}}
}
\rput(-1.7,0){
\rput(0,0){\makebox(0,0){$=$}}}
\rput(1.3,0){
\rput(0,0){\makebox(0,0){$-$}}}
}
\end{pspicture}
%$$
%\tcb{Quantum skein relation}
%$$
\begin{pspicture}(-5,-1.5)(5,1.5){\psset{unit=1}
\rput(-4,0){
\psclip{\pscircle[linewidth=1.5pt, linestyle=dashed](0,0){1}}
\rput(0,0){\psline[linewidth=1.5pt,linecolor=red, linestyle=dashed]{<-}(0.7,-0.7)(-0.7,0.7)}
\rput(0,0){\psline[linewidth=3pt,linecolor=white](-1,-1)(1,1)}
\rput(0,0){\psline[linewidth=1.5pt,linecolor=blue,  linestyle=dashed]{->}(-0.7,-0.7)(0.7,0.7)}
\endpsclip
\rput(0,-0.2){\makebox(0,0)[ct]{$p$}}
\rput(0,-1.2){\makebox(0,0)[ct]{$G^\hbar_{\gamma_1}G^\hbar_{\gamma_2}$}}
\rput(0.8,0.8){\makebox(0,0)[lb]{$\gamma_1$}}
\rput(0.8,-0.8){\makebox(0,0)[lt]{$\gamma_2$}}
}
\rput(0.2,0){
\psclip{\pscircle[linewidth=1.5pt, linestyle=dashed](0,0){1}}
\rput(0,-1.4){\psarc[linewidth=1.5pt,linecolor=green, linestyle=dashed](0,0){1}{45}{135}}
\rput(0,1.4){\psarc[linewidth=1.5pt,linecolor=green, linestyle=dashed](0,0){1}{225}{315}}
\endpsclip
\rput(0,-1.2){\makebox(0,0)[ct]{$G^\hbar_{\gamma_1\circ_P\gamma_2}$}}
\rput(-1.2,0){\makebox(0,0)[rc]{$q^{1/2}$}}
}
\rput(3.8,0){
\psclip{\pscircle[linewidth=1.5pt, linestyle=dashed](0,0){1}}
\rput(-1.4,0){\psarc[linewidth=1.5pt,linecolor=green, linestyle=dashed](0,0){1}{-45}{45}}
\rput(1.4,0){\psarc[linewidth=1.5pt,linecolor=green, linestyle=dashed](0,0){1}{135}{225}}
\endpsclip
\rput(0,-1.2){\makebox(0,0)[ct]{$G^\hbar_{\gamma_1\circ_P\gamma_2^{-1}}$}}
\rput(-1.2,0){\makebox(0,0)[rc]{$q^{-1/2}$}}
}
\rput(-2.3,0){
\rput(0,0){\makebox(0,0){$=\otimes_P\Biggl[$}}}
\rput(1.3,0){
\rput(0.2,0){\makebox(0,0){$+$}}}
\rput(5,0){
\rput(0,0){\makebox(0,0)[lc]{$\Biggr],\quad q=e^{-i\pi\hbar}$}}}
}
\end{pspicture}
\caption{\small The classical skein relation, the Goldman bracket, and the quantum skein relation.}
\label{fi:skein}
\end{figure}

We have the following pattern: geodesic functions $G_{\gamma_1\circ_P\gamma_2}$ and $G_{\gamma_1\circ_P\gamma_2^{-1}}$ in the right-hand sides correspond to closed geodesics obtained by the corresponding resolutions of the intersection at the point $P$; these closed geodesics become arcs if we replace exactly one of geodesics $\gamma_i$ by an arc $\mathfrak a_i$, and in this case the both geodesic functions in the right-hand sides are to be replaced by $\lambda$-lengths of the corresponding arcs. If we replace both closed geodesics by arcs in the left-hand sides, then, instead of a single closed geodesics, in the right-hand side terms we have products of two lambda lengths of arcs obtained by the corresponding resolutions of the crossing at the point $P$. Note also that both geodesic functions and lambda lengths are insensitive to the choice of direction of the closed geodesic/arc. The rules for multiple crossings are as follows:

\begin{itemize}
\item
The classical skein relation holds at {\em any} intersection point $P$; in particular, when we replace both $G_{\gamma_1}$ and $G_{\gamma_2}$ by $\lambda_{\mathfrak a_1}$ and $\lambda_{\mathfrak a_2}$, in the right-hand side we obtain the celebrated Ptolemy relation of \cite{Penn1}. We can apply the skein relation recurrently; every time we obtain an empty loop (a contractible closed curve), we set its ``geodesic function'' to be $G_\varnothing:=-2$ and every time we have an empty arc contractible to a cusp, we set $\lambda_\varnothing:=0$.
\item
The Poisson bracket between two geodesic functions/lambda length is given by the sum over all intersection points of their local resolutions:
\be\label{Gold1}
\{G_{\gamma_1},G_{\gamma_2}\} =\sum_{P\in \gamma_1\# \gamma_2} \{G_{\gamma_1},G_{\gamma_2}\}_P,
\ee
with $\{G_{\gamma_1},G_{\gamma_2}\}_P$ depicted in the figure. We can then apply classical skein relations to terms in the right-hand side at will.
For $\lambda$-lengths of two arcs starting or terminating at the same bordered cusp we have to supplement (\ref{Goldman}) by a homogeneous Poisson relation
$$
\{\lambda_{\mathfrak a_1}, \lambda_{\mathfrak a_2}\}=\frac14 \lambda_{\mathfrak a_1}\lambda_{\mathfrak a_2},
$$
where the arc $\mathfrak a_1$ is to the left from the arc $\mathfrak a_2$ when looking from the cusp; we then add these relations to (\ref{Gold1}) considering cusps shared by arcs as additional intersection points.
\item
For the completeness, we also present the quantum skein relation between two quantum geodesic functions/quantum lambda lengths (the latter can be identified with quantum cluster variables). When we have several intersections, we have to apply the quantum skein relation {\em simultaneously} at all intersection points thus obtaining in the right-hand side a linear combination of quantum laminations---sets of non(self)-intersecting geodesic functions and arcs. Recall that quantum geodesic functions/quantum lambda lengths of nonintersecting geodesics/arcs commute. The quantum geodesic function $G^\hbar_\varnothing$ corresponding to a contractible loop is set to be $G^\hbar_\varnothing:=-q-q^{-1}$ and the quantum lambda length of an arc contractible to a cusp $\lambda^\hbar_\varnothing:=0$; for more details, see \cite{ChMaz2}.
\end{itemize}

\subsection{Fat graph description for Riemann surfaces $\Sigma_{g,s,n}$ and Teichm\"uller spaces $\mathfrak T_{g,s,n}$}\label{ss:fat}
We briefly recall the combinatorial description based on ideal-triangle decompositions of Riemann surfaces $\Sigma_{g,s,n}$, the corresponding (extended) shear coordinates, and related geodesic functions and lambda lengths.

\begin{definition}\label{def-pend}
A fat graph (a connected labelled graph with the fixed cyclic ordering of edges incident to each vertex) ${\mathcal G}_{g,s,n}$ is called a {\em spine} of $\Sigma_{g,s,n}$ with $s_h>0$ holes each containing $n_i>0$ bordered cusps ($\sum_{i} n_i=n>0$), and $s_o\ge 0$ holes/orbifold points without bordered cusps $(s=s_h+s_o)$ if this graph can be embedded  without self-intersections in $\Sigma_{g,s,n}$, its vertices are three-valent except exactly $n$ one-valent vertices---endpoints of $n$ {\it pending edges} with $n_i$ pending edges oriented towards the interior of the corresponding boundary component (a hole) and corresponding to $n_i$ marked points on the hole boundary. Furthermore, we require all holes without marked points and all orbifold points to be confined in monogons (loops) of the spine. We label by $\alpha,j,i$ all $6g-6+3s+2n$ edges of the graph; a real numbers $Z_\alpha$ (a shear coordinate)  corresponds to the $\alpha$th edge that is neither a pending edge, nor a loop. Every pending edge carries a real number $\pi_j$ (an extended shear coordinate), and loops carry numbers (coefficients) $\omega_i= 2\cosh(P_i/2)\ge 2$ for loops circumnavigating holes with perimeters $P_i\ge 0$ and $\omega_i= 2\cos(\pi/p_i)$ for loops with orbifold points of orders  $p_i\ge 2$ inside these loops.
\end{definition}

We identify $Z_\alpha$ in Definition~\ref{def-pend} with {\em (Thurston) shear coordinates} (see \cite{ThSh},\cite{Bon2}) and $\pi_j$ with extended shear coordinates \cite{ChSh}; the coordinate set $\{Z_\alpha, \pi_j,w_i\}$ parameterizes the decorated Teichm\"uller space ${\mathfrak T}_{g,s,n}$ and it was proved in \cite{ChSh} that these sets parameterize all metrizable Riemann surfaces modulo a discretely acting groupoid of flip morphisms and vice versa, every such set corresponds to a metrizable Riemann surface.

Fat graphs ${\mathcal G}_{g,s,n}$ are in bijection with \emph{ideal triangle decompositions} of $\Sigma_{g,s,n}$  \cite{Penn1}: monogons containing either holes without marked points or orbifold points are considered elements of this decomposition, and remaining ideal triangles correspond to three-valent vertices of the graph ${\mathcal G}_{g,s,n}$. Edges of this decomposition are {\em arcs}: geodesic lines starting and terminating at bordered cusps, these arcs are in bijection with edges of ${\mathcal G}_{g,s,n}$ that are not loops. We set into correspondence to an arc $\mathfrak a$ its $\lambda$-length, $\lambda_{\mathfrak a} =e^{\ell_{\mathfrak a}/2}$---the exponential of a half of the signed length of a part of $\mathfrak a$ stretched between the horocycles decorating the end cusps of the arc (the sign is negative if these horocycles intersect). We then have nondegenerate monoidal relations between $\lambda$-lengths of these arcs and exponentiated extended shear coordinates $\{Z_\alpha, \pi_j\}$ (see \cite{All}, \cite{ChMaz2}). Note that coefficients $\omega_i$ do not contribute to these relations.

\subsection{The Fuchsian group $\Delta_{g,s}$, geodesic functions, and $\lambda$-lengths}\label{ss:geodesic}

A metrizable surface $\Sigma_{g,s,n}$ is the quotient of the Poincar\'e hyperbolic upper half plane under a discrete action of a Fuchsian group $\Delta_{g,s}\subset PSL(2,\mathbb R)$. Note that this group is ``almost'' insensitive to the presence of bordered cusps: adding these cusps result only in that the action of elements corresponding to boundaries of holes at which these cusps are added becomes nontrivial. The standard fact in the hyperbolic geometry is that we have sets of 1-1 correspondences
$$
{\psset{unit=1}
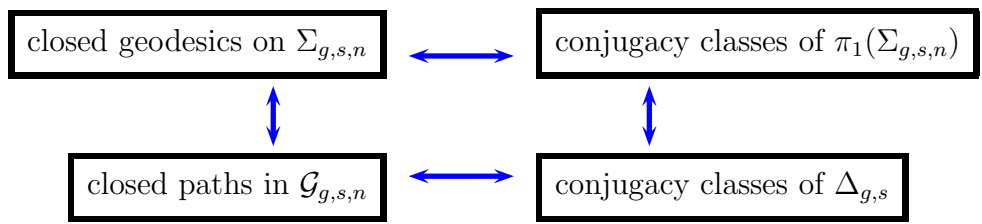
\begin{pspicture}(-4,-1.4)(4,1.4)
\setlength{\fboxsep}{5pt}
\setlength{\fboxrule}{2pt}
\rput(-1,0.5){\makebox(0,0)[rb]{\framebox{{closed geodesics} on $\Sigma_{g,s,n}$}}}
\rput(-1,-0.5){\makebox(0,0)[rt]{\framebox{{closed paths} in ${\mathcal G}_{g,s,n}$}}}
\rput(1,0.5){\makebox(0,0)[lb]{\framebox{conjugacy classes of $\pi_1(\Sigma_{g,s,n})$}}}
\rput(1,-0.5){\makebox(0,0)[lt]{\framebox{conjugacy classes of $\Delta_{g,s}$}}}
\rput(0,0){
\psline[linewidth=2pt,linecolor=blue]{<->}(-0.7,0.8)(0.7,0.8)
\psline[linewidth=2pt,linecolor=blue]{<->}(-0.7,-0.8)(0.7,-0.8)
\psline[linewidth=2pt,linecolor=blue]{<->}(-2.5,-0.4)(-2.5,0.4)
\psline[linewidth=2pt,linecolor=blue]{<->}(2.5,-0.4)(2.5,0.4)
}
\end{pspicture} }
$$
The principal advantage of the fat-graph description of Teichm\"uller spaces $\mathfrak T_{g,s,n}$ are very simple and explicit formulas expressing main algebraic objects, geodesic functions and lambda lengths, in terms of shear coordinates and a very simple Poisson bracket on the set of $Z_\alpha$, $\pi_j$. We begin with describing groupoid of paths.

The {\em groupoid of paths} is the set of homotopy classes of directed paths starting and terminating at the bordered cusps: we denote $\mathfrak a_{i\to j}$ a path starting at the $i$th cusp and terminating at the $j$th cusp. We endow this set with the natural partial composition law:  $\mathfrak a_{j\to k}\circ \mathfrak a_{i\to j}=\mathfrak a_{i\to k}$. Note that for any path $\mathfrak a_{i\to j}$ on $\Sigma_{g,s,n}$ we have a unique path without backtrackings in ${\mathcal G}_{g,s,n}$ lying in the same homotopy class; we therefore use the same notation for paths in $\Sigma_{g,s,n}$ and paths in  ${\mathcal G}_{g,s,n}$; to each such path we set into a correspondence an element from $PSL(2,\mathbb R)$, which we construct as a product of elementary $2\times 2$ matrices (all matrix products go from right to left). Note that elements of the path groupoid are those of {\em decorated character variety} $(SL(2,\mathbb R))^{2g+s+n-2}/U^n$ (see \cite{CMR2}).

Every time the path in a graph ${\mathcal G}_{g,s,n}$ goes along $\alpha$th inner edge or starts or terminates at a pending edge, we insert~\cite{Fock1} the so-called {\it edge matrix}  $X_{Z_\alpha}$ or $ X_{\pi_j}$,
\be
\label{XZ} X_{Y}=\left(
\begin{array}{cc} 0 & -e^{Y/2}\\
                e^{-Y/2} & 0\end{array}\right),\quad Y=Z_\alpha\ \hbox{or}\ Y=\pi_j.
\ee
When a path makes right or left turn at a three-valent vertex, we insert the corresponding ``right'' and ``left'' turn matrices
\be
\label{R}
R=\left(\begin{array}{cc} 1 & 1\\ -1 & 0\end{array}\right), \qquad
L= R^2=\left(\begin{array}{cc} 0 & 1\\ -1 &
-1\end{array}\right),
\ee
and finally, when a path is going along an $i$th loop clockwise (counterclockwise), we insert the matrix $F_{\omega_i}$ (or $-F^{-1}_{\omega_i}$):
\be
\label{F-p}
F_{\omega_i}:=
\left(\begin{array}{cc} 0 & 1\\ -1 & -w_i\end{array}\right),\qquad -F^{-1}_{\omega_i}:=
\left(\begin{array}{cc} w_i & 1\\ -1 & 0\end{array}\right)
\ee
sandwiched between two edge matrices $X_{Z_\alpha}$ of a unique edge incident to the loop.

An element $P_{\mathfrak a}\in PSL(2,\mathbb R)$ in the groupoid of paths has then the typical structure:
\be
\label{Pgamma}
P_{\mathfrak a_{j_1\to j_2}}=X_{\pi_{j_2}}LX_{Z_n}RX_{Z_{n-1}}\cdots RX_{Z_{k+1}}LX_{Z_k}F_{\omega_i}X_{Z_k} L X_{Z_{k-1}}R\dots X_{Z_1}L X_{\pi_{j_1}}
\ee
for a path starting at cusp $j_1$ and terminating at cusp $j_2$. The $\lambda$-length of this path is then the {\em upper-right matrix element} of $P_{\mathfrak a}$ (and it is easy to see that it does not depend on the path direction). For paths $\mathfrak a_{j\to j}$ starting and terminating at the same cusp, we obtain the {\em geodesic functions} as traces of the corresponding path matrices,
\be
\label{G}
G_{\gamma}\equiv \tr P_{\mathfrak a_{j\to j}}=2\cosh(\ell_\gamma/2),
\ee
where $\ell_\gamma$ is the length of the closed geodesic that is homeomorphic to the arc $\mathfrak a_{j\to j}$ upon identifying its endpoints and erasing the thus obtained marked point. Note that the backtracking part that appears in this procedure in the matrix product is then cancelled under the trace sign, so we can always consider only closed paths without backtrackings when evaluating geodesic functions.

This construction implies the following fundamental property: For any graph $\mathcal G_{g,s,n}$ all matrix elements of $P_{\mathfrak a}$ for every arc ${\mathfrak a}$ (and, correspondingly, all $\lambda$-lengths and geodesic functions $G_\gamma$) are polynomials with sign-definite integer coefficients of exponentiated shear coordinates and coefficients: $\lambda_{\mathfrak a}\in \mathbb Z_+[[e^{\pi_j/2}, e^{\pm Z_\alpha/2}, \omega_i]]$ and $G_{\gamma}\in \mathbb Z_+[[e^{\pm Z_\alpha/2}, \omega_i]]$.

%If we have at least one bordered cusp, we can invert formulas (\ref{cross-l}) and express $\lambda$-lengths in terms of extended shear coordinates. The polynomial expressions obtained from (\ref{Pgamma}) for $\lambda$-lengths of arcs become monoidal for arcs constituting the ideal triangle decomposition dual to a given fat graph  ${\mathcal G}_{g,s,n}$.

%%% BEGIN REMARK
\rem{
%%%%%%%

\begin{lemma}\label{lem:lambda}\cite{ChMaz2}
Paths $P_{\mathfrak a_\alpha}$ corresponding to arcs $\mathfrak a_\alpha$ from the ideal triangle decomposition dual to a given fat graph  ${\mathcal G}_{g,s,n}$ and labelled by edges $\alpha,j$ (that are not loops, and the first and the last edge are pending edges) of this graph have the product structure  $P_{\mathfrak a_\alpha}=X_{\pi_2} R X_{Z_{\beta_n}}\cdot R\cdots R F_w R\cdots R X_{Z_\alpha} L X_{Z_{\beta_j}} \cdots X_{Z_{\beta_1}} L X_{\pi_1}$, and their upper-right elements are merely
$$
e^{l_{\mathfrak a_\alpha}/2}=e^{(\pi_1+\pi_2+Z_\alpha +\sum Z_\beta)/2},
$$
where the sum ranges (with multiplicities) all edges along the path $\mathfrak a_\alpha$. Note that coefficients $w_i$ do not contribute to these $\lambda$-lengths.
\end{lemma}
%%%% END REMARK
}
%%%%%

\subsection{Poisson and symplectic structures}\label{ss:Poisson}

One of the most attractive properties of the fat graph description is a very simple Poisson algebra on the set of coordinates $Z_\alpha$, $\pi_j$ \cite{F97}: let $Y_k$, $k=1,2,3 \mod 3$, denote either $Z$- or $\pi$-variables of cyclically ordered edges incident to a three-valent vertex. The Poisson (Weil--Petersson) bi-vector field  is then
\be
\label{WP-PB}
w_{\text{WP}}:=\sum_{{\hbox{\small 3-valent} \atop \hbox{\small vertices} }}
\,\sum_{k=1}^{3} \partial_{Y_k}\wedge \partial_{Y_{k+1}},
\ee

\begin{theorem}\label{th-WP}\cite{ChF2} \cite{ChMaz2}
The bracket \eqref{WP-PB} induces the {\em Goldman
bracket}  \cite{Gold} (see  Fig.~\ref{fi:skein}) on the set of geodesic functions and on the set of $\lambda$-lengths of arcs.
\end{theorem}

The center of the Poisson algebra {\rm(\ref{WP-PB})} is generated by sums of $Z_\alpha$ and $\pi_j$ (taken with multiplicities) of all edges incident to a given hole, so together with the coefficients $w_i$ we have exactly $s$ independent central elements.

A mapping-class group invariant  {\em symplectic structure} on the set of $\lambda$-lengths of arcs was introduced by Penner \cite{Penn92}:
\be
\label{WP-SS}
\Omega_{\text{WP}}:=\sum_{{\hbox{\small ideal} \atop \hbox{\small triangles} }}
\,\sum_{k=1}^{3}d\log \lambda_k\wedge d\log \lambda_{k+1}.
\ee
In the case of $\Sigma_{g,s,n}$ with $n>0$, the bracket (\ref{WP-PB}) induces homogeneous Poisson relations on the set of $\lambda$-lengths of arcs from the same ideal-triangle decomposition amounting to Berenstein--Zelevinsky quantum cluster algebras \cite{BerZel}, and this Poisson structure was shown  \cite{Ch20} to be inverse to the symplectic structure (\ref{WP-SS}).

%Despite the simplicity of the above formulas, before introducing bordered cusps all attempts of constructing Poisson brackets for $\lambda$ lengths or symplectic structure for shear coordinates were marred with discrepancies. But now, armed with explicit formulas from Lemma~\ref{lem:lambda}, we can solve this problem. First, in \cite{ChMaz2} the author together with Marta Mazzocco found Poisson and quantum structures on the set of $\lambda$-lengths; the result reproduced Berenstein--Zelevinsky quantum cluster algebras \cite{BerZel}. The last missing item was a symplectic form on the set of extended shear coordinates; we complete it in Theorem~\ref{th:Omega} in the next section.

\subsection{Flip morphisms of fat graphs}\label{ss:flip}

Different fat-graph parameterizations of Teichm\"uller spaces are related by mapping-class morphisms generated by sequences of {\em flip morphisms} (mutations) of edges: any two spines from a given topological class are related by a finite sequence of flips, and every time this sequence results in a graph homotopically identical to the original graph (note that we have infinitely many copies of the moduli space $\mathcal M_{g,s,n}$ in the Teichm\"uller space $\mathfrak T_{g,s,n}$ but only a finite number of homotopically different spines $\mathcal G_{g,s,n}$), transformations of the variables $\{Z_\alpha,\pi_j\}$ describe a transformation from the mapping class group. It is also a classic result that all such transformations are produced by Dehn twists along closed geodesics, which we use in the next section.

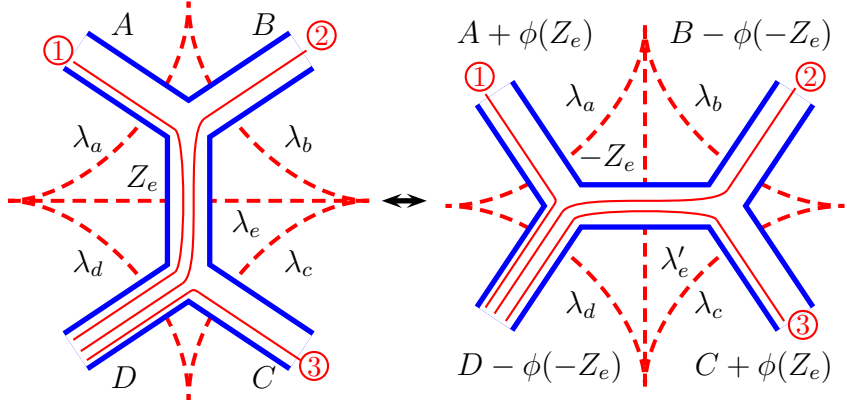
\begin{figure}[tb]
\begin{pspicture}(-3,-3)(4,3){
\newcommand{\FLIP}{%
{\psset{unit=1}
\psline[linewidth=18pt,linecolor=blue](0,-1)(0,1)
\psline[linewidth=18pt,linecolor=blue](0,1)(1.5,2)
\psline[linewidth=18pt,linecolor=blue](0,1)(-1.5,2)
\psline[linewidth=18pt,linecolor=blue](0,-1)(1.5,-2)
\psline[linewidth=18pt,linecolor=blue](0,-1)(-1.5,-2)
\psline[linewidth=14pt,linecolor=white](0,-1)(0,1)
\psline[linewidth=14pt,linecolor=white](0,1)(1.5,2)
\psline[linewidth=14pt,linecolor=white](0,1)(-1.5,2)
\psline[linewidth=14pt,linecolor=white](0,-1)(1.5,-2)
\psline[linewidth=14pt,linecolor=white](0,-1)(-1.5,-2)
}
}
\rput(-2.5,0){\psbezier[linecolor=red, linestyle=dashed, linewidth=1.5pt](2.4,0)(1.4,0)(0,-1.3)(0,-2.65)
\psbezier[linecolor=red, linestyle=dashed, linewidth=1.5pt](-2.4,0)(-1.4,0)(0,-1.3)(0,-2.65)
\psbezier[linecolor=red, linestyle=dashed, linewidth=1.5pt](2.4,0)(1.4,0)(0,1.3)(0,2.65)
\psbezier[linecolor=red, linestyle=dashed, linewidth=1.5pt](-2.4,0)(-1.4,0)(0,1.3)(0,2.65)
\psline[linecolor=red, linestyle=dashed, linewidth=1.5pt](2.4,0)(-2.4,0)
\FLIP}
\rput(.3,0){\psline[linewidth=2pt]{<->}(-0.3,0)(0.3,0)}
\rput{90}(3.5,0){\psbezier[linecolor=red, linestyle=dashed, linewidth=1.5pt](2.4,0)(1.4,0)(0,-1.3)(0,-2.65)
\psbezier[linecolor=red, linestyle=dashed, linewidth=1.5pt](-2.4,0)(-1.4,0)(0,-1.3)(0,-2.65)
\psbezier[linecolor=red, linestyle=dashed, linewidth=1.5pt](2.4,0)(1.4,0)(0,1.3)(0,2.65)
\psbezier[linecolor=red, linestyle=dashed, linewidth=1.5pt](-2.4,0)(-1.4,0)(0,1.3)(0,2.65)
\psline[linecolor=red, linestyle=dashed, linewidth=1.5pt](2.4,0)(-2.4,0)
\FLIP}
\rput(-2.5,0){
\rput(-1.1,2.2){\makebox(0,0)[lb]{$A$}}
\rput(1.1,2.2){\makebox(0,0)[rb]{$B$}}
\rput(-0.9,0.3){\makebox(0,0)[lc]{$Z_e$}}
\rput(1.1,-2.2){\makebox(0,0)[rt]{$C$}}
\rput(-1.1,-2.2){\makebox(0,0)[lt]{$D$}}
\rput(0.9,-0.3){\makebox(0,0)[rc]{$\lambda_e$}}
\rput(-1.2,0.65){\makebox(0,0)[rb]{$\lambda_a$}}
\rput(1.2,0.65){\makebox(0,0)[lb]{$\lambda_b$}}
\rput(1.2,-0.65){\makebox(0,0)[lt]{$\lambda_c$}}
\rput(-1.2,-0.65){\makebox(0,0)[rt]{$\lambda_d$}}
\psline[linecolor=red](-1.5,-2)(-.15,-1.1)
\psline[linecolor=red](1.5,2)(.15,1.1)
\psbezier[linecolor=red](-.15,-1.1)(0.15,-0.9)(-0.15,0.9)(.15,1.1)
\rput(-.1,.15){\psline[linecolor=red](-1.5,-2)(-.15,-1.1)}
\rput(-.1,-.15){\psline[linecolor=red](-1.5,2)(-.15,1.1)}
\psbezier[linecolor=red](-.25,-.95)(-0.1,-0.85)(-0.1,0.85)(-.25,.95)
\rput(.08,-.12){\psline[linecolor=red](-1.5,-2)(-.15,-1.1)}
\rput(-.08,-.12){\psline[linecolor=red](1.5,-2)(.15,-1.1)}
\psbezier[linecolor=red](-.07,-1.22)(0,-1.176)(0,-1.176)(.07,-1.22)
}
\rput(3.5,0){
\rput(-2.5,-2){\makebox(0,0)[lt]{$D - \phi(-Z_e)$}}
\rput(2.5,-2){\makebox(0,0)[rt]{$C+\phi(Z_e)$}}
\rput(2.5,2){\makebox(0,0)[rb]{$B-\phi(-Z_e)$}}
\rput(-2.5,2){\makebox(0,0)[lb]{$A+\phi(Z_e)$}}
\rput(-0.5,0.4){\makebox(0,0)[cb]{$-Z_e$}}
\rput(0.4,-0.6){\makebox(0,0)[ct]{$\lambda'_{e}$}}
\rput(-0.65,1.2){\makebox(0,0)[rb]{$\lambda_a$}}
\rput(0.65,1.2){\makebox(0,0)[lb]{$\lambda_b$}}
\rput(0.65,-1.2){\makebox(0,0)[lt]{$\lambda_c$}}
\rput(-0.65,-1.2){\makebox(0,0)[rt]{$\lambda_d$}}
}
\rput{90}(3.5,0){
\psline[linecolor=red](1.5,-2)(.15,-1.1)
\psline[linecolor=red](-1.5,2)(-.15,1.1)
\psbezier[linecolor=red](.15,-1.1)(-0.15,-0.9)(0.15,0.9)(-.15,1.1)
\rput(-.1,.15){\psline[linecolor=red](-1.5,-2)(-.15,-1.1)}
\rput(-.1,-.15){\psline[linecolor=red](-1.5,2)(-.15,1.1)}
\psbezier[linecolor=red](-.25,-.95)(-0.1,-0.85)(-0.1,0.85)(-.25,.95)
}
\rput{270}(3.5,0){
\rput(.08,-.12){\psline[linecolor=red](-1.5,-2)(-.15,-1.1)}
\rput(-.08,-.12){\psline[linecolor=red](1.5,-2)(.15,-1.1)}
\psbezier[linecolor=red](-.07,-1.22)(0,-1.176)(0,-1.176)(.07,-1.22)
}
% added geodesic lines
\rput(-2.5,0){
\put(-1.8,1.9){\makebox(0,0)[cc]{\hbox{\tcr{\small$1$}}}}
\put(-1.8,1.9){\pscircle[linecolor=red]{.2}}
\put(1.7,2.1){\makebox(0,0)[cc]{\hbox{\tcr{\small$2$}}}}
\put(1.7,2.1){\pscircle[linecolor=red]{.2}}
\put(1.6,-2.3){\makebox(0,0)[cc]{\hbox{\tcr{\small$3$}}}}
\put(1.6,-2.3){\pscircle[linecolor=red]{.2}}
}
%right half
\rput(3.5,0){
\put(-2.2,1.7){\makebox(0,0)[cc]{\hbox{\tcr{\small$1$}}}}
\put(-2.2,1.7){\pscircle[linecolor=red]{.2}}
\put(2.2,1.7){\makebox(0,0)[cc]{\hbox{\tcr{\small$2$}}}}
\put(2.2,1.7){\pscircle[linecolor=red]{.2}}
\put(2.1,-1.6){\makebox(0,0)[cc]{\hbox{\tcr{\small$3$}}}}
\put(2.1,-1.6){\pscircle[linecolor=red]{.2}}
}
}
\end{pspicture}
\caption{\small Flip on an inner edge (labeled ``$e$'') that is neither a loop nor adjacent to a loop. We indicate the correspondences between paths in the graph undergoing the flip. Dashed lines are arcs of the dual ideal triangle decomposition.}
\label{fi:flip}
\end{figure}

We distinguish between two types of flip morphisms: those induced by flips of ``typical'' inner edges (see Fig. \ref{fi:flip}) and those induced by flips of edges that are adjacent to a loop (see Fig. \ref{fi:interchange-p-dual}); no flips can be performed on pending edges and loops.

%\subsubsection{Flipping inner edges}

\begin{lemma} \label{lem-abc}~\cite{ChF1,ChF2}
In the notation of Fig.~\ref{fi:flip}, the transformation
$$
({\tilde A},{\tilde B},{\tilde C},{\tilde D},{\tilde Z_e})=(A+\phi(Z_e), B-\phi(-Z_e), C+\phi(Z_e), D-\phi(-Z_e), -Z_e),
%\label{abc}
$$
where $\phi (Z)={\rm log}(1+e^Z)$, preserves path products {\rm(\ref{G})} (thus preserving both geodesic functions and $\lambda$-lengths)  simultaneously preserving Poisson structure {\rm(\ref{WP-PB})} on the shear coordinates. The dual Ptolemy transformation of $\lambda$-lengths (cluster mutation), $\lambda_e\lambda'_{e}=\lambda_a\lambda_c+\lambda_b\lambda_d$ preserves the symplectic structure (\ref{WP-SS}).
\end{lemma}

Because the proof of the lemma is local w.r.t. the graph $\mathcal G_{g,s,n}$ and follows from matrix equalities $X_DRX_{Z_e}RX_{A}=X_{\tilde A}RX_{\tilde D}$, $X_DRX_{Z_e}LX_B=X_{\tilde D}LX_{\tilde Z_e}RX_{\tilde B}$, and $X_CLX_D=X_{\tilde C}LX_{\tilde Z_e}LX_{\tilde D}$, each pertaining to the corresponding path pattern in Fig.~\ref{fi:flip}, it can be extended to the whole groupoid of $SL(2,\mathbb R)$ monodromies. The same statement is therefore valid for the $\lambda$-lengths of the corresponding arcs. We have a similar statement for flips of inner edges incident to loops:

%\subsubsection{Flipping the edge incident to a loop}\label{sss:pending}

\begin{figure}[tb]
\begin{pspicture}(-3,-2.5)(4,2.5){
\newcommand{\FLIP}{%
{\psset{unit=1}
\psbezier[linewidth=1.5pt, linestyle=dashed, linecolor=red](-5,0)(-3,2.5)(0,2.5)(2,0)
\psbezier[linewidth=1.5pt, linestyle=dashed, linecolor=red](-5,0)(-3,-2.5)(0,-2.5)(2,0)
\psbezier[linewidth=1.5pt, linestyle=dashed, linecolor=red](-5,0)(-3.5,1.6)(-1,1.6)(-1,0)
\psbezier[linewidth=1.5pt, linestyle=dashed, linecolor=red](-5,0)(-3.5,-1.6)(-1,-1.6)(-1,0)
\psline[linewidth=18pt,linecolor=blue](-2,0)(0,0)
\psline[linewidth=18pt,linecolor=blue](0,0)(1,1.5)
\psline[linewidth=18pt,linecolor=blue](0,0)(1,-1.5)
\pscircle[linewidth=2pt,linecolor=blue,fillstyle=solid,fillcolor=white](-2.8,0){1}
\psline[linewidth=14pt,linecolor=white](-2,0)(0,0)
\psline[linewidth=14pt,linecolor=white](0,0)(1,1.5)
\psline[linewidth=14pt,linecolor=white](0,0)(1,-1.5)
\rput(-2.8,0){\pscircle[linecolor=blue,linewidth=2pt,fillstyle=solid,fillcolor=white](0,0){0.5}}
%
%\psline[linewidth=1.5pt, linestyle=dashed, linecolor=red](-5,0)(-3.8,0)
%\psline[linewidth=1.5pt, linestyle=dashed, linecolor=red](-2.8,0)(-3.3,0)
%\rput(-2,0){\pscircle*{0.05}}
\psarc[linecolor=red](-2,0.3){0.1}{200}{270}
\psarc[linecolor=red](-2,0.3){0.2}{200}{270}
\psarc[linecolor=red](-2,-0.3){0.1}{90}{160}
\psarc[linecolor=red](-2,-0.3){0.2}{90}{160}
\psline[linecolor=red](-2,0.1)(-0.01,0.1)
\psline[linecolor=red](-2,0.2)(-0,0.2)
\psline[linecolor=red](-2,-0.1)(-0.01,-0.1)
\psline[linecolor=red](-2,-0.2)(-0,-0.2)
\psarc[linecolor=red](-2.8,0){0.75}{20}{340}
\psarc[linecolor=red](-2.8,0){0.65}{20}{340}
}
}
\rput(-2.5,0){\FLIP}
\psline[linewidth=2pt]{<->}(-0.3,0)(0.3,0)
\rput{180}(2.5,0){\FLIP}
\rput(-2.5,0){
\rput(0.2,1.5){\makebox(0,0)[lb]{$A$}}
\rput(-0.8,0.5){\makebox(0,0)[lb]{$Z_e$}}
\rput(0.2,-1.5){\makebox(0,0)[lt]{$B$}}
\rput(1.3,0.8){\makebox(0,0)[lb]{$\lambda_a$}}
\rput(1.3,-0.8){\makebox(0,0)[lt]{$\lambda_b$}}
\rput(-1,-0.5){\makebox(0,0)[lt]{$\lambda_e$}}
\rput(-1.2,2){\makebox(0,0)[rb]{$\displaystyle w=e^\xi+e^{-\xi},\ \xi\in\mathbb C$}}
\rput(-2.8,0){\makebox(0,0)[cc]{$\omega$}}
}
\rput(3.5,0){
\rput(-2.3,0.8){\makebox(0,0)[rb]{$\lambda_a$}}
\rput(-2.3,-0.8){\makebox(0,0)[rt]{$\lambda_b$}}
\rput(0,-0.5){\makebox(0,0)[rt]{$\lambda'_e$}}
\rput(-1.4,-2){\makebox(0,0)[lt]{$\displaystyle B - \phi(-Z_e+\xi)-\phi(-Z_e-\xi)$}}
\rput(-1.4,2){\makebox(0,0)[lb]{$\displaystyle A + \phi(Z_e+\xi)+\phi(Z_e-\xi)$}}
\rput(-0.1,0.5){\makebox(0,0)[rb]{$-Z_e$}}
\rput(1.8,0){\makebox(0,0)[cc]{$\omega$}}
}
\rput{90}(3.5,0){
%\rput(-.1,.15){\psline[linecolor=red](-1.5,-2)(-.15,-1.1)}
\rput(-.1,-.15){\psline[linecolor=red](-1.5,2)(-.09,1.06)}
\rput(.1,-.15){\psline[linecolor=red](1.5,2)(.09,1.06)}
%\psbezier[linecolor=red](-.25,-.95)(-0.1,-0.85)(-0.1,0.85)(-.25,.95)
}
\rput{270}(-3.5,0){
%\rput(-.1,.15){\psline[linecolor=red](-1.5,-2)(-.15,-1.1)}
\rput(-.1,-.15){\psline[linecolor=red](-1.5,2)(-.09,1.06)}
\rput(.1,-.15){\psline[linecolor=red](1.5,2)(.09,1.06)}
}
\rput{270}(3.5,0){
\rput(.08,-.12){\psline[linecolor=red](-1.5,-2)(-.15,-1.1)}
\rput(-.08,-.12){\psline[linecolor=red](1.5,-2)(.15,-1.1)}
\psbezier[linecolor=red](-.07,-1.22)(0,-1.176)(0,-1.176)(.07,-1.22)
}
\rput{90}(-3.5,0){
\rput(.08,-.12){\psline[linecolor=red](-1.5,-2)(-.15,-1.1)}
\rput(-.08,-.12){\psline[linecolor=red](1.5,-2)(.15,-1.1)}
\psbezier[linecolor=red](-.07,-1.22)(0,-1.176)(0,-1.176)(.07,-1.22)
}
\rput(-2.5,0){
\rput(-.06,.04){\psline[linecolor=red](1,1.5)(0.1,.15)}
\rput(.06,-.04){\psline[linecolor=red](1,1.5)(0.1,.15)}
\psbezier[linecolor=red](.04,.19)(-0.02,0.1)(-0.02,0.1)(-0.1,0.1)
\psbezier[linecolor=red](.16,.11)(0.02,-0.1)(0.02,-.1)(-0.1,-.1)
}
\rput(2.5,0){
\rput(.06,.04){\psline[linecolor=red](-1,1.5)(-0.1,.15)}
\rput(-.06,-.04){\psline[linecolor=red](-1,1.5)(-0.1,.15)}
\psbezier[linecolor=red](-.04,.19)(0.02,0.1)(0.02,0.1)(0.1,0.1)
\psbezier[linecolor=red](-.16,.11)(-0.02,-0.1)(-0.02,-.1)(0.1,-.1)
}
% added geodesic lines
\rput(-3.5,0){
\put(2.3,-1.6){\makebox(0,0)[cc]{\hbox{\tcr{\small$3$}}}}
\put(2.3,-1.6){\pscircle[linecolor=red]{.2}}
\put(2.2,1.7){\makebox(0,0)[cc]{\hbox{\tcr{\small$2$}}}}
\put(2.2,1.7){\pscircle[linecolor=red]{.2}}
\put(1.95,-1.8){\makebox(0,0)[cc]{\hbox{\tcr{\small$1$}}}}
\put(1.95,-1.8){\pscircle[linecolor=red]{.2}}
}
%right half
\rput(3.5,0){
\put(-2.3,-1.6){\makebox(0,0)[cc]{\hbox{\tcr{\small$1$}}}}
\put(-2.3,-1.6){\pscircle[linecolor=red]{.2}}
\put(-2.2,1.7){\makebox(0,0)[cc]{\hbox{\tcr{\small$2$}}}}
\put(-2.2,1.7){\pscircle[linecolor=red]{.2}}
\put(-1.95,-1.8){\makebox(0,0)[cc]{\hbox{\tcr{\small$3$}}}}
\put(-1.95,-1.8){\pscircle[linecolor=red]{.2}}
}
}
\end{pspicture}
\caption{\small
The transformation of shear coordinates under the flip of an edge incident to a loop; we indicate how paths change upon flipping the edge. Dashed lines are edges of the dual ideal triangle decomposition.
}
\label{fi:interchange-p-dual}
\end{figure}
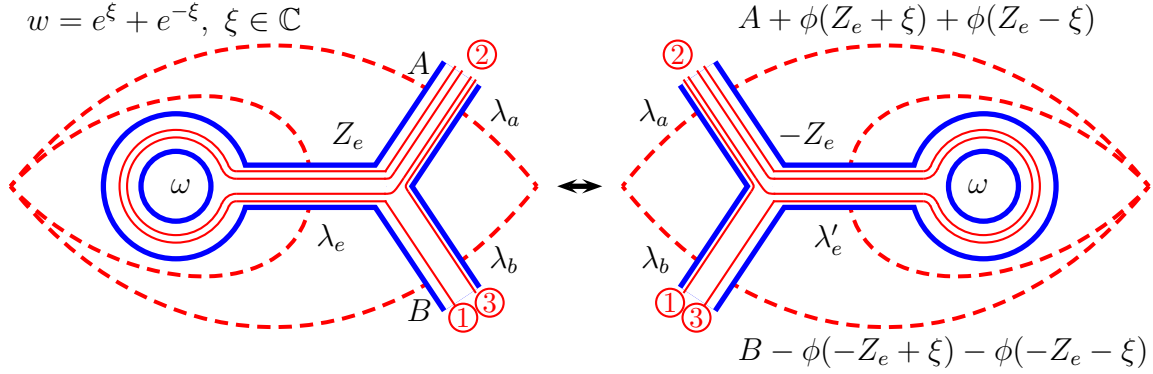

\begin{lemma} \label{lem-pending1}{\rm (\cite{ChSh},\cite{ChMaz2})}
The transformation~in Fig.~\ref{fi:interchange-p-dual}
$$
\{\tilde A,\tilde B,\tilde Z_e\}:= \{A+\phi(Z_e+\xi)+\phi(Z_e-\xi), B-\phi(-Z_e+\xi)-\phi(-Z_e-\xi),-Z_e\},\quad w=e^\xi+e^{-\xi},
$$
where $\phi(x)=\log(1+e^x)$ and $\xi\in\mathbb C$ is a morphism of the space ${\mathfrak T}_{g,s,n}$ that  preserves both Poisson structures {\rm(\ref{WP-PB})} and the path elements from $SL(2,\mathbb R)$. The dual transformation (generalized cluster transformation) $\lambda_e\lambda'_{e}=\lambda_a^2+w\lambda_a\lambda_b+\lambda_b^2$ preserves the symplectic structure  (\ref{WP-SS}).
\end{lemma}

\section{Fenchel--Nielsen brackets and shear coordinates}\label{s:FN-PB}
\setcounter{equation}{0}

\subsection{Fenchel--Nielsen coordinates for Riemann surfaces with holes}\label{ss:FN}

\subsubsection{Fenchel--Nielsen coordinates for $\Sigma_{1,1}$}\label{ss:FN-11}

We first consider the case of a torus with one hole of perimeter $p$. Choosing a closed geodesic $A$ with the length $\ell_A$ and a dual closed geodesic $B$ having a single intersection point with $A$ we are about to construct a twist coordinate $\tau_B$ that is a function of $G_B$, $\ell_A$ and $G_P:=e^{p/2}+e^{-p/2}$. It is useful to consider together with $B$ all geodesics obtained by Dehn rotations along the geodesics $A$: $A^nB$, $n\in \mathbb Z$. The Poisson bracket is
\be
\{G_A,G_{A^nB}\}=\frac12 (G_{A^{n-1}B}-G_{A^{n+1}B}),\quad n\in \mathbb Z,
\label{AB}
\ee
where $G_{A^{n-1}B}$ and $G_{A^{n+1}B}$ are two solutions of the quadratic equations generated by a {\em Markov triple},
\be
G_AG_{A^nB}G_{A^{n-1}B}-G_A^2-G_{A^nB}^2-G_{A^{n-1}B}^2=G_P-2, \quad n\in\mathbb Z.
\label{ABC}
\ee
We also have the classical skein relations
\be
G_AG_{A^nB}=G_{A^{n+1}B}+G_{A^{n-1}B},\quad n\in \mathbb Z.
\label{sk-torus}
\ee

\begin{proposition}\label{pr:torus}
A twist coordinate $\tau_B$ having a constant unit bracket with the length $\ell_A$ of a closed geodesic $A$, $\{\tau_B,\ell_A\}=1$, is
\be
\tau_{A^nB}=\log\bigl(G_{A^{n-1}B}-e^{-\ell_A/2}G_{A^nB} \bigr).
\label{twist-torus}
\ee
Then $\{\ell_A,\tau_{A^nB}\}=1$ for all $n\in\mathbb Z$ and all $\tau_{A^nB}$ are related by constant shifts:
\be
\tau_{A^nB}-\tau_{A^{n+1}B}=\ell_A/2.
\label{mod-torus}
\ee
In particular, a Dehn twist along $A$ transforms $B$ into $A^2B$, and the corresponding modular transformation is $\tau_{B}\to \tau_{A^2B}=\tau_B-\ell_A$, so, as expected, the twist coordinate (\ref{twist-torus}) takes values between $0$ and $\ell_A$ in a single copy of the modular space labelled by $A$ in Mirzakhani's terminology \cite{Mir06}.
\end{proposition}

The {\em proof} is a direct calculation. We have to choose the sign of $\ell_A$: applying (\ref{sk-torus}) we obtain, say, for $n=-1$,
\be
G_{A^{-2}B}-e^{-\ell_A/2}G_{A^{-1}B}=(e^{\ell_A/2}+e^{-\ell_A/2})G_{A^{-1}B} - e^{-\ell_A/2}G_{A^{-1}B}=e^{\ell_A/2}(G_{A^{-1}B}-e^{-\ell_A/2}G_{B})
\label{sk-2}
\ee
Checking that $\{\tau_B,\ell_A\}=1$ is a straightforward calculation:
\begin{align*}
&\{e^{\tau_B},G_A\}=\frac12 e^{\tau_B}(e^{\ell_A/2}-e^{-\ell_A/2})\{\tau_B,\ell_A\}
=\frac12\Bigl[ G_{A^{-2}B}-G_B- e^{-\ell_A/2}\bigl( G_{A^{-1}B}-G_{AB}\bigr) \Bigr]\\
&=\frac12\Bigl[ G_AG_{A^{-1}B}-2G_B- e^{-\ell_A/2}\bigl( 2G_{A^{-1}B}-G_{A}G_{B}\bigr) \Bigr]
=\frac12\bigl[e^{\ell_A/2}-e^{-\ell_A/2} \Bigr] \bigl( G_{A^{-1}B}-e^{-\ell_A/2}G_{B}\bigr).
\end{align*}

What if we choose $\tau_B$ in another form, say, $\tau'_B=\log\bigl(G_{A^{-1}B}-e^{\ell_A/2}G_{B} \bigr)$? It is easy to see, using the Markov triple and skein relations, that
\be
\bigl(G_{A^{-1}B}-e^{-\ell_A/2}G_{B} \bigr)\bigl(G_{A^{-1}B}-e^{\ell_A/2}G_{B} \bigr)=e^{\ell_A}+e^{-\ell_A}+e^{p/2}+e^{-p/2},
\label{minus-tau}
\ee
so, since $\tau_B$ is defined up to a constant shift by $f(\ell_A,p)$, we can identify $\tau'_B$ with negative $\tau_B$.

We still have however a substantial ambiguity in choosing $\tau_B$: we can add to it any function $f(\ell_A,p)$ (normalization) without breaking both the Poisson brackets with $\ell_A$ and shift symmetries. Note however that adding this function (i) changes the symmetry properties under the inversion $\ell_A\to -\ell_A$ and (ii) affects Poisson relations between twist coordinates $\tau_B$ and $\tau_{B'}$ dual to different cycles $\gamma_A$ and $\gamma_{A'}$ from the same pair-of-pant decomposition of $\Sigma_{g,s}$. Below we see that choosing the normalization such that $\tau_B\to -\tau_B$ under the changing of the sign of $\ell_A$ allows us to solve these two problems simultaneously. Moreover, the thus constructed canonical twist coordinates have a clear geometrical sense.

\begin{figure}[tb]
\begin{pspicture}(-3,-2.5)(4,4){
\newcommand{\PATTERN}{%
{\psset{unit=1}
\rput{90}(0,-1){\psellipse[linecolor=blue, linewidth=1pt](0,0)(1,.3)\psframe[linecolor=white, fillstyle=solid, fillcolor=white](-1.1,0)(1.1,.6)\psellipse[linecolor=blue, linestyle=dashed, linewidth=1pt](0,0)(1,.3)}
\rput(0,0){\psarc[linecolor=blue](0,1){3}{210}{330}}
\rput(0,0){\psarc[linecolor=blue](0,1){1}{210}{330}}
\rput(0,0){\psarc[linecolor=blue](0,-0.74){1}{60}{120}}
\rput(0,0){\psbezier[linecolor=blue](-2.6,-0.5)(-3.6,1.236)(-1,2)(-1,3)}
\rput(0,0){\psbezier[linecolor=blue](2.6,-0.5)(3.6,1.236)(1,2)(1,3)}
\psellipse[linecolor=blue, linewidth=1pt](0,3)(1,.5)
}
}
\rput(-2.5,0){\PATTERN}
\rput{-5}(-2.7,0.2){\psellipse[linecolor=red, linewidth=1pt](0,0)(2.08,1.05)}
%\rput{75}(-2.1,-0.05){\parametricplot[linecolor=red,linestyle=solid,linewidth=1pt]{0}{180}{2.05 t cos mul  0.45 t sin  mul}
%\parametricplot[linecolor=red,linestyle=dashed,linewidth=1pt]{180}{360}{2.05 t cos mul  0.45 t sin  mul}}
%
\rput{0}(-2.6,0){\parametricplot[linecolor=green,linestyle=solid,linewidth=1pt]{-81.5}{276.5}{2.08 t cos mul 1 t sin mul 0.002 t mul add}}
%
%
%
%\psline[linewidth=2pt]{<->}(-0.5,0)(0.5,0)
%\rput{180}(2.5,0){\FLIP}
%
\rput(3,-5){{\psset{unit=5}
\psarc[linewidth=1pt,linecolor=magenta](-2,0){2.5}{23}{37}
\psarc[linewidth=1pt,linecolor=green](-3,0){3.354}{16.5}{26.5}
\psarc[linewidth=1pt,linecolor=blue](0.5,0){1}{101.5}{120}
\psline[linewidth=1pt,linecolor=red](0,0.87)(0,1.5)
}}
\rput(3,0){
\put(0,5.5){\makebox(0,0)[bc]{\hbox{{$S$}}}}
\put(0,2){\makebox(0,0)[tc]{\hbox{{$P$}}}}
\put(1.1,2.45){\makebox(0,0)[tc]{\hbox{{$Q$}}}}
\put(1.07,2.53){\psarc[linewidth=0.5pt](0,0){0.3}{105}{200}}
\put(0.92,2.61){\makebox(0,0)[cc]{\hbox{{$\cdot$}}}}
\put(1.5,2.55){\makebox(0,0)[tc]{\hbox{{$R$}}}}
}
\rput(5.5,-1){
\rput(0,3){\makebox(0,0)[lc]{$\displaystyle |SP|=\ell_{B}/2$}}
\rput(0,2.2){\makebox(0,0)[lc]{$\displaystyle |PR|=\ell_{A}/2$}}
\rput(0,1.4){\makebox(0,0)[lc]{$\displaystyle |SR|=\ell_{A^{-1}B}/2$}}
\rput(0,0.6){\makebox(0,0)[lc]{$\displaystyle |PQ|=\widehat\tau_{B}$}}
}

% added geodesic lines
\rput(-3.5,-0.5){
\put(1.8,-1.5){\makebox(0,0)[cc]{\hbox{\tcg{\small$H$}}}}
%\put(2.4,-1.6){\pscircle[linecolor=green]{.2}}
\put(0.3,-1.2){\makebox(0,0)[cc]{\hbox{\tcr{\small$B$}}}}
%\put(2.3,1.7){\pscircle[linecolor=red]{.2}}
\put(1.3,-0.5){\makebox(0,0)[cl]{\hbox{\tcb{\small$A$}}}}
%\put(0,-0.5){\pscircle[linecolor=blue]{.2}}
}
%right half
}
\end{pspicture}
\caption{\small
A torus with a hole as a pair of pants
}
\label{fi:torus-hole}
\end{figure}
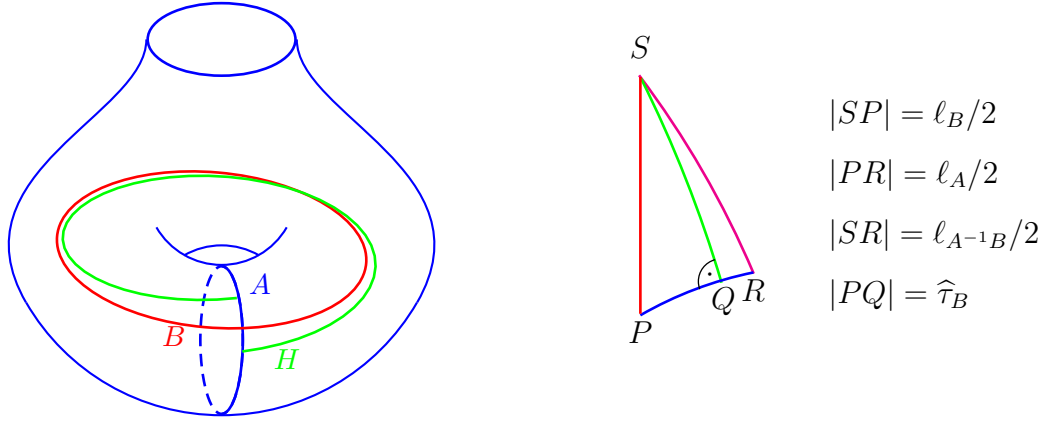

We may address the geometrical interpretation of the twist coordinate $\tau_B$. In the right side of Fig.~\ref{fi:torus-hole} we depict the triangle constituted by halves of geodesics $A$, $B$, and $A^{-1}B$. It is well-known fact from hyperbolic geometry (see the Appendix) that for a two curves $A$ and $B$ with the single intersection point $P$, the curve $A^{-1}B$ intersects $A$ and $B$ at the respective points $R$ and $S$ lying exactly at the respective distances  $\ell_A/2$ and $\ell_B/2$ from $P$; the distance $|SR|$ is then exactly half of the length of the geodesics $A^{-1}B$.

Using the standard formulas for a triangle in hyperbolic geometry (see, e.g., \cite{Buser}), we have
\[
\cosh(|PR|) \cosh(|SP|)-\cosh(|SR|)=\sinh(|SP|)\cosh(|PR|) \tanh(|PQ|),    %\o\o\o\o
\]
or, upon identification of the triangle sides as in Fig.~\ref{fi:torus-hole}, we obtain
\[
e^{2|PQ|}=\frac{\frac 12 G_AG_B -G_{A^{-1}B} -G_B\sinh (\ell_A/2)}{\frac 12 G_AG_B -G_{A^{-1}B} +G_B\sinh (\ell_A/2)}
=\frac{G_{A^{-1}B} -e^{-\ell_A/2}G_B}{G_{A^{-1}B} -e^{\ell_A/2}G_B},
\]
so multiplying the numerator and denominator by $G_A$ and applying again skein relation (\ref{sk-torus}), we obtain
\be
e^{2|PQ|}=\frac{e^{\ell_A/2}G_{A^{-1}B}-e^{-\ell_A/2}G_{AB} }{e^{-\ell_A/2}G_{A^{-1}B}-e^{\ell_A/2}G_{AB} }=\frac{e^{2\tau_B}}{e^{\ell_A}+e^{-\ell_A}+e^{p/2}+e^{-p/2}},
\label{PQ}
\ee
so $|PQ|=\tau_B+f(\ell_A,p)$, and the rates of change of the both quantities at fixed $\ell_A$ and $p$ are the same. We then declare the signed  $|PQ|$ to be a {\em canonical twist coordinate}.

\begin{lemma}\label{lm:torus}
The canonical twist coordinate that change its sign under changing the orientation ($\ell_A\to -\ell_A$) for a torus with one hole having the perimeter $p$ is
\be\label{tau-torus}
\widehat \tau_B:=\log \Bigl[ \frac{G_{A^{-1}B} - e^{-\ell_A/2}G_B}{(e^{\ell_A}+e^{-\ell_A}+e^{p/2}+e^{-p/2})^{1/2}} \Bigr]
\ee
Geometrically, $\widehat \tau_B$ is the half of the signed geodesic length $2|PQ|$ along the geodesic $\gamma_A$ between endpoints of a geodesic $H$ perpendicular to $\gamma_A$ and homeomorphic to the $B$-cycle (see Fig.~\ref{fi:torus-hole}).
\end{lemma}

\subsubsection{Fenchel--Nielsen coordinates for $\Sigma_{0,4}$}\label{ss:FN-04}

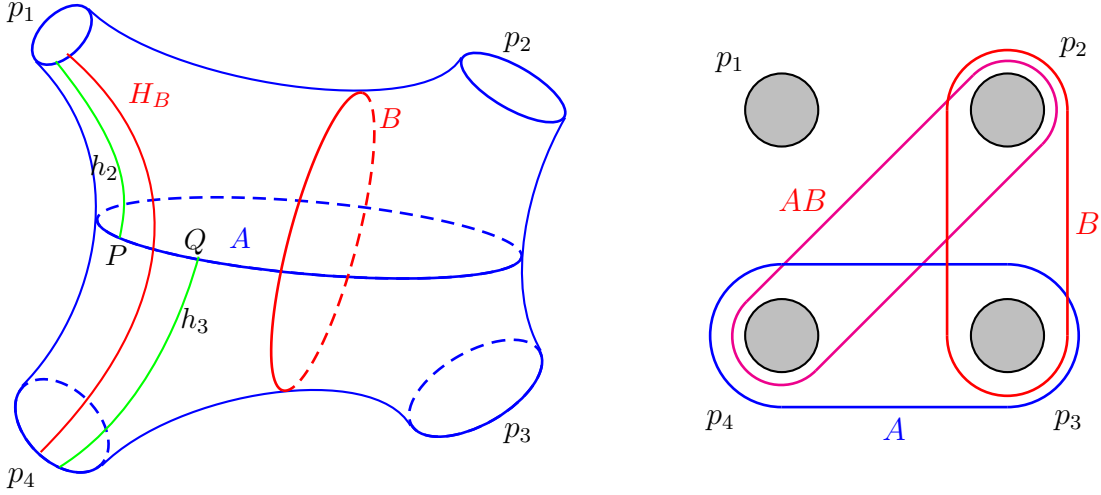
\begin{figure}[tb]
\begin{pspicture}(-3,-4)(4,4){
\newcommand{\PATTERN}{%
{\psset{unit=1}
\rput{-45}(-3,-2.5){\psellipse[linecolor=blue, linewidth=1pt](0,0)(.75,.5)\psframe[linecolor=white, fillstyle=solid, fillcolor=white](-1,0)(1,.6)\psellipse[linecolor=blue, linestyle=dashed,  linewidth=1pt](0,0)(.75,.5)}
\rput{30}(2.5,-2){\psellipse[linecolor=blue, linewidth=1pt](0,0)(1,.5)\psframe[linecolor=white, fillstyle=solid, fillcolor=white](-1.1,0)(1.1,.6)\psellipse[linecolor=blue, linestyle=dashed, linewidth=1pt](0,0)(1,.5)}
\rput{-5}(0.29,0){\psellipse[linecolor=blue,  linewidth=1pt](0,0)(2.85,0.5)\psframe[linecolor=white, fillstyle=solid, fillcolor=white](-2.9,0)(2.9,.6)\psellipse[linecolor=blue, linestyle=dashed,  linewidth=1pt](0,0)(2.85,.5)}
\rput{45}(-3,2.7){\psellipse[linecolor=blue,  linewidth=1pt](0,0)(.5,.3)}
\rput{-30}(3,2){\psellipse[linecolor=blue,  linewidth=1pt](0,0)(.8,.3)}
\psbezier[linecolor=blue](-3.53,-1.97)(-2.56,-1)(-2,1)(-3.35,2.35)
\psbezier[linecolor=blue](3.366,-1.5)(2.866,-0.634)(3.19,0.734)(3.69,1.6)
\psbezier[linecolor=blue](-2.47,-3.03)(-1.47,-2.03)(1.134,-1.634)(1.634,-2.5)
\psbezier[linecolor=blue](-2.65,3.05)(-1.65,2.05)(1.81,1.534)(2.31,2.4)
}
}
\rput(-2.5,0){\PATTERN}
\rput{75}(-2.1,-0.05){\parametricplot[linecolor=red,linestyle=solid,linewidth=1pt]{0}{180}{2.05 t cos mul  0.45 t sin  mul}
\parametricplot[linecolor=red,linestyle=dashed,linewidth=1pt]{180}{360}{2.05 t cos mul  0.45 t sin  mul}}
\rput(-2.5,0){\psbezier[linecolor=red](-3,2.45)(-1.3,1)(-1.5,-1)(-3.35,-2.85)
\psbezier[linecolor=green](-3.15,2.35)(-2.1,1.05)(-2.2,0.5)(-2.3,0)
\psbezier[linecolor=green](-3.1,-3.05)(-1.9,-2.3)(-1.4,-0.8)(-1.25,-0.25)}
\rput(-3.5,0.4){
\put(-0.9,2){\makebox(0,0)[cc]{\hbox{\tcr{\small$H_B$}}}}
\put(-1.7,1.2){\makebox(0,0)[tl]{\hbox{{\small$h_2$}}}}
\put(-0.3,-1){\makebox(0,0)[cc]{\hbox{{\small$h_3$}}}}
\put(-1.5,0){\makebox(0,0)[tl]{\hbox{{\small$P$}}}}
\put(-0.3,0.05){\makebox(0,0)[cc]{\hbox{{\small$Q$}}}}
%\put(2.4,-1.6){\pscircle[linecolor=green]{.2}}
\put(2.3,1.7){\makebox(0,0)[cc]{\hbox{\tcr{\small$B$}}}}
%\put(2.3,1.7){\pscircle[linecolor=red]{.2}}
\put(0.3,0.1){\makebox(0,0)[cc]{\hbox{\tcb{\small$A$}}}}
%\put(0,-0.5){\pscircle[linecolor=blue]{.2}}
\rput(-2.6,3.1){\makebox(0,0)[cc]{\hbox{{$p_1$}}}}
\rput(4,2.7){\makebox(0,0)[cc]{\hbox{{$p_2$}}}}
\rput(4,-2.5){\makebox(0,0)[cc]{\hbox{{$p_3$}}}}
\rput(-2.6,-3.1){\makebox(0,0)[cc]{\hbox{{$p_4$}}}}
}
%right half
\rput(5.5,0){
\put(-1.5,-1.5){\pscircle[fillstyle=solid,fillcolor=lightgray]{.5}}
\put(-1.5,1.5){\pscircle[fillstyle=solid,fillcolor=lightgray]{.5}}
\put(1.5,-1.5){\pscircle[fillstyle=solid,fillcolor=lightgray]{.5}}
\put(1.5,1.5){\pscircle[fillstyle=solid,fillcolor=lightgray]{.5}}
\rput{45}(0,0){\psarc[linewidth=1pt,linecolor=magenta](-2.13,0){.65}{90}{270}
\psarc[linewidth=1pt,linecolor=magenta](2.13,0){.65}{-90}{90}
\psline[linewidth=1pt,linecolor=magenta](-2.13,0.65)(2.13,0.65)
\psline[linewidth=1pt,linecolor=magenta](-2.13,-0.65)(2.13,-0.65)
}
\psarc[linewidth=1pt,linecolor=blue](1.5,-1.5){0.95}{-90}{90}
\psarc[linewidth=1pt,linecolor=blue](-1.5,-1.5){0.95}{90}{270}
\psline[linewidth=1pt,linecolor=blue](-1.5,-.55)(1.5,-.55)
\psline[linewidth=1pt,linecolor=blue](-1.5,-2.45)(1.5,-2.45)
\psarc[linewidth=1pt,linecolor=red](1.5,1.5){0.8}{0}{180}
\psarc[linewidth=1pt,linecolor=red](1.5,-1.5){0.8}{180}{360}
\psline[linewidth=1pt,linecolor=red](0.7,1.5)(0.7,-1.5)
\psline[linewidth=1pt,linecolor=red](2.3,1.5)(2.3,-1.5)
\put(-2,2){\makebox(0,0)[br]{\hbox{{$p_1$}}}}
\put(2.2,2.2){\makebox(0,0)[bl]{\hbox{{$p_2$}}}}
\put(2.5,-2.5){\makebox(0,0)[tr]{\hbox{{$p_3$}}}}
\put(-2.5,-2.5){\makebox(0,0)[tl]{\hbox{{$p_4$}}}}
\put(0,-2.6){\makebox(0,0)[tc]{\hbox{\tcb{$A$}}}}
\put(2.4,0){\makebox(0,0)[cl]{\hbox{\tcr{$B$}}}}
\put(-0.9,0.3){\makebox(0,0)[cr]{\hbox{\tcr{$AB$}}}}
}
}
\end{pspicture}
\caption{\small
A sphere with four holes: gluing two pair of pants. Here $H_B$, $h_2$, and $h_3$ are perpendiculars between respectively $\gamma_{p_1}$ and $\gamma_{p_4}$ (a unique perpendicular between these two cycles that has no intersections with $\gamma_B$), $\gamma_{p_1}$ and $\gamma_A$, and $\gamma_{p_4}$ and $\gamma_A$. The signed geodesic length (along the geodesic $\gamma_A$) $|PQ|$ is identified with the canonical twist coordinate $\widehat\tau_B$ (\ref{tau-can}).
}
\label{fi:four-hole}
\end{figure}

After a ``warming-up'' case of the torus with one hole let us proceed to a more laborous case of a four-holed sphere, which can be considered as two pairs of pants glued along a geodesic $A$ (see Fig.~\ref{fi:four-hole}). We let $G_i:=e^{p_i/2}+e^{-p_i/2}$, $i=1,\dots, 4$, denote the geodesic function for perimeters of four holes, the geodesics $A$ and $B$ have now two intersection points

The bracket is exactly the same as in the torus case (note the absence of the factor $1/2$):
\be
\{G_A,G_B\}=G_{A^{-1}B}-G_{AB},
\label{AB-M4}
\ee
where $G_{A^{-1}B}$ and $G_{AB}$ are two solutions of the quadratic equations generated by a Markov triple,
\be
G_AG_BG_{A^{\pm 1}B}-G_A^2-G_B^2-G_{A^{\pm 1}B}^2-G_{A^{\pm 1}B}S_{AB}-G_AS_A-G_{B}S_B=R,
\label{ABC-M4}
\ee
where
\be
S_A:=G_1G_2+G_3G_4,\quad S_B:=G_1G_4+G_2G_3,\quad S_{AB}:=G_1G_3+G_2G_4,
\label{S-4}
\ee
and
\be
R:=G_1G_2G_3G_4+\sum_{i=1}^4 G_i^2-4,
\label{S-R}
\ee
and the same relation (\ref{ABC-M4}) holds true upon the replacement of the indices $B\to A^{2k}B$ and $A^{\pm 1}B\to A^{2k\pm 1}$ for any integer $k$.

We can also use the classical skein relations
\be
G_AG_{A^{2k}B}=G_{A^{2k+1}B}+G_{A^{2k-1}B}+S_{AB},\quad G_AG_{A^{2k+1}B}=G_{A^{2k+2}B}+G_{A^{2k}B}+S_{B},\quad k\in \mathbb Z.
\label{sk-M4}
\ee

\begin{proposition}\label{lm-M4}
A twist coordinate having a constant unit bracket with the half-length $\ell_A/2$ of a closed geodesic $A$ for a sphere with four holes, $\{\tau_B,\ell_A/2\}=1$, is
\be
\tau_{A^{2k}B}=\log\Bigl(G_{A^{2k-1}B}-e^{-\ell_A/2}G_{A^{2k}B}-\frac{e^{-\ell_A/2} S_{AB}+S_B}{2\sinh(\ell_A/2)} \Bigr)
\label{twistB-M4}
\ee
and
\be
\tau_{A^{2k+1}B}=\log\Bigl(G_{A^{2k}B}-e^{-\ell_A/2}G_{A^{2k+1}B}-\frac{e^{-\ell_A/2} S_{B}+S_{AB}}{2\sinh(\ell_A/2)} \Bigr).
\label{twistAB-M4}
\ee
all these choices of a dual coordinate are related by constant shifts:
\be
\tau_{A^nB}-\tau_{A^{n+1}B}=\ell_A/2.  %%%% o\(
\label{mod-M4}
\ee
In particular, a Dehn twist along $A$ transforms $B$ into $A^2B$, and the corresponding modular transformation is again $\tau_{B}\to \tau_{A^2B}=\tau_B-\ell_A$, so, again, the twist coordinate (\ref{twistB-M4}) assumes values between $0$ and $\ell_A$ in a single copy of the modular space labelled by $A$.
\end{proposition}

The {\em proof} is a direct calculation. Checking that $\{\tau_B,\ell_A/2\}=1$ is straightforward using only (\ref{sk-M4}), which we leave to the reader.

Note that shifting a twist coordinate $\tau_B$ by any function $f(G_A,G_i)$ preserves Poisson brackets with all length coordinates.

Take, say, $\tau_{A^{2k}B}$ and use that $G_{A^{2k}B}=G_AG_{A^{2k-1}B}-G_{A^{2k-2}B}-S_B$:
\begin{align*}
\tau_{A^{2k}B}&=\log\Bigl(G_{A^{2k-1}B}-e^{-\ell_A/2}(G_AG_{A^{2k-1}B}-G_{A^{2k-2}B}-S_B)-\frac{e^{-\ell_A/2} S_{AB}+S_{B}}{2\sinh(\ell_A/2)} \Bigr)\\
&=\log\Bigl( e^{-\ell_A/2} (G_{A^{2k-2}B}-e^{-\ell_A/2} G_{A^{2k-1}B})-\frac{e^{-\ell_A/2} S_{AB}+e^{-\ell_A}S_{B}}{2\sinh(\ell_A/2)} \Bigr)=-\ell_A/{2}+\tau_{A^{2k-1}B}.
\end{align*}
We therefore again have that $\tau_{A^nB}-\tau_{A^{n+1}B}=\ell_A/2$ for any $n\in \mathbb Z$.

Let us see now what happens if we choose another sign of  $\ell_A$, i.e., we are about to compare $\tau_{A^{2k}B}$ (\ref{twistB-M4}) with
\[
\tau'_{A^{2k}B}:=\log\Bigl(G_{A^{2k-1}B}-e^{\ell_A/2}G_{A^{2k}B}+\frac{S_{B}+e^{\ell_A/2}S_{AB}}{2\sinh(\ell_A/2)} \Bigr).
\]
Then using (\ref{ABC-M4}), after some algebra, we obtain that
\be\label{tau-11}
\tau_{A^kB}+\tau'_{A^kB}=\log\Bigl( \frac{S_{AB}+G_AS_BS_{AB} +S_B^2}{4\sinh^2(\ell_A/2)}+G_A^2+S_AG_A+R\Bigr),
\ee
i.e., it does not depend on $G_B$, so $\tau'_{A^kB}$ and $\tau_{A^kB}$ have opposite rates of change upon varying $\ell_B$.

A more convenient way of writing the expression in (\ref{tau-11}) is due to the identity
\begin{align}
\label{tau-1}
&S_{AB}+G_AS_BS_{AB} +S_B^2+(G^2_A-4)(G_A^2+S_AG_A+R)\\
&\qquad=(G_AG_1G_2+G_A^2+G_1^2+G_2^2-4)(G_AG_3G_4+G_A^2+G_3^2+G_4^2-4).\nonumber
\end{align}
From this we can guess the proper normalization for the twist coordinate:

\begin{lemma}\label{lm:twist}
The canonical twist coordinate that change its sign under changing the orientation ($\ell_A\to -\ell_A$) is
\be\label{tau-can}
\widehat \tau_B:=\log \Bigl[ \frac{(e^{\ell_A/2}-e^{-\ell_A/2})(G_{A^{-1}B} - e^{-\ell_A/2}G_B) -e^{-\ell_A/2} S_{AB}-S_B}{(G_AG_1G_2+G_A^2+G_1^2+G_2^2-4)^{1/2}(G_AG_3G_4+G_A^2+G_3^2+G_4^2-4)^{1/2}} \Bigr]
\ee
Geometrically, $\widehat \tau_B$ is the signed geodesic length along the geodesic $\gamma_A$ between endpoints of perpendiculars to $\gamma_A$ from the holes $P_1$ and $P_4$.
\end{lemma}

Even before proving this lemma, let us formulate the main statement of this section.

\begin{theorem}\label{th:twist}
Given any pair-of-pant decomposition of a Riemann surface $\Sigma_{g,s}$ and defining the twist variable $\widehat \tau_B$ for each of $3g-3+s$ inner geodesics $\gamma_A$ of this decomposition by formulas (\ref{tau-can}) and (\ref{tau-torus}) we have that $\{ \widehat \tau_B, \widehat \tau_{B'} \}=0$, that is, all canonical twist coordinates Poisson commute.
\end{theorem}

\begin{remark}
The coordinates $\widehat \tau_B$ (\ref{tau-can}) and (\ref{tau-torus}) were introduced by Nekrasov, Rosly, and Shatashvili in \cite{NRS} as the canonical coordinates having unit brackets with the corresponding length coordinates $\ell_A$.
\end{remark}

\subsubsection{Proof of Lemma~\ref{lm:twist}}
Since the commutation relations between $\ell_A$ and $\widehat \tau_B$ were proven above, it only remains to identify the canonical twist with a geometrical object, namely, the (signed) geodesic distance $|PQ|$ in Fig.~\ref{fi:four-hole}. We evaluate it using hyperbolic geometry identities (a list of very useful identities can be found on p.454 of Buser's monograph \cite{Buser}). First, note that considering any pair of pants, that is, a sphere with three holes with, say, perimeters $p_1$, $p_2$, and $\ell_A$, cutting along perpendiculars between the hole boundaries turns it into a union of two mirror-symmetrical right-angled hexagons with cyclically ordered boundary lengths $\{ h_2, p_1/2, h_A, p_2/2,h_1,\ell_A/2\}$. We can then express $h_2$ via $p_1$, $p_2$, and $\ell_A$ by the formula
\be
\cosh(h_2)\sinh(p_1/2)\sinh(\ell_A/2)=\cosh(p_2/2)+\cosh(p_1/2)\cosh(\ell_A/2).
\ee
An analogous formula for another pair of pants bounded by holes with perimeters $p_3$, $p_4$, and $\ell_A$ expresses $h_3$:
\be
\cosh(h_3)\sinh(p_4/2)\sinh(\ell_A/2)=\cosh(p_3/2)+\cosh(p_4/2)\cosh(\ell_A/2).
\ee
The third formula expresses $H_B$: just note that it is a perpendicular in the three-holes sphere bounded by $p_1$, $p_4$ and $B$:
\be
\cosh(H_B)\sinh(p_1/2)\sinh(p_4/2)=\frac12 G_B+\cosh(p_1/2)\cosh(p_4/2).
\ee
Finally, we need the formula for $|PQ|$: note that it is the length of a side of a right-angled self-crossing hexagon (in Fig.~\ref{fi:four-hole} it is constituted by $h_2$, $h_3$, $H_B$, $PQ$ and parts of holes $p_1$ and $p_4$); the relevant relation from \cite{Buser} is
\be
\cosh(H_B)=\sinh(h_2)\sinh(h_3)\cosh(|PQ|)+\cosh(h_2)\cosh(h_3).
\ee
Using these formulas we are able to express $e^{|PQ|}$; it happens that some of square-root expressions appearing in the answer coincide with the expressions appearing when we solve equation (\ref{ABC-M4}) w.r.t. $G_{A^{-1}B}$; this enable us to simplify the expression and finally obtain the formula (\ref{tau-can}) in which $\widehat\tau_B=|PQ|$.

\subsubsection{Proof of Theorem~\ref{th:twist}}
\begin{figure}[tb]
\begin{pspicture}(-3,-3)(3,3){
\newcommand{\PATTERN}{%
{\psset{unit=1.5}
\rput(0,1.2){\pscircle[linecolor=black, linewidth=1pt, fillcolor=lightgray, fillstyle=solid](0,0){0.07}}
\rput(0,-1.2){\pscircle[linecolor=black, linewidth=1pt, fillcolor=lightgray, fillstyle=solid](0,0){0.07}}
\rput(0.8,1.2){\pscircle[linecolor=black, linewidth=1pt, fillcolor=lightgray, fillstyle=solid](0,0){0.07}}
\rput(0.8,-1.2){\pscircle[linecolor=black, linewidth=1pt, fillcolor=lightgray, fillstyle=solid](0,0){0.07}}
\rput(1.6,0){\pscircle[linecolor=black, linewidth=1pt, fillcolor=lightgray, fillstyle=solid](0,0){0.07}}
\rput(0.4,1.2){\psellipse[linecolor=red, linewidth=1pt](0,0)(.75,.3)
\psline[linecolor=red, linewidth=1pt]{<-}(-0.05,-0.29)(0.05,-0.29)}
\rput(0.4,-1.2){\psellipse[linecolor=red, linewidth=1pt](0,0)(.75,.3)
\psline[linecolor=red, linewidth=1pt]{<-}(-0.05,0.29)(0.05,0.29)}
\rput{-57}(1.2,0.6){\psellipse[linecolor=blue, linewidth=1pt](0,0)(0.9,.3)}
\rput{57}(1.2,-0.6){\psellipse[linecolor=blue, linewidth=1pt](0,0)(0.9,.3)}
\psarc[linewidth=1pt,linestyle=dashed](2.878,0){2.5}{150}{210}
\psarc[linewidth=1pt,linestyle=dashed](2.878,0){2.3}{150}{210}
\psarc[linewidth=1pt,linestyle=dashed](0.8,1.2){0.1}{-30}{150}
\psarc[linewidth=1pt,linestyle=dashed](0.8,-1.2){0.1}{-150}{30}
\psarc[linewidth=1pt,linestyle=dashed](0.32,0){1.34}{105}{255}
\psarc[linewidth=1pt,linestyle=dashed](0.32,0){1.14}{105}{255}
\psarc[linewidth=1pt,linestyle=dashed](0,1.2){0.1}{-75}{105}
\psarc[linewidth=1pt,linestyle=dashed](0,-1.2){0.1}{-105}{75}
\psarc[linewidth=1pt,linestyle=dashed](-0.4,0){1.56}{45}{315}
\psarc[linewidth=1pt,linestyle=dashed](-0.4,0){1.84}{45}{315}
\psarc[linewidth=1pt,linestyle=dashed](0.8,1.2){0.14}{-135}{45}
\psarc[linewidth=1pt,linestyle=dashed](0.8,-1.2){0.14}{-45}{135}
%
%\psarc[linewidth=1pt,linecolor=blue,linestyle=solid](0,0){0.5}{150}{390}
%
%
%\psbezier[linecolor=blue](-3.53,-1.97)(-2.56,-1)(-2,1)(-3.35,2.35)
%\psbezier[linecolor=blue](3.366,-1.5)(2.866,-0.634)(3.19,0.734)(3.69,1.6)
%\psbezier[linecolor=blue](-2.47,-3.03)(-1.47,-2.03)(1.134,-1.634)(1.634,-2.5)
%\psbezier[linecolor=blue](-2.65,3.05)(-1.65,2.05)(1.81,1.534)(2.31,2.4)
}
}
\rput(0,0){\PATTERN}
%\rput{75}(-2.1,-0.05){\parametricplot[linecolor=red,linestyle=solid,linewidth=1pt]{0}{180}{2.05 t cos mul  0.45 t sin  mul}
%\parametricplot[linecolor=red,linestyle=dashed,linewidth=1pt]{180}{360}{2.05 t cos mul  0.45 t sin  mul}}
%
\rput(0,0){\psset{unit=1.5}
%\put(2.4,-1.6){\pscircle[linecolor=green]{.2}}
%\put(2.3,1.7){\makebox(0,0)[cc]{\hbox{\tcr{\small$B$}}}}
%\put(2.3,1.7){\pscircle[linecolor=red]{.2}}
\put(0.05,0.9){\makebox(0,0)[tc]{\hbox{\tcr{\small$A'$}}}}
\put(1.6,0.6){\makebox(0,0)[cl]{\hbox{\tcb{\small${A'}^nB'$}}}}
\put(0.05,-0.9){\makebox(0,0)[bc]{\hbox{\tcr{\small$A$}}}}
\put(1.6,-0.6){\makebox(0,0)[cl]{\hbox{\tcb{\small${A}^m B$}}}}
%\put(0,-0.5){\pscircle[linecolor=blue]{.2}}
\put(-0.8,0){\makebox(0,0)[cl]{\hbox{{\small$D_{m{-}1,n{-}1}$}}}}
\put(-1.7,0){\makebox(0,0)[cl]{\hbox{{\small$D_{m,n}$}}}}
\put(-3.4,0){\makebox(0,0)[cl]{\hbox{{\small$D_{m{+}1,n{+}1}$}}}}
\put(0.2,1.2){\makebox(0,0)[cc]{\hbox{{$p_\alpha$}}}}
\put(1,1.4){\makebox(0,0)[bl]{\hbox{{$p_\beta$}}}}
\put(1.8,0){\makebox(0,0)[cl]{\hbox{{$p_3$}}}}
\put(0.2,-1.2){\makebox(0,0)[cc]{\hbox{{$p_\gamma$}}}}
\put(1,-1.4){\makebox(0,0)[tl]{\hbox{{$p_\delta$}}}}
%right half
}
}
\end{pspicture}
\caption{\small
A sphere with five holes split into three pair of pants: one pair bounded by $\gamma_{p_\alpha}$, $\gamma_{p_\beta}$, and $\gamma_{A'}$, another by $\gamma_{A'}$, $\gamma_{p_3}$, and $\gamma_{A}$, and the third by $\gamma_{p_\gamma}$, $\gamma_{p_\delta}$, and $\gamma_{A}$.
}
\label{fi:five}
\end{figure}
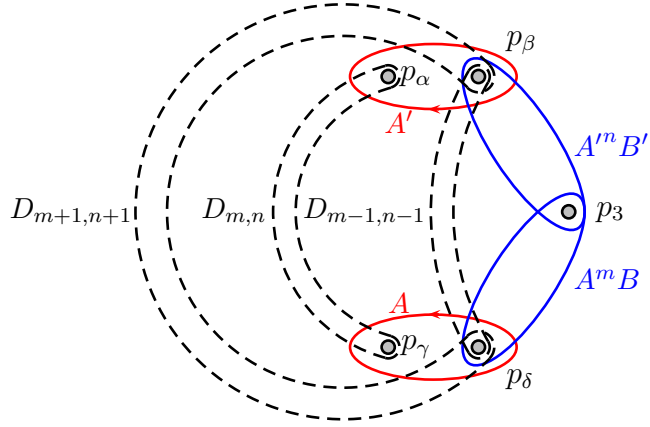

We begin with proving a technically most difficult case in which both $\widehat\tau_B$ and  $\widehat\tau_{B'}$ are twist coordinates of type (\ref{tau-can}). First, these coordinates obviously commute if four-holed spheres determining these coordinates do not share a common three-holed sphere (a pair of pants). So, our basic pattern is a five-holed sphere depicted in Fig.~\ref{fi:five}.

We use the following convention in Fig.~\ref{fi:five}: $\alpha=1$, $\beta=2$ if $n \in 2\mathbb Z$ and  $\alpha=2$, $\beta=1$ if $n \in 2\mathbb Z+1$; correspondingly,  $\gamma=5$, $\delta=4$ if $m \in 2\mathbb Z$ and  $\gamma=4$, $\delta=5$ if $m \in 2\mathbb Z+1$. Denoting by $D_{m,n}$ the geodesic functions of geodesics encircling corresponding holes (as indicated in the figure) and lying in the proper homotopy class, we apply the standard Goldman brackets and skein relations to obtain a set of useful formulas:
\be
\label{G1}
\{ G_{A^mB},G_{{A'}^n B'}\}=D_{m+1,n+1}-D_{m-1,n-1},
\ee
\be
\label{G2}
G_{A^mB} G_{{A'}^n B'}=D_{m+1,n+1} + D_{m-1,n-1} +G_3 D_{m,n}+G_\beta G_\delta,
\ee
\be
\label{G3}
G_{A'} D_{m,n}=D_{m,n+1} + D_{m,n-1} +G_{A^m B} G_\alpha+G_\beta G_\gamma,
\ee
\be
\label{G4}
G_{A} D_{m,n}=D_{m+1,n} + D_{m-1,n} +G_{{A'}^nB'} G_\gamma+G_\alpha G_\delta.
\ee
The first (and crucial) step is to find the brackets between non-normalized twists
\be\label{tB}
e^{\tau_B}:= G_{A^{-1}B}-e^{-\ell_A/2} G_B -\frac{G_5G_{A'} +G_3G_4 +e^{-\ell_A/2}(G_4G_{A'} +G_3G_5) }{e^{\ell_A/2}-e^{-\ell_A/2}}
\ee
and
\be\label{tBB}
e^{\tau_{B'}}:= G_{{A'}^{-1}B'}-e^{-\ell_{A'}/2} G_{B'} -\frac{G_1G_{A} +G_3G_2 +e^{-\ell_{A'}/2}(G_2G_{A} +G_3G_1) }{e^{\ell_{A'}/2}-e^{-\ell_{A'}/2}}.
\ee
Note that for the sake of a mirror-like symmetry in Fig.~\ref{fi:five}, we choose opposite orientations of $A$ and $A'$ cycles; correspondingly, in our convention, we have that $\{\tau_B,\ell_A\}=-\{\tau_{B'},\ell_{A'}\}=-2$.
A direct calculation gives
\begin{align*}
\{ e^{\tau_B}, e^{\tau_{B'}}\}&=\bigl[D_{0,0}-e^{-\ell_{A'}/2}D_{0,1}- e^{-\ell_{A}/2}D_{1,0}+e^{-\ell_A/2-\ell_{A'}/2}D_{1,1}\bigr]\\
&-\bigl[D_{-2,-2}-e^{-\ell_{A'}/2}D_{-2,-1}- e^{-\ell_{A}/2}D_{-1,-2}+e^{-\ell_A/2-\ell_{A'}/2}D_{-1,-1}\bigr]\\
&+(G_1+e^{-\ell_{A'}/2}G_2)\frac{e^{\ell_A/2}-e^{-\ell_A/2}}{e^{\ell_{A'}/2}-e^{-\ell_{A'}/2}}e^{\tau_B} +(G_5+e^{-\ell_{A}/2}G_4)\frac{e^{\ell_{A'}/2}-e^{-\ell_{A'}/2}}{e^{\ell_A/2}-e^{-\ell_A/2}} e^{\tau_{B'}}.
\end{align*}
The next set of formulas pertains to index ``shifts'' of blocks containing $D$-terms in the above formula:
\begin{align*}
&D_{-1,-1}-e^{-\ell_{A'}/2}D_{-1,0}- e^{-\ell_{A}/2}D_{0,-1}+e^{-\ell_A/2-\ell_{A'}/2}D_{0,0}\\
&\quad = e^{\ell_A/2 +\ell_{A'}/2}\bigl[D_{0,0}-e^{-\ell_{A'}/2}D_{0,1}- e^{-\ell_{A}/2}D_{1,0}+e^{-\ell_A/2-\ell_{A'}/2}D_{1,1}\bigr]\\
&\quad -G_1\bigl[G_{A^{-1}B} -e^{-\ell_{A}/2}G_B\bigr] - G_5\bigl[G_{{A'}^{-1}B'} - e^{-\ell_{A'}/2}G_{B'}\bigr] -G_1G_3G_5 -e^{\ell_{A}/2} G_2G_5-e^{\ell_{A'}/2}G_1G_4
\end{align*}
and
\begin{align*}
&D_{-2,-2}-e^{-\ell_{A'}/2}D_{-2,-1}- e^{-\ell_{A}/2}D_{-1,-2}+e^{-\ell_A/2-\ell_{A'}/2}D_{-1,-1}\\
&\quad = e^{\ell_A/2 +\ell_{A'}/2}\bigl[D_{-1,-1}-e^{-\ell_{A'}/2}D_{-1,0}- e^{-\ell_{A}/2}D_{0,-1}+e^{-\ell_A/2-\ell_{A'}/2}D_{0,0}\bigr]\\
&\quad -G_2\bigl[e^{\ell_{A}/2}G_{A^{-1}B} -G_B\bigr] - G_4\bigl[e^{\ell_{A'}/2}G_{{A'}^{-1}B'} - G_{B'}\bigr] +G_2G_3G_4 +e^{-\ell_{A}/2} G_1G_4 +e^{-\ell_{A'}/2}G_2G_5,
\end{align*}
using which we find that
\begin{align}
\{ e^{\tau_B}, e^{\tau_{B'}}\}&=(e^{-\ell_A/2-\ell_{A'}/2}-e^{\ell_A/2+\ell_{A'}/2})\Bigl(
\bigl[D_{-1,-1}-e^{-\ell_{A'}/2}D_{-1,0}- e^{-\ell_{A}/2}D_{0,-1}+e^{-\ell_A/2-\ell_{A'}/2}D_{0,0}\bigr] \Bigr. \nonumber\\
&-\bigl[G_{A^{-1}B} -e^{-\ell_{A}/2}G_B\bigr] \frac{e^{-\ell_{A'}/2}G_1+G_2}{e^{\ell_{A'}/2}-e^{-\ell_{A'}/2}} - \bigl[G_{{A'}^{-1}B'} -e^{-\ell_{A'}/2} G_{B'}\bigr] \frac{e^{-\ell_{A}/2}G_5+G_4}{e^{\ell_{A}/2}-e^{-\ell_{A}/2}}\nonumber\\
\Bigl. +&\frac{G_3\bigl(e^{-\ell_{A'}/2}G_1+G_2\bigr) \bigl( e^{-\ell_{A}/2}G_5+G_4\bigr) +2\bigl(e^{-\ell_{A'}/2}G_2+G_1\bigr) \bigl( e^{-\ell_{A}/2}G_4+G_5\bigr)  }{\bigl(e^{\ell_{A}/2}-e^{-\ell_{A}/2} \bigr) \bigl(e^{\ell_{A'}/2}-e^{-\ell_{A'}/2} \bigr)} \Bigr).\label{t-brak}
\end{align}
We obtain the last necessary formula from the skein relation (\ref{G2}): it turns out that its form exactly repeats (\ref{t-brak}), just with a different pre-factor:
\begin{align}
e^{\tau_B} e^{\tau_{B'}}&=(e^{-\ell_A/2-\ell_{A'}/2}+G_3+e^{\ell_A/2+\ell_{A'}/2})\Bigl(
\bigl[D_{-1,-1}-e^{-\ell_{A'}/2}D_{-1,0}- e^{-\ell_{A}/2}D_{0,-1}+e^{-\ell_A/2-\ell_{A'}/2}D_{0,0}\bigr] \Bigr. \nonumber\\
&-\bigl[G_{A^{-1}B} -e^{-\ell_{A}/2}G_B\bigr] \frac{e^{-\ell_{A'}/2}G_1+G_2}{e^{\ell_{A'}/2}-e^{-\ell_{A'}/2}} - \bigl[G_{{A'}^{-1}B'} -e^{-\ell_{A'}/2} G_{B'}\bigr] \frac{e^{-\ell_{A}/2}G_5+G_4}{e^{\ell_{A}/2}-e^{-\ell_{A}/2}}\nonumber\\
\Bigl. +&\frac{G_3\bigl(e^{-\ell_{A'}/2}G_1+G_2\bigr) \bigl( e^{-\ell_{A}/2}G_5+G_4\bigr) +2\bigl(e^{-\ell_{A'}/2}G_2+G_1\bigr) \bigl( e^{-\ell_{A}/2}G_4+G_5\bigr)  }{\bigl(e^{\ell_{A}/2}-e^{-\ell_{A}/2} \bigr) \bigl(e^{\ell_{A'}/2}-e^{-\ell_{A'}/2} \bigr)} \Bigr).\label{t-prod}
\end{align}
We have therefore proved a technical proposition.
\begin{proposition}\label{pr:tau}
The Poisson bracket between non-normalized twist coordinates $\tau_B$ and $\tau_{B'}$ defined by (\ref{tB}) and (\ref{tBB}) is
\be
\{\tau_B, \tau_{B'}\}=\frac{-\sinh(\ell_A/2+\ell_{A'}/2)}{\cosh(\ell_A/2+\ell_{A'}/2)+\cosh(p_3/2)}=\frac{e^{-\ell_A/2-\ell_{A'}/2-p_3/2}-e^{\ell_A/2 +\ell_{A'}/2-p_3/2} }{\bigl(e^{-\ell_A/2-\ell_{A'}/2-p_3/2}+1\bigr)\bigl(e^{\ell_A/2+\ell_{A'}/2-p_3/2}+1\bigr)}.
\ee
\end{proposition}
To complete the proof of the first case of the theorem, we present one more useful relation: For any three geodesic functions $G_{A_i}:=e^{\ell_{A_i}/2}+e^{-\ell_{A_i}/2}$, $i=1,2,3$, we have
\begin{align}
\label{TTTT}
G_{A_1}G_{A_2}G_{A_3}+G^2_{A_1}+G^2_{A_2}+G^2_{A_3}{-}4&=\bigl(e^{-\ell_{A_1}/2-\ell_{A_2}/2-\ell_{A_3}/2} +1\bigr)\bigl(e^{\ell_{A_1}/2+\ell_{A_2}/2-\ell_{A_3}/2} +1\bigr)\\
\times& \bigl(e^{\ell_{A_1}/2+\ell_{A_3}/2-\ell_{A_2}/2} +1\bigr)\bigl(e^{\ell_{A_2}/2+\ell_{A_3}/2-\ell_{A_1}/2} +1\bigr)\nonumber
\end{align}
And therefore
\begin{align*}
&\Bigl[\frac{\partial}{\partial \ell_A}+ \frac{\partial}{\partial \ell_{A'}}\Bigr] \log \bigl( G_{A}G_{A'}G_{p_3}+G^2_{A}+G^2_{A'}+G^2_{p_3}{-}4\bigr)\\
&=\frac{e^{\ell_{A}/2+\ell_{A'}/2-p_3/2}}{e^{\ell_{A}/2+\ell_{A'}/2-p_3/2} +1}-\frac{e^{-\ell_{A}/2-\ell_{A'}/2-p_3/2}}{e^{-\ell_{A}/2-\ell_{A'}/2-p_3/2} +1}
=\frac{e^{\ell_{A}/2+\ell_{A'}/2-p_3/2}-e^{-\ell_{A}/2-\ell_{A'}/2-p_3/2}}{\bigl(e^{\ell_{A}/2+\ell_{A'}/2-p_3/2} +1\bigr)\bigl(e^{-\ell_{A}/2-\ell_{A'}/2-p_3/2} +1\bigr)},
\end{align*}
which is exactly the expression in the right-hand side of Proposition~\ref{pr:tau}. We have therefore proved the first case of the theorem statement.

\begin{figure}[tb]
\begin{pspicture}(-2.5,-2.5)(2.5,2.5)
\newcommand{\PATTERN}{%
{\psset{unit=1}
\put(-1.5,-1.5){\pscircle[fillstyle=solid,fillcolor=lightgray]{.5}}
\put(-1.5,1.5){\pscircle[linecolor=red]{.5}}
\put(1.5,-1.5){\pscircle[fillstyle=solid,fillcolor=lightgray]{.5}}
\put(1.5,1.5){\pscircle[linecolor=red]{.5}}
\psline[linewidth=1.5pt,linecolor=blue](-1,1.5)(1,1.5)
\rput{45}(0,0){\psarc[linewidth=1pt,linecolor=magenta](-2.13,0){.65}{90}{270}
\psarc[linewidth=1pt,linecolor=magenta](2.13,0){.65}{-90}{90}
\psline[linewidth=1pt,linecolor=magenta](-2.13,0.65)(2.13,0.65)
\psline[linewidth=1pt,linecolor=magenta](-2.13,-0.65)(2.13,-0.65)
}
\psarc[linewidth=1pt,linecolor=blue](1.5,-1.5){0.75}{-90}{90}
\psarc[linewidth=1pt,linecolor=blue](-1.5,-1.5){0.75}{90}{270}
\psline[linewidth=1pt,linecolor=blue](-1.5,-.75)(1.5,-.75)
\psline[linewidth=1pt,linecolor=blue](-1.5,-2.25)(1.5,-2.25)
%
%\psarc[linewidth=1pt,linecolor=red](1.5,1.5){0.8}{0}{180}
%\psarc[linewidth=1pt,linecolor=red](1.5,-1.5){0.8}{180}{360}
%\psline[linewidth=1pt,linecolor=red](0.7,1.5)(0.7,-1.5)
%\psline[linewidth=1pt,linecolor=red](2.3,1.5)(2.3,-1.5)
%
\put(-1.5,1.5){\makebox(0,0)[cc]{\hbox{{$A$}}}}
\put(1.5,1.5){\makebox(0,0)[cc]{\hbox{{$A$}}}}
\put(-0.6,1.6){\makebox(0,0)[bc]{\hbox{{$B$}}}}
\put(1.5,-1.5){\makebox(0,0)[cc]{\hbox{{$p_3$}}}}
\put(-1.5,-1.5){\makebox(0,0)[cc]{\hbox{{$p_4$}}}}
\put(0,-2){\makebox(0,0)[bc]{\hbox{$A'$}}}
%\put(2.4,0){\makebox(0,0)[cl]{\hbox{\tcr{$B'{A'}^{-1}$}}}}
\put(-0.9,0.3){\makebox(0,0)[cr]{\hbox{$B'$}}}
\put(0,-2.4){\makebox(0,0)[tc]{\hbox{$\{G_B,G_{B'}\}$}}}
%\psbezier[linecolor=blue](3.366,-1.5)(2.866,-0.634)(3.19,0.734)(3.69,1.6)
%\psbezier[linecolor=blue](-2.47,-3.03)(-1.47,-2.03)(1.134,-1.634)(1.634,-2.5)
%\psbezier[linecolor=blue](-2.65,3.05)(-1.65,2.05)(1.81,1.534)(2.31,2.4)
}
}
\rput(0,0){\PATTERN
\put(0.58,1.5){\pscircle[linestyle=dashed, linewidth=1pt]{.3}}
}
\end{pspicture}
\begin{pspicture}(-2.5,-2.5)(2.5,2.5)
\newcommand{\PATTERN}{%
{\psset{unit=1}
\put(-1.5,-1.5){\pscircle[fillstyle=solid,fillcolor=lightgray]{.5}}
\put(-1.5,1.5){\pscircle[linecolor=red]{.5}}
\put(1.5,-1.5){\pscircle[fillstyle=solid,fillcolor=lightgray]{.5}}
\put(1.5,1.5){\pscircle[linecolor=red]{.5}}
\psline[linewidth=1.5pt,linecolor=blue](-1,1.5)(1,1.5)
\rput{45}(0,0){\psarc[linewidth=1pt,linecolor=magenta](-2.13,0){.65}{90}{270}
\psarc[linewidth=1pt,linecolor=magenta](2.13,0){.65}{-90}{90}
\psline[linewidth=1pt,linecolor=magenta](-2.13,0.65)(2.13,0.65)
\psline[linewidth=1pt,linecolor=magenta](-2.13,-0.65)(2.13,-0.65)
}
\psarc[linewidth=1pt,linecolor=blue](1.5,-1.5){0.75}{-90}{90}
\psarc[linewidth=1pt,linecolor=blue](-1.5,-1.5){0.75}{90}{270}
\psline[linewidth=1pt,linecolor=blue](-1.5,-.75)(1.5,-.75)
\psline[linewidth=1pt,linecolor=blue](-1.5,-2.25)(1.5,-2.25)
%
%\psarc[linewidth=1pt,linecolor=red](1.5,1.5){0.8}{0}{180}
%\psarc[linewidth=1pt,linecolor=red](1.5,-1.5){0.8}{180}{360}
%\psline[linewidth=1pt,linecolor=red](0.7,1.5)(0.7,-1.5)
%\psline[linewidth=1pt,linecolor=red](2.3,1.5)(2.3,-1.5)
%
\put(-1.5,1.5){\makebox(0,0)[cc]{\hbox{{$A$}}}}
\put(1.5,1.5){\makebox(0,0)[cc]{\hbox{{$A$}}}}
%\put(-0.6,1.6){\makebox(0,0)[bc]{\hbox{{$B$}}}}
\put(1.5,-1.5){\makebox(0,0)[cc]{\hbox{{$p_3$}}}}
\put(-1.5,-1.5){\makebox(0,0)[cc]{\hbox{{$p_4$}}}}
\put(0,-2){\makebox(0,0)[bc]{\hbox{$A'$}}}
\put(-2,0){\makebox(0,0)[cr]{\hbox{$=1/2$}}}
\put(2.3,0){\makebox(0,0)[cl]{\hbox{$-1/2$}}}
%\put(2.4,0){\makebox(0,0)[cl]{\hbox{\tcr{$B'{A'}^{-1}$}}}}
%\put(-0.9,0.3){\makebox(0,0)[cr]{\hbox{$B'$}}}
%\psbezier[linecolor=blue](3.366,-1.5)(2.866,-0.634)(3.19,0.734)(3.69,1.6)
%\psbezier[linecolor=blue](-2.47,-3.03)(-1.47,-2.03)(1.134,-1.634)(1.634,-2.5)
%\psbezier[linecolor=blue](-2.65,3.05)(-1.65,2.05)(1.81,1.534)(2.31,2.4)
}
}
\rput(0,0){\PATTERN
\put(0.58,1.5){\pscircle[linestyle=dashed, linewidth=1pt,fillstyle=solid, fillcolor=white]{.3}
\psarc[linewidth=1pt,linecolor=blue](-0.3,0.724){0.724}{270}{315}
\psarc[linewidth=1pt,linecolor=blue](0.3,-0.724){0.724}{90}{135}}
}
\end{pspicture}
\begin{pspicture}(-2.5,-2.5)(2.5,2.5)
\newcommand{\PATTERN}{%
{\psset{unit=1}
\put(-1.5,-1.5){\pscircle[fillstyle=solid,fillcolor=lightgray]{.5}}
\put(-1.5,1.5){\pscircle[linecolor=red]{.5}}
\put(1.5,-1.5){\pscircle[fillstyle=solid,fillcolor=lightgray]{.5}}
\put(1.5,1.5){\pscircle[linecolor=red]{.5}}
\psline[linewidth=1.5pt,linecolor=blue](-1,1.5)(1,1.5)
\rput{45}(0,0){\psarc[linewidth=1pt,linecolor=magenta](-2.13,0){.65}{90}{270}
\psarc[linewidth=1pt,linecolor=magenta](2.13,0){.65}{-90}{90}
\psline[linewidth=1pt,linecolor=magenta](-2.13,0.65)(2.13,0.65)
\psline[linewidth=1pt,linecolor=magenta](-2.13,-0.65)(2.13,-0.65)
}
\psarc[linewidth=1pt,linecolor=blue](1.5,-1.5){0.75}{-90}{90}
\psarc[linewidth=1pt,linecolor=blue](-1.5,-1.5){0.75}{90}{270}
\psline[linewidth=1pt,linecolor=blue](-1.5,-.75)(1.5,-.75)
\psline[linewidth=1pt,linecolor=blue](-1.5,-2.25)(1.5,-2.25)
%
%\psarc[linewidth=1pt,linecolor=red](1.5,1.5){0.8}{0}{180}
%\psarc[linewidth=1pt,linecolor=red](1.5,-1.5){0.8}{180}{360}
%\psline[linewidth=1pt,linecolor=red](0.7,1.5)(0.7,-1.5)
%\psline[linewidth=1pt,linecolor=red](2.3,1.5)(2.3,-1.5)
%
\put(-1.5,1.5){\makebox(0,0)[cc]{\hbox{{$A$}}}}
\put(1.5,1.5){\makebox(0,0)[cc]{\hbox{{$A$}}}}
%\put(-0.6,1.6){\makebox(0,0)[bc]{\hbox{{$B$}}}}
\put(1.5,-1.5){\makebox(0,0)[cc]{\hbox{{$p_3$}}}}
\put(-1.5,-1.5){\makebox(0,0)[cc]{\hbox{{$p_4$}}}}
\put(0,-2){\makebox(0,0)[bc]{\hbox{$A'$}}}
%\put(2.4,0){\makebox(0,0)[cl]{\hbox{\tcr{$B'{A'}^{-1}$}}}}
%\put(-0.9,0.3){\makebox(0,0)[cr]{\hbox{$B'$}}}
%\psbezier[linecolor=blue](3.366,-1.5)(2.866,-0.634)(3.19,0.734)(3.69,1.6)
%\psbezier[linecolor=blue](-2.47,-3.03)(-1.47,-2.03)(1.134,-1.634)(1.634,-2.5)
%\psbezier[linecolor=blue](-2.65,3.05)(-1.65,2.05)(1.81,1.534)(2.31,2.4)
}
}
\rput(0,0){\PATTERN
\put(0.58,1.5){\pscircle[linestyle=dashed, linewidth=1pt,fillstyle=solid, fillcolor=white]{.3}
\psarc[linewidth=1pt,linecolor=blue](-0.3,-0.125){0.125}{-45}{90}
\psarc[linewidth=1pt,linecolor=blue](0.3,0.125){0.125}{135}{270}}
}
\end{pspicture}
\begin{pspicture}(-2.5,-2.5)(2.5,2.5)
%%% EMPTY PICTURE
\end{pspicture}
\begin{pspicture}(-2.5,-2.5)(2.5,2.5)
\newcommand{\PATTERN}{%
{\psset{unit=1}
\put(-1.5,-1.5){\pscircle[fillstyle=solid,fillcolor=lightgray]{.5}}
\put(-1.5,1.5){\pscircle[linecolor=red]{.5}}
\put(1.5,-1.5){\pscircle[fillstyle=solid,fillcolor=lightgray]{.5}}
\put(1.5,1.5){\pscircle[linecolor=red]{.5}}
\psline[linewidth=1.5pt,linecolor=blue](-1,1.5)(1,1.5)
\rput{45}(0,0){\psarc[linewidth=1pt,linecolor=magenta](-2.13,0){.65}{90}{270}
\psarc[linewidth=1pt,linecolor=magenta](2.13,0){.65}{-90}{90}
\psline[linewidth=1pt,linecolor=magenta](-2.13,0.65)(2.13,0.65)
\psline[linewidth=1pt,linecolor=magenta](-2.13,-0.65)(2.13,-0.65)
}
\psarc[linewidth=1pt,linecolor=blue](1.5,-1.5){0.75}{-90}{90}
\psarc[linewidth=1pt,linecolor=blue](-1.5,-1.5){0.75}{90}{270}
\psline[linewidth=1pt,linecolor=blue](-1.5,-.75)(1.5,-.75)
\psline[linewidth=1pt,linecolor=blue](-1.5,-2.25)(1.5,-2.25)
%
%\psarc[linewidth=1pt,linecolor=red](1.5,1.5){0.8}{0}{180}
%\psarc[linewidth=1pt,linecolor=red](1.5,-1.5){0.8}{180}{360}
%\psline[linewidth=1pt,linecolor=red](0.7,1.5)(0.7,-1.5)
%\psline[linewidth=1pt,linecolor=red](2.3,1.5)(2.3,-1.5)
%
\put(-1.5,1.5){\makebox(0,0)[cc]{\hbox{{$A$}}}}
\put(1.5,1.5){\makebox(0,0)[cc]{\hbox{{$A$}}}}
%\put(-0.6,1.6){\makebox(0,0)[bc]{\hbox{{$B$}}}}
\put(1.5,-1.5){\makebox(0,0)[cc]{\hbox{{$p_3$}}}}
\put(-1.5,-1.5){\makebox(0,0)[cc]{\hbox{{$p_4$}}}}
\put(0,-2){\makebox(0,0)[bc]{\hbox{$A'$}}}
\put(-2,0){\makebox(0,0)[cr]{\hbox{$=1/2$}}}
\put(2.3,0){\makebox(0,0)[cl]{\hbox{$-1/2$}}}
%\put(2.4,0){\makebox(0,0)[cl]{\hbox{\tcr{$B'{A'}^{-1}$}}}}
%\put(-0.9,0.3){\makebox(0,0)[cr]{\hbox{$B'$}}}
%\psbezier[linecolor=blue](3.366,-1.5)(2.866,-0.634)(3.19,0.734)(3.69,1.6)
%\psbezier[linecolor=blue](-2.47,-3.03)(-1.47,-2.03)(1.134,-1.634)(1.634,-2.5)
%\psbezier[linecolor=blue](-2.65,3.05)(-1.65,2.05)(1.81,1.534)(2.31,2.4)
}
}
\rput(0,0){\PATTERN
\put(0.58,1.5){\pscircle[linestyle=dashed, linewidth=1pt,fillstyle=solid, fillcolor=white]{.3}
\psarc[linewidth=1pt,linecolor=blue](-0.3,0.724){0.724}{270}{315}
\psarc[linewidth=1pt,linecolor=blue](0.3,-0.724){0.724}{90}{135}}
}
\end{pspicture}
%
%right half
\begin{pspicture}(-2.5,-2.5)(2.5,2.5)
\newcommand{\PATTERN}{%
{\psset{unit=1}
\put(-1.5,-1.5){\pscircle[fillstyle=solid,fillcolor=lightgray]{.5}}
\put(-1.5,1.5){\pscircle[linecolor=red]{.5}}
\put(1.5,-1.5){\pscircle[fillstyle=solid,fillcolor=lightgray]{.5}}
\put(1.5,1.5){\pscircle[linecolor=red]{.5}}
\psline[linewidth=1.5pt,linecolor=blue](-0.4,1.5)(0.4,1.5)
%
%\rput{45}(0,0){\psarc[linewidth=1pt,linecolor=magenta](-2.13,0){.65}{90}{270}
%\psarc[linewidth=1pt,linecolor=magenta](2.13,0){.65}{-90}{90}
%\psline[linewidth=1pt,linecolor=magenta](-2.13,0.65)(2.13,0.65)
%\psline[linewidth=1pt,linecolor=magenta](-2.13,-0.65)(2.13,-0.65)
%}
\psarc[linewidth=1pt,linecolor=blue](1.5,-1.5){0.75}{-90}{90}
\psarc[linewidth=1pt,linecolor=blue](-1.5,-1.5){0.75}{90}{270}
\psline[linewidth=1pt,linecolor=blue](-1.5,-.75)(1.5,-.75)
\psline[linewidth=1pt,linecolor=blue](-1.5,-2.25)(1.5,-2.25)
\psarc[linewidth=1pt,linecolor=red](1.5,1.5){0.9}{0}{180}
\psarc[linewidth=1pt,linecolor=red](1.5,-1.5){0.9}{180}{360}
\psline[linewidth=1pt,linecolor=red](0.6,1.5)(0.6,-1.5)
\psline[linewidth=1pt,linecolor=red](2.4,1.5)(2.4,-1.5)
\put(-1.5,1.5){\makebox(0,0)[cc]{\hbox{{$A$}}}}
\put(1.5,1.5){\makebox(0,0)[cc]{\hbox{{$A$}}}}
\put(-0.8,1.4){\makebox(0,0)[tl]{\hbox{{$A^{-2}B$}}}}
\put(1.5,-1.5){\makebox(0,0)[cc]{\hbox{{$p_3$}}}}
\put(-1.5,-1.5){\makebox(0,0)[cc]{\hbox{{$p_4$}}}}
\put(0,-2){\makebox(0,0)[bc]{\hbox{$A'$}}}
\put(2.5,0){\makebox(0,0)[cl]{\hbox{${A'}^{-1}B'$}}}
%\put(-0.9,0.3){\makebox(0,0)[cr]{\hbox{$B'$}}}
%\psbezier[linecolor=blue](3.366,-1.5)(2.866,-0.634)(3.19,0.734)(3.69,1.6)
%\psbezier[linecolor=blue](-2.47,-3.03)(-1.47,-2.03)(1.134,-1.634)(1.634,-2.5)
%\psbezier[linecolor=blue](-2.65,3.05)(-1.65,2.05)(1.81,1.534)(2.31,2.4)
}
}
\rput(0,0){\PATTERN
%\put(0.58,1.5){\pscircle[linestyle=dashed, linewidth=1pt]{.3}}
\rput(-1.5,1.5){\parametricplot[linecolor=blue,linestyle=solid,linewidth=1pt]{0}{1}{0.5 0.6 t mul add 360 t t t mul mul t -2 add mul  t -2 add mul  t -2 add mul  mul cos mul  0.5 0.6 t mul add 360 t t t mul mul t -2 add mul  t -2 add mul  t -2 add mul  mul sin mul}}
\rput{180}(1.5,1.5){\parametricplot[linecolor=blue,linestyle=solid,linewidth=1pt]{0}{1}{0.5 0.6 t mul add 360 t t t mul mul t -2 add mul  t -2 add mul  t -2 add mul  mul cos mul  0.5 0.6 t mul add 360 t t t mul mul t -2 add mul  t -2 add mul  t -2 add mul  mul sin mul}}
\put(0.6,1.1){\pscircle[linestyle=dashed, linewidth=1pt,fillstyle=solid, fillcolor=white]{.3}
\psarc[linewidth=1pt,linecolor=blue](-0.07,0.3){0.07}{210}{360}
\psarc[linewidth=1pt,linecolor=blue](0.1,-0.3){0.1}{30}{180}}
}
\end{pspicture}
\caption{\small
A torus with two holes: we cut it along an $A$-cycle obtaining a sphere with four holes. The $B$-cycle is then a line between two copies of the $A$-cycle. Cycles $B$ and $B'$ have a single intersection point at which we have two resolutions of the crossing. We denote $D_{n,m}$ the first-type resolution of the crossing between cycles $A^nB$ and ${A'}^mB'$. It turns out that the second-type resolution of the crossing between $B$ and $B'$  is homotopically equivalent to the first-type resolution of the crossing between $A^2B$ and ${A'}^{-1}B'$, so the second term is $D_{n-2,m-1}$.
}
\label{fi:torus-2holes}
\end{figure}
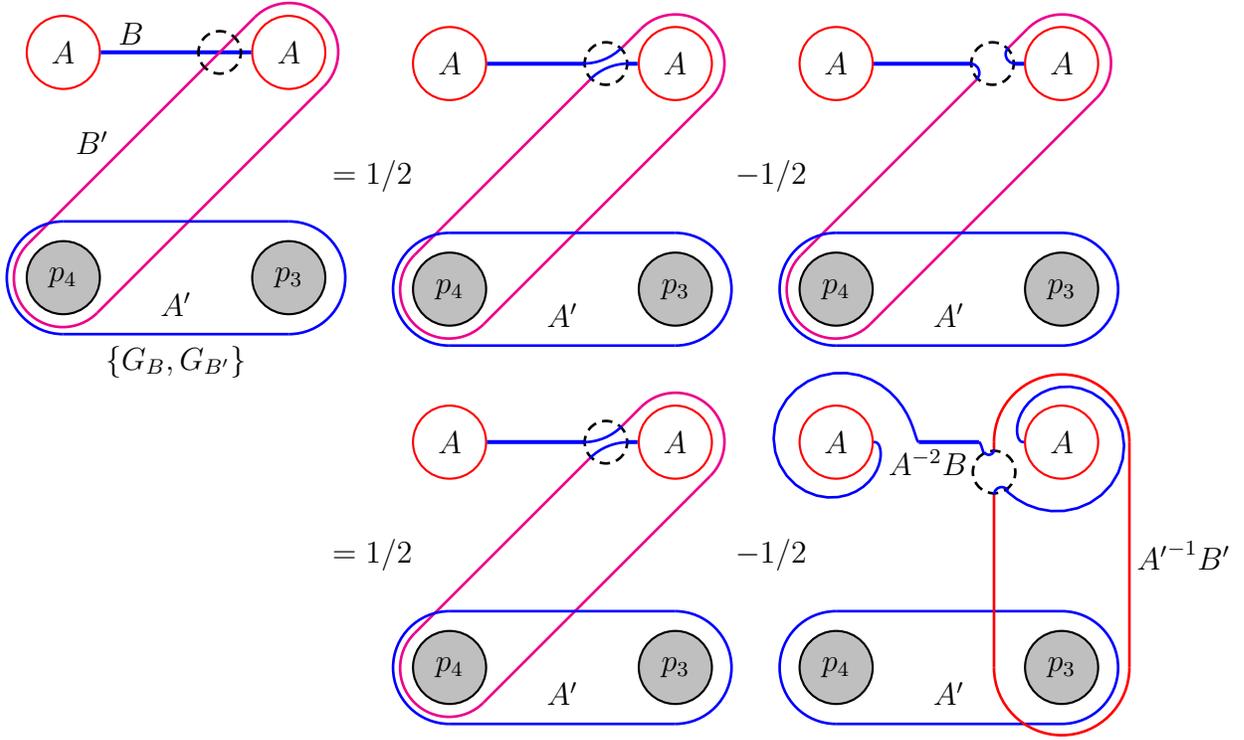

The second case is where $\widehat\tau_B$ is a twist coordinate of type  (\ref{tau-torus}) and  $\widehat\tau_{B'}$ is a twist coordinates of type (\ref{tau-can}). We then cut the torus along an $A$-cycle thus obtaining a sphere with four holes, two of which are copies of the $A$-cycle; the $B$-cycle is then an interval joining these two copies, see Fig.~\ref{fi:torus-2holes}. In this case, cycles $A^nB$ and ${A'}^mB'$ have a single intersection point, and as in the first case, we introduce cycles $D_{n,m}$ obtained by the first-type resolution at this point. In Fig.~\ref{fi:torus-2holes} we demonstrate that the second-type resolution of the crossing between $\gamma_{A^nB}$ and $\gamma_{{A'}^mB'}$ is homotopically equivalent to the first-type resolution of the crossing between $\gamma_{A^{n-2}B}$ and $\gamma_{{A'}^{m-1}B'}$, so
\be
\{G_{A^nB},G_{{A'}^mB'}\}=\frac12 \bigl(D_{n,m}-D_{n-2,m-1} \bigr),
\ee
whereas the skein relation reads
\be
G_{A^nB}G_{{A'}^mB'}=D_{n,m} + D_{n-2,m-1}.
\ee
Note that in this case $S_{AB}=S_B=G_A(G_3+G_4)$, so we have uniform formulas for all $\tau_{{A'}^mB'}$.

Omitting details, we have the following proposition.
\begin{proposition}\label{pr:torus-tau}
The Poisson bracket between non-normalized twist coordinates $\tau_{A^nB}$ and $\tau_{{A'}^mB'}$ such that
\begin{align*}
e^{\tau_{A^nB}}:= &G_{{A}^{n-1}B}-e^{-\ell_{A}/2} G_{A^nB'},\\
e^{\tau_{{A'}^mB'}}:= &G_{{A'}^{m-1}B'}-e^{-\ell_{A'}/2} G_{{A'}^mB'} -\frac{G_{A}(G_3+G_4)}{e^{\ell_{A'}/2}-1},
\end{align*}
is
\be
2\{\tau_{A^nB}, \tau_{{A'}^mB'}\}= \frac{1-e^{\ell_{A'}/2+\ell_A}}{1+e^{\ell_{A'}/2+\ell_A}}.
\ee
\end{proposition}
To come to the canonical twist coordinates (\ref{tau-torus}) and (\ref{tau-can}), we have to replace $p$ by $\ell_{A'}$ in (\ref{tau-torus}) and shift the above variables by
\begin{align*}
\widehat\tau_{A^nB}&=\tau_{A^nB}-\frac12\log\Bigl[(e^{-\ell_{A'}/2+\ell_A}+1)(e^{\ell_{A'}/2+\ell_A}+1)e^{-\ell_A}\Bigr],\\
\widehat\tau_{{A'}^mB'}&=\tau_{{A'}^mB'}-\frac12\log\Bigl[(e^{-\ell_{A'}/2+\ell_A}+1)(e^{-\ell_{A'}/2-\ell_A}+1)(e^{\ell_{A'}/2} +1)^2\Bigr].
\end{align*}
it is now an easy exercise to verify, using the constant brackets between $\tau$- and $\ell$-variables, that adding this terms results in the vanishing commutation relations between the canonical twist coordinates. The theorem is therefore proved.

\section{Fenchel--Nielsen coordinates for $\Sigma_{g,s,n}$}\label{s:arcs}
\setcounter{equation}{0}

We now generalize the Fenchel--Nielsen coordinate setting to the case of surfaces with marked points on boundaries. The standard trick making geometrical lengths finite is to introduce a regularisation by decorating all bordered cusps with horocycles. Note that all geodesic functions are insensitive to these decorations, which therefore affect only $\lambda$-lengths of arcs. We restrict consideration in this paper to subalgebras of $\lambda$-lengths combinations that are decoration-independent. On the language of Teichm\"uller spaces ${\mathfrak T}_{g,s,n}$, this corresponds to considering subalgebras of shear coordinates $Z_\alpha$ removing the extended shear coordinates $\pi_j$ from consideration. The basic example of this construction is a sphere with three holes and bordered cusps situated on the boundary of one of the holes.

\subsection{Fenchel--Nielsen coordinates for $\Sigma_{0,3,m}$}\label{ss:FN-03n}

Consider the case of a sphere with three holes and with $m$ bordered cusps located at the boundary of one of holes. All these cusps are endowed with horocycle decorations. If we restrict the phase space to a subspace of objects independent on decorations, then, likewise the cases of $\Sigma_{0,4}$ and $\Sigma_{1,1}$, we can single out a canonical twist coordinate $\widehat\tau_B$ dual to $\ell_A$, where $\ell_A$ is the length of a unique closed geodesic (the hole perimeter) separating the boundary component with cusps from the rest of the Riemann surface.

We show below that the remaining $m-1$ coordinates can be chosen to have homogeneous constant brackets between themselves commuting with all canonical length and twist coordinates for the rest of the surface. Their algebra is nondegenerate for odd $m$ and has exactly one Casimir element for even $m$.

\begin{figure}[tb]
\begin{pspicture}(-3,-4)(4,4){
\newcommand{\PATTERN}{%
{\psset{unit=1}
\rput{-45}(-3,-2.5){\psellipse[linecolor=blue, linewidth=1pt](0,0)(.75,.5)}%\psframe[linecolor=white, fillstyle=solid, fillcolor=white](-1,0)(1,.6)\psellipse[linecolor=blue, linestyle=dashed,  linewidth=1pt](0,0)(.75,.5)}
\rput{30}(2.5,-2){\psellipse[linecolor=blue, linewidth=1pt](0,0)(1,.5)}%\psframe[linecolor=white, fillstyle=solid, fillcolor=white](-1.1,0)(1.1,.6)\psellipse[linecolor=blue, linestyle=dashed, linewidth=1pt](0,0)(1,.5)}
\rput{-5}(0.29,0){\psellipse[linecolor=blue,  linewidth=1pt](0,0)(2.85,0.3)\psframe[linecolor=white, fillstyle=solid, fillcolor=white](-2.9,0)(2.9,.6)\psellipse[linecolor=blue, linestyle=dashed,  linewidth=1pt](0,0)(2.85,.3)}
%\rput{45}(-3,2.7){\psellipse[linecolor=blue,  linewidth=1pt](0,0)(.5,.3)}
%\rput{-30}(3,2){\psellipse[linecolor=blue,  linewidth=1pt](0,0)(.8,.3)}
\psbezier[linecolor=blue](-2.47,-3.03)(-1.47,-2.03)(1.134,-1.634)(1.634,-2.5)
\psbezier[linecolor=blue](-3.53,-1.97)(-2.56,-1)(-2,1)(-3.35,2.35)
\psbezier[linecolor=blue](3.35,-1.5)(2.85,-0.6)(3.19,1.2)(3.69,2)
% the crown
\psbezier[linecolor=blue](-1.5,1.6)(-1.56,0.6)(-2.35,1.35)(-3.35,2.35)
\psbezier[linecolor=blue](1.5,1.3)(1.56,0.4)(3.19,1.2)(3.69,2)
\psbezier[linecolor=blue](1.5,1.3)(1.46,0.4)(-1.46,0.6)(-1.5,1.6)
\psbezier[linecolor=blue](-0.7,2.8)(-0.76,1.8)(-2.35,1.35)(-3.35,2.35)
\psbezier[linecolor=blue](2,2.6)(2.06,1.6)(3.19,1.2)(3.69,2)
\psbezier[linecolor=blue](2,2.6)(2,1.6)(-0.66,1.8)(-0.7,2.8)
% decorations
\pscircle[linewidth=1pt,linestyle=dashed](-3.15,2.15){0.3}
\pscircle[linewidth=1pt,linestyle=dashed](-1.5,1.45){0.15}
\pscircle[linewidth=1pt,linestyle=dashed](-0.7,2.5){0.3}
\pscircle[linewidth=1pt,linestyle=dashed](1.5,1.1){0.2}
\pscircle[linewidth=1pt,linestyle=dashed](2,2.3){0.3}
\pscircle[linewidth=1pt,linestyle=dashed](3.55,1.86){0.2}
%\psbezier[linecolor=blue](-2.65,3.05)(-1.65,2.05)(1.81,1.534)(2.31,2.4)
}
}
\rput(-3,0){\rput(0,0){\PATTERN}
\rput{0}(0,0){\psbezier[linecolor=red](1.44,1.27)(1.42,0.4)(.75,-1.95)(1.15,-2.15)
\psbezier[linecolor=red](1.44,1.27)(1.46,0.4)(3,-0.4)(3.1,-1)
\psbezier[linecolor=red](1.15,-2.15)(1.55,-2.35)(3.3,-2.2)(3.1,-1)
\rput(1.44,-2.23){\makebox(0,0)[cc]{\hbox{\tcw{\tiny$\bullet$}}}}
\rput(3.1,-1.25){\makebox(0,0)[cc]{\hbox{\tcw{\tiny$\bullet$}}}}
% perpendiculars
\rput(0,0){\psbezier[linecolor=black](1.44,1.27)(1.47,0.7)(1.47,0.7)(1.6,-0.35)
\psbezier[linecolor=red](1.44,1.27)(1.5,0.7)(2,-0.8)(2.4,-1.5)
\psbezier[linecolor=black](2.3,-0.37)(2.4,-0.7)(2.4,-0.7)(2.65,-1.38)
\psbezier[linecolor=black](-1.6,-0.05)(-1.9,-0.8)(-2.3,-1.8)(-2.7,-2.2)
\psline[linecolor=black](2.6,-2)(2.8,-2.35)
\psline[linecolor=black, linewidth=0.5pt, linestyle=dashed](-2,0.37)(-2,0.6)}
}
\rput{-5}(0.29,0.3){\parametricplot[linecolor=green,linestyle=dashed,linewidth=1pt]{0}{180}{2.8 t cos mul  0.15 t sin  mul}
\psbezier[linecolor=green](1.1,1.05)(1.2,0.3)(2.75,-0.5)(2.8,0)
\psbezier[linecolor=green](1.05,1.05)(1.1,0.2)(-2.95,-0.5)(-2.9,0)
}
\rput(-1,0){
\rput(2.4,-1.6){\makebox(0,0)[cc]{\hbox{\tcr{\small$B$}}}}
%\put(2.4,-1.6){\pscircle[linecolor=green]{.2}}
\rput(-0.4,0.8){\makebox(0,0)[cc]{\hbox{\tcg{\small$\lambda_0$}}}}
%\put(2.3,1.7){\pscircle[linecolor=red]{.2}}
\rput(0.5,-0.5){\makebox(0,0)[cc]{\hbox{\tcb{\small$A$}}}}
\rput(-2,-2.5){\makebox(0,0)[cc]{\hbox{\small$p_1$}}}
\rput(4.4,-2){\makebox(0,0)[cl]{\hbox{\small$p_2$}}}
\rput(2.6,-0.45){\makebox(0,0)[tc]{\hbox{\tiny$P$}}}
\rput(3.2,-0.33){\makebox(0,0)[bc]{\hbox{\tiny$Q$}}}
\rput(-0.5,-0.33){\makebox(0,0)[bc]{\hbox{\tiny$R$}}}
\rput(2.35,0.6){\makebox(0,0)[cr]{\hbox{\tiny$h_1$}}}
\rput(3.5,-0.8){\makebox(0,0)[cl]{\hbox{\tiny$h_2$}}}
\rput(-1.2,-1.2){\makebox(0,0)[cr]{\hbox{\tiny$h'_2$}}}
\rput(3.8,-2.45){\makebox(0,0)[tc]{\hbox{\tiny$h_3$}}}
\rput(-1,0.7){\makebox(0,0)[br]{\hbox{\tiny$h_4$}}}
\rput(3,-1){\makebox(0,0)[cr]{\tcr{\hbox{\tiny$h_B$}}}}
%\put(0,-0.5){\pscircle[linecolor=blue]{.2}}
}
}
}
%right half
\rput(2,-0.3){
\psline[linecolor=red](0,3)(0,-2)
\psline[linecolor=blue](-1,0.4)(1,-0.4)
\psbezier[linecolor=blue](0,-2)(0.2,-2)(0.3,-2.04)(0.5,-2.1)
\psbezier[linecolor=black](0.5,-2.1)(0.8,-1.1)(0.8,-0.9)(1,-0.4)
\psbezier[linecolor=black](0.5,-2.1)(0.8,-1.1)(0.8,-0.9)(1,-0.4)
\psbezier[linecolor=black](0,3)(0,2.5)(-0.6,1.4)(-1,0.4)
\pscircle[linewidth=1pt,linestyle=dashed](0,2.8){0.2}
\rput(-1,0.35){\makebox(0,0)[tr]{\hbox{\tiny$P$}}}
\rput(1,-0.45){\makebox(0,0)[tl]{\hbox{\tiny$Q$}}}
\rput(-0.55,1.6){\makebox(0,0)[cr]{\hbox{\tiny$h_1$}}}
\rput(0.75,-1.4){\makebox(0,0)[cl]{\hbox{\tiny$h_2$}}}
\rput(-0.1,-1){\makebox(0,0)[cr]{\tcr{\hbox{\tiny$h_B$}}}}
}
\rput(6,0){
\psline[linecolor=red](0,3)(0,0)
\psline[linecolor=blue](-1.5,0)(0,0)
\psline[linecolor=black](-1.5,0)(-1.5,0.5)
\psbezier[linecolor=red](-1.5,0.5)(-0.7,0.5)(0,2)(0,3)
\pscircle[linewidth=1pt,linestyle=dashed](0,2.8){0.2}
%\rput(-1,0.35){\makebox(0,0)[tr]{\hbox{\tiny$P$}}}
\rput(-0.6,1.3){\makebox(0,0)[cr]{\hbox{\small$\ell_{B}/2$}}}
\rput(-0.75,-0.1){\makebox(0,0)[tc]{\hbox{\small$p_2/2$}}}
\rput(-1.6,0.25){\makebox(0,0)[cr]{\hbox{\small$h_3$}}}
\rput(0.1,1.5){\makebox(0,0)[cl]{\tcr{\hbox{\small$h_B$}}}}
}
\rput(6,-2.6){
\psline[linecolor=black](0,2)(0,0)
\psline[linecolor=blue](-2,0)(0,0)
\psline[linecolor=black](-2,0)(-2,0.5)
\psbezier[linecolor=green](-2,0.5)(-1.2,0.5)(0,1)(0,2)
\pscircle[linewidth=1pt,linestyle=dashed](0,1.8){0.2}
%\rput(-1,0.35){\makebox(0,0)[tr]{\hbox{\tiny$P$}}}
\rput(-1,1){\makebox(0,0)[cr]{\hbox{\small$\ell_{\lambda_0}/2$}}}
\rput(-0.75,-0.1){\makebox(0,0)[tc]{\hbox{\small$\ell_A/2$}}}
\rput(-2.1,0.25){\makebox(0,0)[cr]{\hbox{\small$h_4$}}}
\rput(0.1,0.8){\makebox(0,0)[cl]{\hbox{\small$h_1$}}}
}
\end{pspicture}
\caption{\small
A sphere with three holes and with a ``crown'' of $m$ decorated bordered cusps}
\label{fi:n-cusps}
\end{figure}
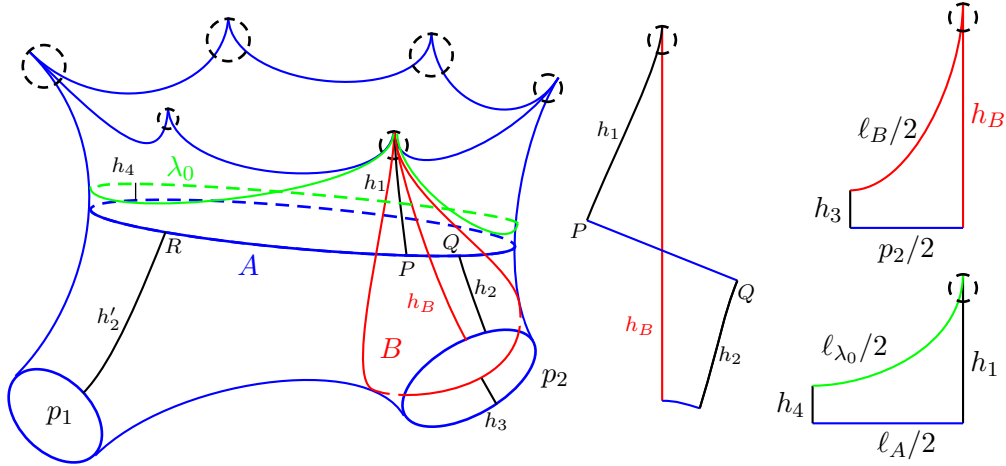

\begin{figure}[tb]
\begin{pspicture}(-6,-4)(6,5)
\newcommand{\PATTERNone}{%
{\psset{unit=1}
\pscircle[linecolor=lightgray, linewidth=1pt, fillcolor=lightgray, fillstyle=solid](0,0){0.87}
\rput{30}(0,0){\pscircle[linecolor=lightgray, linewidth=1pt, fillcolor=lightgray, fillstyle=solid](0,-1){0.5}\psarc[linewidth=1pt,linecolor=blue,linestyle=solid](0,-1){0.5}{150}{390}}
\rput{90}(0,0){\pscircle[linecolor=lightgray, linewidth=1pt, fillcolor=lightgray, fillstyle=solid](0,-1){0.5}\psarc[linewidth=1pt,linecolor=blue,linestyle=solid](0,-1){0.5}{150}{390}}
\rput{150}(0,0){\pscircle[linecolor=lightgray, linewidth=1pt, fillcolor=lightgray, fillstyle=solid](0,-1){0.5}\psarc[linewidth=1pt,linecolor=blue,linestyle=solid](0,-1){0.5}{150}{390}}
\rput{210}(0,0){\pscircle[linecolor=lightgray, linewidth=1pt, fillcolor=lightgray, fillstyle=solid](0,-1){0.5}\psarc[linewidth=1pt,linecolor=blue,linestyle=solid](0,-1){0.5}{150}{390}}
\rput{270}(0,0){\pscircle[linecolor=lightgray, linewidth=1pt, fillcolor=lightgray, fillstyle=solid](0,-1){0.5}\psarc[linewidth=1pt,linecolor=blue,linestyle=solid](0,-1){0.5}{150}{390}}
\rput{330}(0,0){\pscircle[linecolor=lightgray, linewidth=1pt, fillcolor=lightgray, fillstyle=solid](0,-1){0.5}\psarc[linewidth=1pt,linecolor=blue,linestyle=solid](0,-1){0.5}{150}{390}}
\rput(3.05,-4.4){\pscircle[linecolor=blue, linewidth=1pt, fillcolor=lightgray, fillstyle=solid](0,0){0.2}
\psarc[linewidth=1pt,linecolor=green,linestyle=dashed](-0.02,0.02){0.27}{-40}{140}
\psarc[linewidth=1pt,linecolor=green](0,0){0.35}{-120}{60}}
\rput(-2.4,-5.53){\pscircle[linecolor=blue, linewidth=1pt, fillcolor=lightgray, fillstyle=solid](0,0){0.2}
\psarc[linewidth=1pt,linecolor=red,linestyle=dashed](0,0){0.25}{-120}{60}
\psarc[linewidth=1pt,linecolor=red](0,0){0.35}{120}{300}}
\rput(0,0){\pscircle[linecolor=blue, linewidth=1pt](0,0){2}}
\put(-1.5,-5.7){\makebox(0,0)[cc]{\hbox{$P_1$}}}
\put(3.7,-4.8){\makebox(0,0)[cc]{\hbox{$P_2$}}}
\put(-2.05,0){\makebox(0,0)[cr]{\hbox{\small$\ell_A$}}}
\psline[linewidth=0.5pt,linestyle=dashed](0,0.87)(0,2.7)
% the crown
%\psbezier[linecolor=blue](-2.65,3.05)(-1.65,2.05)(1.81,1.534)(2.31,2.4)
}
}
\newcommand{\PATTERNtwo}{%
{\psset{unit=1}
\pscircle[linecolor=lightgray, linewidth=1pt, fillcolor=lightgray, fillstyle=solid](0,0){0.87}
\rput{30}(0,0){\pscircle[linecolor=lightgray, linewidth=1pt, fillcolor=lightgray, fillstyle=solid](0,-1){0.5}\psarc[linewidth=1pt,linecolor=blue,linestyle=solid](0,-1){0.5}{150}{390}
\psbezier[linecolor=blue, linestyle=dashed](-0.435,-0.75)(-0.675,-1.185)(-0.4,-1.6)(0,-1.6)}
\rput{90}(0,0){\pscircle[linecolor=lightgray, linewidth=1pt, fillcolor=lightgray, fillstyle=solid](0,-1){0.5}\psarc[linewidth=1pt,linecolor=blue,linestyle=solid](0,-1){0.5}{150}{390}}
\rput{150}(0,0){\pscircle[linecolor=lightgray, linewidth=1pt, fillcolor=lightgray, fillstyle=solid](0,-1){0.5}\psarc[linewidth=1pt,linecolor=blue,linestyle=solid](0,-1){0.5}{150}{390}}
\rput{210}(0,0){\pscircle[linecolor=lightgray, linewidth=1pt, fillcolor=lightgray, fillstyle=solid](0,-1){0.5}\psarc[linewidth=1pt,linecolor=blue,linestyle=solid](0,-1){0.5}{150}{390}
\psbezier[linecolor=blue, linestyle=dashed](0.435,-0.75)(0.675,-1.185)(0.4,-1.6)(0,-1.6)}
\rput{270}(0,0){\pscircle[linecolor=lightgray, linewidth=1pt, fillcolor=lightgray, fillstyle=solid](0,-1){0.5}\psarc[linewidth=1pt,linecolor=blue,linestyle=solid](0,-1){0.5}{150}{390}
\psbezier[linecolor=blue, linestyle=dashed](-0.435,-0.75)(-0.675,-1.185)(-0.4,-1.6)(0,-1.6)}
\rput{330}(0,0){\pscircle[linecolor=lightgray, linewidth=1pt, fillcolor=lightgray, fillstyle=solid](0,-1){0.5}\psarc[linewidth=1pt,linecolor=blue,linestyle=solid](0,-1){0.5}{150}{390}
\psbezier[linecolor=blue, linestyle=dashed](0.435,-0.75)(0.675,-1.185)(0.4,-1.6)(0,-1.6)}
%cycloide-almost
\rput{270}(0,0){\parametricplot[linecolor=green,linestyle=dashed,linewidth=1pt]{0}{360}{ 1.74 t cos mul -0.7 2 t mul cos mul add 0.6 0.5 t mul sin mul add 1.4 t sin mul -0.7 2 t mul sin mul add}}
\psarc[linewidth=1pt,linecolor=blue,linestyle=dashed](0,0){1.6}{-60}{120}
\psarc[linewidth=1pt,linecolor=blue,linestyle=dashed](0,0){1.6}{180}{240}
%
%cycloide-almost
%\rput{150}(0,0){\parametricplot[linecolor=red,linestyle=dashed,linewidth=1pt]{0}{360}{ 1.8 t cos mul -0.7 2 t mul cos mul add 0.55 0.5 t mul sin mul add 1.4 t sin mul -0.7 2 t mul sin mul add}}
%\rput{150}(0,0){\parametricplot[linecolor=red,linestyle=dashed,linewidth=1pt]{0}{360}{ 1.74 t cos mul -0.7 2 t mul cos mul add 0.55 0.5 t mul sin mul add 1.74 t sin mul -0.87 2 t mul sin mul add}}
%\rput{150}(0,0){\parametricplot[linecolor=red,linestyle=dashed,linewidth=1pt]{0}{360}{ 1.74 t cos mul -0.87 2 t mul cos mul add 0.6 0.5 t mul sin mul add 1.4 t sin mul -0.7 2 t mul sin mul add}}
%
\rput(3.05,-4.4){\pscircle[linecolor=blue, linewidth=1pt, fillcolor=lightgray, fillstyle=solid](0,0){0.2}
%\psarc[linewidth=1pt,linecolor=red](0,0){0.27}{-90}{90}
\psarc[linewidth=1pt,linecolor=green](0,0){0.35}{-120}{60}}
\rput(-2.4,-5.53){\pscircle[linecolor=white, linewidth=1pt, fillcolor=white, fillstyle=solid](0,0){0.2}
%\psarc[linewidth=1pt,linecolor=red,linestyle=dashed](0,0){0.25}{-120}{60}
%\psarc[linewidth=1pt,linecolor=red](0,0){0.35}{120}{300}
}
\rput(0,0){\pscircle[linecolor=blue, linewidth=1pt](0,0){2.7}}
%\put(-1.5,-5.7){\makebox(0,0)[cc]{\hbox{$P_1$}}}
\put(3.7,-4.8){\makebox(0,0)[cc]{\hbox{$P_2$}}}
\put(-2.75,0){\makebox(0,0)[cr]{\hbox{\small$\ell_A$}}}
% the crown
%\psbezier[linecolor=blue](-2.65,3.05)(-1.65,2.05)(1.81,1.534)(2.31,2.4)
}
}
\rput(-2.5,-0.35){\PATTERNone}
\rput{270}(-2.55,1.45){\parametricplot[linecolor=green,linestyle=dashed,linewidth=1pt]{0}{404}{ 0.016 t -0.000002 t 3 exp mul add mul  t cos mul  0.015  t  -0.000002 t 3 exp mul add mul  t sin mul }}
%
%\rput{270}(-2.55,1.45){\parametricplot[linecolor=green,linestyle=dashed,linewidth=1pt]{0}{450}{ 0.017 t -0.000002 t 3 exp mul add mul  t cos mul  0.016  t  -0.000002 t 3 exp mul add mul  t sin mul }}
%
\rput{270}(-2.55,1.45){\parametricplot[linecolor=green,linestyle=dashed,linewidth=1pt]{0}{404}{ 0.018 t -0.000002 t 3 exp mul add mul  t cos mul  0.017  t  -0.000002 t 3 exp mul add mul  t sin mul }}
\rput{270}(-2.55,1.45){\parametricplot[linecolor=red,linestyle=dashed,linewidth=1pt]{0}{330}{ 0.019 t -0.000002 t 3 exp mul add mul  t cos mul  0.018  t  -0.000002 t 3 exp mul add mul  t sin mul }}
%
%\rput{270}(-2.55,1.45){\parametricplot[linecolor=green,linestyle=dashed,linewidth=1pt]{0}{450}{ 0.020 t -0.000002 t 3 exp mul add mul  t cos mul  0.019  t  -0.000002 t 3 exp mul add mul  t sin mul }}
%
\rput{270}(-2.55,1.45){\parametricplot[linecolor=red,linestyle=dashed,linewidth=1pt]{0}{330}{ 0.021 t -0.000002 t 3 exp mul add mul  t cos mul  0.020  t  -0.000002 t 3 exp mul add mul  t sin mul }}
%\rput{270}(-2.55,1.45){\parametricplot[linecolor=red,linestyle=dashed,linewidth=1pt]{0}{450}{0 0.017 t mul add t cos mul 0 0.017  t mul add t sin mul }} %26 57 t add -2 exp mul add}
%
\rput{0}(-2.56,2.5){\psbezier[linecolor=red](-2.57,-5.215)(-1.7,-4.715)(0,-1.9)(0,-0.9)
\psbezier[linecolor=red](-2.23,-5.845)(-1.35,-5.345)(0,-2)(0,-0.9)}
\rput{0}(-2.56,2.5){\psbezier[linecolor=green](3.22,-4.085)(2.35,-3.585)(0,-2.8)(0,-0.9)
\psbezier[linecolor=green](2.88,-4.715)(2,-4.215)(0,-2.9)(0,-0.9)}
%\parametricplot[linecolor=red,linestyle=solid,linewidth=1pt]{0}{180}{2.05 t cos mul  0.45 t sin  mul}
%\parametricplot[linecolor=red,linestyle=dashed,linewidth=1pt]{180}{360}{2.05 t cos mul  0.45 t sin  mul}
%
%\rput{90}(0.1,-0.05){\parametricplot[linecolor=green,linestyle=solid,linewidth=1pt]{103}{180}{2.08 t cos mul 0.45 t sin mul 0.008 t mul add}
%\parametricplot[linecolor=green,linestyle=dashed,linewidth=1pt]{180}{360}{2.08 t cos mul  0.45 t sin  mul  0.008 t mul add}
%\parametricplot[linecolor=green,linestyle=solid,linewidth=1pt]{360}{455}{2.08 t cos mul 0.45 t sin mul 0.008 t mul add}}
%
%
%
%\psline[linewidth=2pt]{<->}(-0.5,0)(0.5,0)
%\rput{180}(2.5,0){\FLIP}
%
\put(0.3,-1.1){\makebox(0,0)[cc]{\hbox{\tcg{$\lambda_B$}}}}
\put(-4.2,-1.1){\makebox(0,0)[cc]{\hbox{\tcr{$\lambda_{AB}$}}}}
\put(-6.7,-2.3){\makebox(0,0)[cc]{\hbox{\tcr{$\lambda_{A^{-1}B}$}}}}
\put(-2.5,-2.3){\makebox(0,0)[cc]{\hbox{\tcg{$\lambda_{A^{-2}B}$}}}}
\rput(5,0.1){{\PATTERNtwo}
\rput{0}(0,0){\psbezier[linecolor=white](3.05,-4.67)(-3.5,-4.67)(-3.3,2)(-0.75,0.435)
\psbezier[linecolor=white](3.05,-4.13)(-3.15,-4.13)(-3,1.8)(-0.75,0.435)}
}
\rput{0}(5,2.5){\psbezier[linecolor=green](3.22,-4.085)(2.35,-3.585)(0,-2.8)(0,-0.9)
\psbezier[linecolor=green](2.88,-4.715)(2,-4.215)(0,-2.9)(0,-0.9)}
\put(7.8,-1.1){\makebox(0,0)[cc]{\hbox{\tcg{$\lambda_B$}}}}
%\put(6.5,-2.4){\makebox(0,0)[cc]{\hbox{\tcr{$\lambda'_B$}}}}
%\put(6.88,1.25){\makebox(0,0)[cc]{\hbox{\tcr{$\lambda'_0$}}}}
\put(3.4,3.88){\makebox(0,0)[cc]{\hbox{\tcg{$\lambda_0$}}}}
\put(4.5,2){\makebox(0,0)[cc]{\hbox{\small$\lambda^S_i$}}}
\put(5.5,3){\makebox(0,0)[cc]{\hbox{\small$\lambda^N_i$}}}
\put(4.35,3){\makebox(0,0)[cc]{\hbox{\small$i$}}}
\psline[linewidth=0.5pt,linecolor=black]{->}(5.65,3.15)(6.25,3.5)
\psline[linewidth=0.5pt,linecolor=black]{->}(4.35,1.85)(3.75,1.5)
\end{pspicture}
\caption{\small
A set of $\lambda$-lengths on a pair of pants with $m$ decorated bordered cusps (decorations are not shown); in the left side we present arcs contributing to $\tau_{B}$ and $\tau_{A^{-1}B}$, in the right side we present arcs $\lambda^N_i$ and $\lambda^S_i$ used in the construction of remaining $m-1$ canonical variables.}
\label{fi:tau-cusps}
\end{figure}
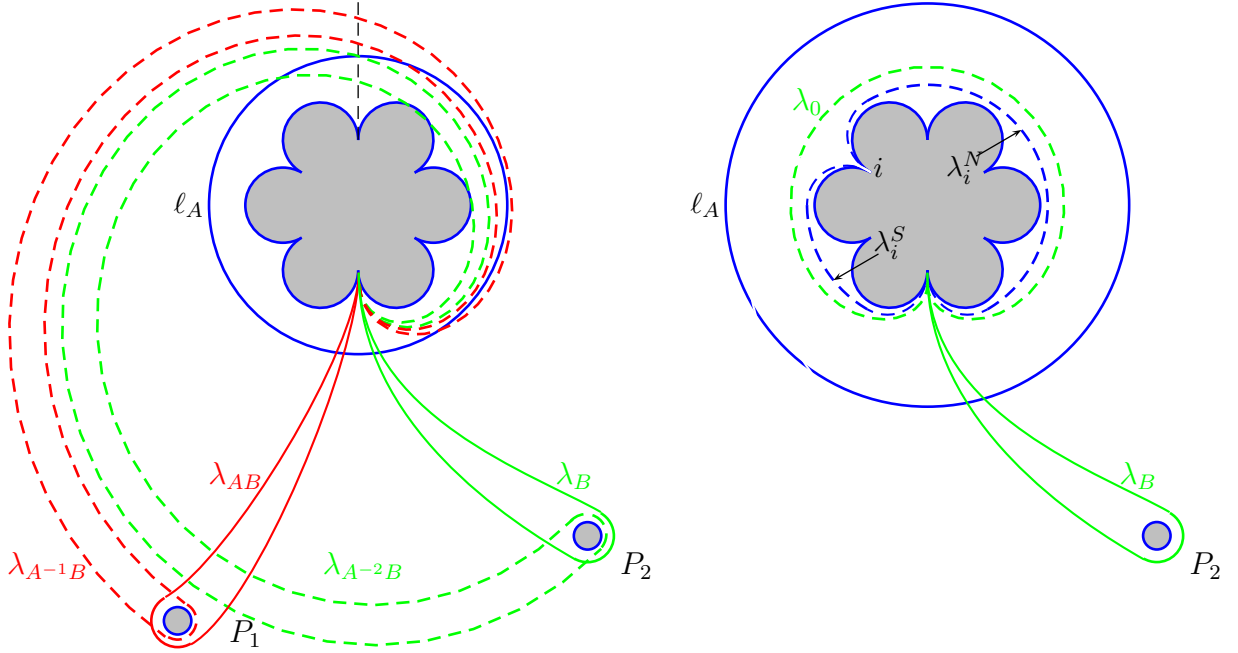

We first single out one bordered cusp assigning a label ``$0$'' to it and identify the twist coordinate. This coordinate must have a constant bracket with $\ell_A$ and has to depend only on $\ell_A$, $p_1$, $p_2$, $\lambda_B/\lambda_0$ (see Fig.~\ref{fi:n-cusps}) where $\lambda_B$ is the signed exponentiated half-length (a {\em $\lambda$-length}) of the part of a geodesic arc $B$ confined between two points of intersection with the decorating horocycle and $\lambda_0$ is the $\lambda$-length of the geodesic arc that separates the cusp ``crown'' region from the rest of the surface, as shown in the figure; the quotient $X_B:=\lambda_B/\lambda_0$ is a natural decoration-independent variable.

As in two previous cases, the choice of $\lambda_B$ is by no way unique: besides a possibility to choose $\lambda_B/\lambda_0$ at any bordered cusp out of $m$ cusps on the boundary of the hole, we have an infinite family of $\lambda_{A^kB}$, $k\in \mathbb Z$, and we present several terms of arcs from this sequence in the left part of Fig.~\ref{fi:tau-cusps}. We however show that, as for other twist coordinates, all the corresponding $\tau$-variables are related by constant shifts by integer multiples of $\ell_A/2$.

Let us define arcs $A^kB$, $k\in\mathbb Z$, as follows (see the left side of Fig.~\ref{fi:tau-cusps}): all these arcs start and terminate at the same boundary cusp; arcs $A^{2l}B$ go around hole $P_2$, arcs $A^{2l-1}B$ go around hole $P_1$ and on their way to the corresponding hole they intersect $l$ times the dashed vertical line in the figure (intersections are counted with signs: intersections in a clockwise direction come with the negative sign, those in a counterclockwise direction come with the positive sign).  We let $G_1=e^{p_1/2}+e^{-p_1/2}$ and $G_2=e^{p_2/2}+e^{-p_2/2}$ be the geodesic functions of perimeters of the corresponding holes.

In these notations, the bracket has the same form as in the two previous cases,
\be
\{\lambda_{A^kB},G_A\}=\lambda_{A^{k-1}B}-\lambda_{A^{k+1}B},
\label{AB-n}
\ee
where $\lambda_{A^{k-1}B}$ and $\lambda_{A^{k+1}B}$ are two solutions of the quadratic equations generated by a Markov triple,
\be
G_A\lambda_{A^{2l}B}\lambda_{A^{2l\pm 1}B}-\lambda_{A^{2l}B}^2-\lambda_{A^{2l\pm 1}B}^2-\lambda_{A^{2l\pm 1}B}\lambda_0G_{2}-\lambda_{A^{2l}B}\lambda_0G_{1}-\lambda^2_0=0,
\label{ABC-n}
\ee
Note that the variable $\lambda_{0}$ Poisson commute with $G_A$ and with all $\lambda_{A^{k}B}$.

We have the classical skein relations
\be
G_A\lambda_{A^{2l}B}=\lambda_{A^{2l+1}B}+\lambda_{A^{2l-1}B}+\lambda_0G_{2},\quad G_A\lambda_{A^{2l+1}B}=\lambda_{A^{2l+2}B}+\lambda_{A^{2l}B}+\lambda_0G_{1},\quad l\in \mathbb Z.
\label{sk-n}
\ee
Let us introduce the variable that is independent on the decoration of the cusp:
\be
X_{A^kB}:=\lambda_{A^{k}B}/\lambda_{0},\quad k\in\mathbb Z.
\label{X-n}
\ee
We can easily rewrite all relations \eqref{AB-n}--\eqref{sk-n} in terms of variables $X_{A^kB}$; for instance, skein relation \eqref{sk-n} becomes
\[
G_A X_{A^{2l}B}=X_{A^{2l+1}B}+X_{A^{2l-1}B}+G_{2},\quad G_AX_{A^{2l+1}B}=X_{A^{2l+2}B}+X_{A^{2l}B}+G_{1},\quad l\in \mathbb Z,
\]
and the Markov triple takes the form
\be
G_A X_{A^{2l}B}X_{A^{2l\pm 1}B}-X_{A^{2l}B}^2-X_{A^{2l\pm 1}B}^2-X_{A^{2l\pm 1}B}G_{2}-X_{A^{2l}B}G_{1}-1=0.
\label{Markov-X}
\ee
Since $X_{A^{k+1}B}$ and $X_{A^{k-1}B}$ are two solutions of the same quadratic equation (\ref{Markov-X}), their sum and product are
\be\label{sum-t}
\begin{array}{l}X_{A^{k+1}B}+X_{A^{k-1}B}=G_AX_{A^nB}-G_\beta,\phantom{\Bigm|} \\
X_{A^{k+1}B}X_{A^{k-1}B}=X_{A^kB}^2+X_{A^kB}G_\alpha+1,
\end{array} \quad \{\alpha,\beta\}=\left\{\begin{array}{ll}\{1,2\} & k\in 2{\mathbb Z}\\ \{2,1\} & k\in 2{\mathbb Z}+1\end{array}\right.
\ee

\begin{proposition}\label{lm-n1}
A decoration independent twist coordinate having a unit bracket with the half length $\ell_A/2$ of a closed geodesic $A$ for a sphere with three holes and with the boundary cusp is
\be
\tau_{A^{k}B}=\log\bigl((e^{\ell_A/2}-e^{-\ell_A/2})(X_{A^{k-1}B}-e^{-\ell_A/2}X_{A^{k}B})-e^{-\ell_A/2}G_{\alpha}-G_\beta\bigr),
\label{twistB-n}
\ee
where  $\{\alpha,\beta\}=\{1,2\}$ for even $k$ and  $\{\alpha,\beta\}=\{2,1\}$ for odd $k$.
All these twist coordinates are related by constant shifts:
\be
\tau_{A^kB}-\tau_{A^{k+1}B}=\ell_A/2.
\label{mod-n}
\ee
In particular, a Dehn twist along $A$ transforms $B$ into $A^2B$, and the corresponding modular transformation is $\tau_{B}\to \tau_{A^2B}=\tau_B-\ell_A$, so the twist coordinate (\ref{twistB-M4}) assumes values between $0$ and $\ell_A$ in a single copy of the modular space labelled by $A$.
\end{proposition}

The {\em proof} is a direct calculation. Checking that $\{\tau_B,\ell_A/2\}=1$ is straightforward using (\ref{sum-t}) and (\ref{sk-n}):
\begin{align*}
\{e^{\tau_B},G_A\}&=(e^{\ell_A/2}-e^{-\ell_A/2})\bigl(X_{A^{-2}B} - X_B-e^{-\ell_A/2}(X_{A^{-1}B} - X_{AB})\bigr)\\
&=(e^{\ell_A/2}-e^{-\ell_A/2})\bigl(G_AX_{A^{-1}B} - G_1- 2X_B-e^{-\ell_A/2}(2X_{A^{-1}B} - G_AX_{B}+G_2)\bigr)\\
&=(e^{\ell_A/2}-e^{-\ell_A/2})\bigl((e^{\ell_A/2}-e^{-\ell_A/2})(X_{A^{-1}B} - e^{-\ell_A/2}X_B)- G_1-e^{-\ell_A/2}G_2\bigr)\\
&=(e^{\ell_A/2}-e^{-\ell_A/2})e^{\tau_B}.
\end{align*}
The proof of (\ref{mod-n}) we leave to the reader as an exercise.

Choosing another sign of  $\ell_A$, we obtain
\[
\tau'_{A^{k}B}=\log\bigl((e^{-\ell_A/2}-e^{\ell_A/2})(X_{A^{k-1}B}-e^{\ell_A/2}X_{A^{k}B})-e^{\ell_A/2}G_{\alpha}-G_\beta\bigr),
\]
Then, for the sum $\tau_{A^kB}+\tau'_{A^kB}$, we obtain that all terms containing $X$-terms combine into a Markov triple, and finally,
\be
\tau_{A^kB}+\tau'_{A^kB}=\log\bigl(G_{A}G_1G_2+G^2_A +G_1^2+G^2_2-4\bigr).
\ee
This immediately implies the analogue of Lemma~\ref{lm:twist}.

\begin{lemma}\label{lm:twist-X}
The canonical twist coordinate that change its sign under changing the orientation ($\ell_A\to -\ell_A$) for the perimeter of a hole containing a bordered cusp is
\be\label{tau-X}
\widehat \tau_{A^kB}:=\log \Bigl[ \frac{(e^{\ell_A/2}-e^{-\ell_A/2})(X_{A^{k-1}B} - e^{-\ell_A/2}X_{A^kB}) -e^{-\ell_A/2}G_{\alpha}-G_\beta}{(G_AG_1G_2+G_A^2+G_1^2+G_2^2-4)^{1/2}} \Bigr]
\ee
Geometrically, $\widehat \tau_B$ is the signed geodesic length along the geodesic $\gamma_A$ between endpoints of perpendiculars to $\gamma_A$ from the selected (``zeroth'') cusp and the hole $P_2$, see Fig.~\ref{fi:n-cusps}.
\end{lemma}

To prove the geometrical component of this statement, we again use identities from hyperbolic geometry; note that all identities in a right-angled pentagon with one cusp can be obtained from those of right-angled hexagon in the limit where the length of one of its sides goes to zero (then lengths of two adjacent sides tend to infinity, and only the ratio of their $\lambda$-lengths enters relations). In the geometry of Figure~\ref{fi:n-cusps}, we have the following identities: for a pentagon, we have that
\be\label{CUSP1}
e^{h_B-h_1}=\sinh(h_2)\cosh(|PQ|)+\cosh(h_2),
\ee
and in two quadrangles we have
\be\label{CUSP2}
e^{\ell_B/2-h_B}=\sinh(p_2)\quad \hbox{and}\quad e^{\ell_{\lambda_0}/2-h_1}=\sinh(\ell_{A}/2).
\ee
Combining these three relations and taking into account that $X_B=e^{\ell_B/2-\ell_{\lambda_0}/2}$, we obtain
\be
X_B \sinh(\ell_A/2)=\sinh(p_2)\bigl( \sinh(h_2)\cosh(|PQ|)+\cosh(h_2)\bigr).
\ee
A similar relation holds for $X_{A^{-1}B}$ except that we have to replace $h_2$ by $h'_2$ (see Fig.~\ref{fi:n-cusps}) and $|PQ|$ by $|PR|$: since $|RQ|=\ell_A/2$,
we have
\be
X_{A^{-1}B} \sinh(\ell_A/2)=\sinh(p_1)\bigl( \sinh(h'_2)\cosh(|PQ|+\ell_A/2)+\cosh(h'_2)\bigr).
\ee
Here $h_2$ and $h'_2$ are lengths of sides of a right-angled hexagon with other sides $\ell_A/2$, $p_1/2$, and $p_2/2$. The hyperbolic sine theorem indicates that
$$
\sinh(p_1) \sinh(h'_2) = \sinh(p_2) \sinh(h_2).
$$
Finally, we have the standard relations for the perpendiculars $h_2$ and $h'_2$:
\begin{align*}
\cosh(h_2)\sinh(p_2/2)\sinh(\ell_A/2)&=\cosh(p_1/2)+\cosh(p_2/2)\cosh(\ell_A/2),\\
\cosh(h'_2)\sinh(p_1/2)\sinh(\ell_A/2)&=\cosh(p_2/2)+\cosh(p_1/2)\cosh(\ell_A/2),
\end{align*}
using which, after a short algebra, we obtain
$$
\sinh(\ell_A/2)(X_{A^{-1}B} - e^{-\ell_A/2}X_{B})-e^{-\ell_A/2}\cosh(p_2/2)-\cosh(p_1/2)=2\sinh(\ell_A/2) \sinh(p_2) \sinh(h_2)e^{|PQ|}
$$
and it only remains to note that $4\sinh^2(\ell_A/2) \sinh^2(p_2) \sinh^2(h_2)=G_AG_1G_2+G_A^2+G_1^2+G_2^2-4$.

We then have an extension of Theorem~\ref{th:twist}.

\begin{theorem}\label{th:twist-X}
Given any pair-of-pant decomposition of a Riemann surface $\Sigma_{g,s,n}$ and defining the canonical twist variables $\widehat \tau_B$ by formulas (\ref{tau-can}), (\ref{tau-torus}), and (\ref{tau-X}) we have that all these canonical twist variables Poisson commute.
\end{theorem}

\begin{remark}\label{rm:twist-X}
Although we omit the proof of commutativity of canonical twist variables of type (\ref{tau-X}) with other types of canonical twist variables and with canonical twist variables of the same type for the brevity of exposition, note that either we can perform this proof directly, using the same template as in the proof of Theorem~\ref{th:twist}, or, in a more advantageous way, we can use the fact that all surfaces with bordered cusps can be obtained as reductions of surfaces with holes using ``chewing gum'' moves from \cite{ChMaz2}, \cite{CMR2} and, correspondingly, there should be a way to obtain canonical twist variables of type (\ref{tau-X}) from canonical twist variables of types (\ref{tau-can}) and (\ref{tau-torus}).
\end{remark}

\subsection{Local variables on $\Sigma_{0,2,m}$ with $m\ge 2$}
After identifying the twist coordinate $\tau_B$ (\ref{tau-X}) we are going to construct the remaining $m-1$ coordinates on the selected boundary component. We construct them in the following way. Fixing the zeroth cusp which we use to construct  $\tau_B$, we enumerate all other cusps in clockwise direction counted from the selected zeroth cusp. Let then $\lambda^S_i$ and $\lambda^N_i$, $i=1,\dots,m-1$, denote the lambda lengths of arcs starting at the zeroth cusp and terminating at the $i$th cusp and going in the respective clockwise and counterclockwise directions (see the right half of Fig.~\ref{fi:tau-cusps}).

Note first that all $\lambda^N_i$ and $\lambda^S_i$ for a given hole commute with $X_B$ and $\ell_A$ constructed for this hole and obviously commute with similar variables constructed for other holes and with all canonical twist--length coordinates, so the algebra of these variables separates completely from the ``big'' algebra of Fenchel--Nielsen coordinates.

We introduce $m-1$ decoration-independent variables
\be\label{x-i}
x_i:=\lambda^N_i/\lambda^S_i.
\ee
It is an easy exercise using skein relations and Goldman brackets to find the Poisson relations between $x_i$:
\be
\{x_i,x_j\}=x_i^2-x_ix_j,\quad 1\le i<j\le m-1.
\ee
We next consider the combinations of these variables:
\be\label{r-i}
r_1:=x_1,\quad r_i:=x_i-x_{i-1},\ 2\le i\le m-1.
\ee
We then have the lemma
\begin{lemma}
The decoration-independent variables $r_i$ (\ref{r-i}) have homogeneous Poisson relations
\be\label{r-r}
\{r_i,r_j\}=r_ir_j,\quad 1\le i<j\le m-1;
\ee
this algebra is nondegenerate for odd $m$ and it has a unique Casimir element
\be\label{r-C}
C=\frac{r_1r_3\cdots r_{m-3} r_{m-1}}{r_2r_4\cdots r_{m-2}}
\ee
for even $m$.
\end{lemma}
Note that the variables $r_i$ take all values in $\mathbb R_+^{m-1}$, and for any given set of values of $r_i$ we have a unique configuration of cusps on the corresponding boundary component.

\begin{remark}
Since the Fenchel--Nielsen Poisson structure was derived solely using the Goldman bracket, and the Goldman bracket, in turn, follows from the Fock bracket on the set of shear coordinates, the constructed Fenchel--Nielsen Poisson structure is derived from the Fock bracket. Since the latter is manifestly mapping-class group invariant, so is the Fenchel--Nielsen bracket. Moreover, since we restrict consideration to decoration-independent variables, all length and twist coordinates and all $r$-variables are functions of shear coordinates $Z_\alpha$ and coefficients $\omega_j$ only. The only Casimirs of subalgebra of $Z_\alpha$ are those for holes containing even numbers $n_i$ of bordered cusps; for each of $n_i$ {\em windows} (in the Kauffman--Penner terminology \cite{KP}) constituting together the boundary of the corresponding hole, we consider sums of shear coordinates $\sum_\alpha Z_{\alpha}^{(j)}$ of edges incident to the $j$th window, and the logarithm of the Casimir element (\ref{r-C}) is then given by the alternating sum
$$
\sum_{j=1}^{n_i} (-1)^j \Bigl[{\textstyle\sum_\alpha }Z_{\alpha}^{(j)}\Bigr].
$$
\end{remark}

\subsection{Special cases}\label{ss:special}
To complete our accounting of Fenchel--Nielsen coordinates, it remains to consider two special cases: $\Sigma_{0,2,n_1+n_2}$---a cylinder with $n_1>0$ and $n_2>0$ bordered cusps on its two boundaries and $\Sigma_{0,1,n}$---a disc with $n\ge 3$ bordered cusps on the boundary.

\subsubsection{The twist coordinate for $\Sigma_{0,2,n_1+n_2}$}
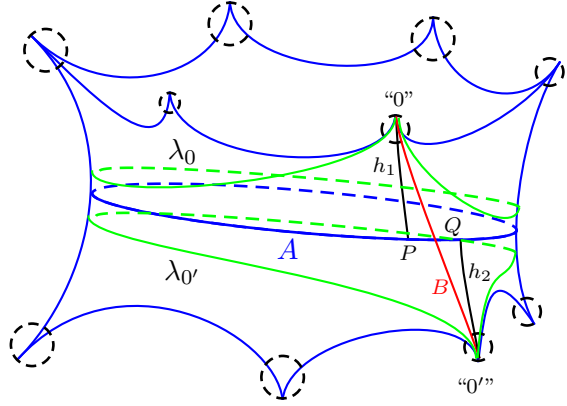
\begin{figure}[tb]
\begin{pspicture}(-3,-3)(4,3){
\newcommand{\PATTERN}{%
{\psset{unit=1}
%\rput{-45}(-3,-2.5){\psellipse[linecolor=blue, linewidth=1pt](0,0)(.75,.5)}\rput{30}(2.5,-2){\psellipse[linecolor=blue, linewidth=1pt](0,0)(1,.5)}
\rput{-5}(0.29,0){\psellipse[linecolor=blue,  linewidth=1pt](0,0)(2.85,0.3)\psframe[linecolor=white, fillstyle=solid, fillcolor=white](-2.9,0)(2.9,.6)\psellipse[linecolor=blue, linestyle=dashed,  linewidth=1pt](0,0)(2.85,.3)}
%\rput{45}(-3,2.7){\psellipse[linecolor=blue,  linewidth=1pt](0,0)(.5,.3)}
%\rput{-30}(3,2){\psellipse[linecolor=blue,  linewidth=1pt](0,0)(.8,.3)}
%\psbezier[linecolor=blue](-2.47,-3.03)(-1.47,-2.03)(1.134,-1.634)(1.634,-2.5)
\psbezier[linecolor=blue](-3.53,-1.97)(-2.56,-1)(-2,1)(-3.35,2.35)
\psbezier[linecolor=blue](3.35,-1.5)(2.85,-0.6)(3.19,1.2)(3.69,2)
\psbezier[linecolor=blue](3.35,-1.5)(2.75,-0.6)(2.7,-1.2)(2.6,-2)
\psbezier[linecolor=blue](0,-2.5)(0.1,-1.6)(2.5,-1.2)(2.6,-2)
\psbezier[linecolor=blue](0,-2.5)(-0.1,-1.6)(-2.2,-0.5)(-3.53,-1.97)
% the crown
\psbezier[linecolor=blue](-1.5,1.6)(-1.56,0.6)(-2.35,1.35)(-3.35,2.35)
\psbezier[linecolor=blue](1.5,1.3)(1.56,0.4)(3.19,1.2)(3.69,2)
\psbezier[linecolor=blue](1.5,1.3)(1.46,0.4)(-1.46,0.6)(-1.5,1.6)
\psbezier[linecolor=blue](-0.7,2.8)(-0.76,1.8)(-2.35,1.35)(-3.35,2.35)
\psbezier[linecolor=blue](2,2.6)(2.06,1.6)(3.19,1.2)(3.69,2)
\psbezier[linecolor=blue](2,2.6)(2,1.6)(-0.66,1.8)(-0.7,2.8)
% decorations: upper part
\pscircle[linewidth=1pt,linestyle=dashed](-3.15,2.15){0.3}
\pscircle[linewidth=1pt,linestyle=dashed](-1.5,1.45){0.15}
\pscircle[linewidth=1pt,linestyle=dashed](-0.7,2.5){0.3}
\pscircle[linewidth=1pt,linestyle=dashed](1.5,1.1){0.2}
\pscircle[linewidth=1pt,linestyle=dashed](2,2.3){0.3}
\pscircle[linewidth=1pt,linestyle=dashed](3.55,1.86){0.2}
% decorations: lower part
\pscircle[linewidth=1pt,linestyle=dashed](-3.33,-1.77){0.3}
\pscircle[linewidth=1pt,linestyle=dashed](0,-2.2){0.3}
\pscircle[linewidth=1pt,linestyle=dashed](2.6,-1.8){0.2}
\pscircle[linewidth=1pt,linestyle=dashed](3.25,-1.33){0.2}
%\psbezier[linecolor=blue](-2.65,3.05)(-1.65,2.05)(1.81,1.534)(2.31,2.4)
}
}
\rput(0,0){\rput(0,0){\PATTERN}
\rput{0}(0,0){%\psbezier[linecolor=red](1.44,1.27)(1.42,0.4)(.75,-1.95)(1.15,-2.15)
%\psbezier[linecolor=red](1.44,1.27)(1.46,0.4)(3,-0.4)(3.1,-1)
%\psbezier[linecolor=red](1.15,-2.15)(1.55,-2.35)(3.3,-2.2)(3.1,-1)
%\rput(1.44,-2.23){\makebox(0,0)[cc]{\hbox{\tcw{\tiny$\bullet$}}}}
%\rput(3.1,-1.25){\makebox(0,0)[cc]{\hbox{\tcw{\tiny$\bullet$}}}}
% perpendiculars
\rput(0,0){\psbezier[linecolor=black](1.44,1.27)(1.47,0.7)(1.47,0.7)(1.6,-0.35)
\psbezier[linecolor=red](1.44,1.27)(1.5,0.7)(2.3,-1.2)(2.55,-1.95)
\psbezier[linecolor=black](2.3,-0.37)(2.3,-0.7)(2.4,-1.2)(2.55,-1.95)
%\psbezier[linecolor=black](-1.6,-0.05)(-1.9,-0.8)(-2.3,-1.8)(-2.7,-2.2)
%\psline[linecolor=black](2.6,-2)(2.8,-2.35)
%\psline[linecolor=black, linewidth=0.5pt, linestyle=dashed](-2,0.37)(-2,0.6)
}
}
\rput{-5}(0.29,0.3){\parametricplot[linecolor=green,linestyle=dashed,linewidth=1pt]{0}{180}{2.8 t cos mul  0.15 t sin  mul}
\psbezier[linecolor=green](1.1,1.05)(1.2,0.3)(2.75,-0.5)(2.8,0)
\psbezier[linecolor=green](1.05,1.05)(1.1,0.2)(-2.95,-0.5)(-2.9,0)
}
\rput{-5}(0.24,-0.3){\parametricplot[linecolor=green,linestyle=dashed,linewidth=1pt]{0}{180}{2.8 t cos mul  0.15 t sin  mul}
\psbezier[linecolor=green](2.42,-1.44)(2.4,0)(2.75,-0.5)(2.8,0)
\psbezier[linecolor=green](2.4,-1.45)(2.3,-0.4)(-2.95,-0.5)(-2.9,0)
}
\rput(-1,0){
%\rput(2.4,-1.6){\makebox(0,0)[cc]{\hbox{\tcr{\small$B$}}}}
%\put(2.4,-1.6){\pscircle[linecolor=green]{.2}}
\rput(-0.4,0.8){\makebox(0,0)[cc]{\hbox{{\small$\lambda_0$}}}}
\rput(-0.4,-0.8){\makebox(0,0)[cc]{\hbox{{\small$\lambda_{0'}$}}}}
%\put(2.3,1.7){\pscircle[linecolor=red]{.2}}
\rput(1,-0.5){\makebox(0,0)[cc]{\hbox{\tcb{\small$A$}}}}
%\rput(-2,-2.5){\makebox(0,0)[cc]{\hbox{\small$p_1$}}}
%\rput(4.4,-2){\makebox(0,0)[cl]{\hbox{\small$p_2$}}}
\rput(2.6,-0.45){\makebox(0,0)[tc]{\hbox{\tiny$P$}}}
\rput(3.2,-0.33){\makebox(0,0)[bc]{\hbox{\tiny$Q$}}}
%\rput(-0.5,-0.33){\makebox(0,0)[bc]{\hbox{\tiny$R$}}}
\rput(2.45,0.6){\makebox(0,0)[cr]{\hbox{\tiny$h_1$}}}
\rput(3.4,-0.8){\makebox(0,0)[cl]{\hbox{\tiny$h_2$}}}
%\rput(-1.2,-1.2){\makebox(0,0)[cr]{\hbox{\tiny$h'_2$}}}
\rput(2.7,1.5){\makebox(0,0)[cr]{\hbox{\tiny``$0$''}}}
\rput(3.5,-2.2){\makebox(0,0)[tc]{\hbox{\tiny``$0'$''}}}
%\rput(-1,0.7){\makebox(0,0)[br]{\hbox{\tiny$h_4$}}}
\rput(3.15,-1){\makebox(0,0)[cr]{\tcr{\hbox{\tiny$B$}}}}
%\put(0,-0.5){\pscircle[linecolor=blue]{.2}}
}
}
}
\end{pspicture}
\caption{\small
A cylinder with $n_1>0$ and $n_2>0$ decorated bordered cusps on its two boundaries.}
\label{fi:cylinder}
\end{figure}

In the case of cylinder containing bordered cusps on both its boundaries, we have a unique closed geodesics $A$, we select two cusps marked ``$0$'' and ``$0'$'' on two different boundary components, and consider a $\mathbb Z$-labeled set of arcs $\mathfrak a_{A^kB}$ joining these two cusps. We let $\lambda_{A^kB}$ denote $\lambda$-lengths of arcs from this set. As above, we also have two special arcs, $\mathfrak a_0$ and $\mathfrak a_{0'}$ (see Fig.~\ref{fi:cylinder}) and we let $\lambda_0$ and $\lambda_{0'}$ denote their $\lambda$-lengths. We now construct the twist coordinate dual to $\ell_A$ and determine its geometric origin.

We have a set of Poisson and skein relations on $G_A$ and $\lambda$-lengths:
\be
\{G_A,\lambda_{A^kB}\}=\frac12 (\lambda_{A^{k-1}B}-\lambda_{A^{k+1}B}),
\ee
\be
G_A\lambda_{A^kB}=\lambda_{A^{k-1}B}+\lambda_{A^{k+1}B};\qquad \lambda_{A^{k-1}B}\lambda_{A^{k+1}B}=\lambda^2_{A^kB}+\lambda_0\lambda_{0'}.
\ee
We now introduce decoration-independent variables
\be\label{X-cyl}
X_{A^kB}:=\frac{\lambda_{A^kB}}{\sqrt{\lambda_0\lambda_{0'}}}=e^{\ell_{A^kB}/2-\ell_{\lambda_0}/4-\ell_{\lambda_{0'}}/4},
\ee
in terms of which the above relations read
\begin{align*}
&\{G_A,X_{A^kB}\}=\frac12 (X_{A^{k-1}B}-X_{A^{k+1}B})\\
&G_AX_{A^kB}=X_{A^{k-1}B}+X_{A^{k+1}B}, \ \hbox{and}\ X_{A^{k-1}B}X_{A^{k+1}B}=X^2_{A^kB}+1,
\end{align*}
which implies the Markov relation
\be\label{Markov-cyl}
G_AX_{A^kB}X_{A^{k-1}B}-X^2_{A^kB}-X^2_{A^{k-1}B}-1=0,\quad k\in\mathbb Z.
\ee

As before, a twist coordinate having a unit bracket with $\ell_A$ is
$$
\tau_{A^kB}:=\log\bigl(X_{A^{k-1}B}-e^{-\ell_A/2}X_{A^kB} \bigr),
$$
and it happens that it is also a canonical twist coordinate; indeed, taking $\ell_A\to -\ell_A$, we obtain
$$
\tau'_{A^kB}:=\log\bigl(X_{A^{k-1}B}-e^{\ell_A/2}X_{A^kB} \bigr),
$$
and
\begin{align}
e^{\tau_{A^kB}}e^{\tau'_{A^kB}}&=\bigl(X_{A^{k-1}B}-e^{-\ell_A/2}X_{A^kB} \bigr)\bigl(X_{A^{k-1}B}-e^{\ell_A/2}X_{A^kB} \bigr)\nonumber\\
&=-G_AX_{A^kB}X_{A^{k-1}B}+X^2_{A^kB}+X^2_{A^{k-1}B}=-1.\label{CC}
\end{align}
So, we have the lemma
\begin{lemma}\label{lm:twist-cyl}
For a cylinder with bordered cusps ``$0$'' and ``$0'$'' on its two boundary components,
the canonical twist coordinate that have the unit bracket with $\ell_A$ and change its sign under changing the orientation ($\ell_A\to -\ell_A$) for the diameter of a cylinder is
\be\label{tau-cyl}
\widehat \tau_{A^kB}:=\log\bigl(X_{A^{k-1}B}-e^{-\ell_A/2}X_{A^kB} \bigr),
\ee
where $X_{A^kB}$ are decoration-independent variables (\ref{X-cyl}).
Geometrically, $\widehat \tau_B$ is half of the signed geodesic length $|PQ|$ along the geodesic $\gamma_A$ between endpoints of perpendiculars to $\gamma_A$ from the selected (``zeroth'') cusps, see Fig.~\ref{fi:cylinder}.
\end{lemma}

Let us prove the geometrical part of the lemma statement. Taking the limit $h_2\to \infty$ in (\ref{CUSP1}) with accounting for decorations of cusps, we obtain
$$
e^{\ell_{B}-h_1-h_2}=\frac12(\cosh(|PQ|)+1)=\cosh^2(|PQ/2|).
$$
Together with the standard relations (\ref{CUSP2}) in quadrangles, which read
$$
e^{\ell_{\lambda_0}/2-h_1}=e^{\ell_{\lambda_{0'}}/2-h_2}=\sinh(\ell_{A}/2),
$$
we obtain
$$
X_B^2=e^{\ell_{B}-\ell_{\lambda_0}/2-\ell_{\lambda_{0'}}/2}=\frac{\cosh^2(|PQ/2|)}{\sinh^2(\ell_{A}/2)},
$$
or
\be
(e^{\ell_A/2}-e^{-\ell_A/2})X_B=e^{|PQ|/2}+e^{-|PQ|/2}.
\ee
We now use (\ref{CC}) that indicates that $e^{-\tau_B}=-\bigl(X_{A^{-1}B}-e^{\ell_A/2}X_{B} \bigr)$ to obtain
$$
e^{\tau_B}+e^{-\tau_B}=X_B(e^{\ell_A/2}-e^{-\ell_A/2})=e^{|PQ|/2}+e^{-|PQ|/2},
$$
which immediately implies the lemma statement.

\subsubsection{Fenchel--Nielsen coordinates for $\Sigma_{0,1,n}$}
The last remaining case is a disc $\Sigma_{0,1,n}$ with $n\ge3$ boundary cusps, which we enumerate from $0$ to $n-1$ moving counterclockwise. The cusps labeled ``$0$'', ``$1$'', and ``$2$'' play a special role in the construction below. We now have unique arcs $\mathfrak a_{i,j}$ joining cusps with labels $i$ and $j$, and since arcs are not directed, we can always assume that $i<j$. Consider arcs $\mathfrak a_{0,i}$ and $\mathfrak a_{1,i}$ with the corresponding lambda lengths $\lambda_{0,i}$ and $\lambda_{1,i}$. We have the Poisson relations
\begin{align*}
&\{\lambda_{0,i},\lambda_{0,j}\}=\frac14 \lambda_{0,i}\lambda_{0,j},\quad \{\lambda_{1,i},\lambda_{1,j}\}=\frac14 \lambda_{1,i}\lambda_{1,j},\quad  i<j\\
&\{\lambda_{1,i},\lambda_{0,j}\}=0,\quad \{\lambda_{0,i},\lambda_{1,j}\}=\frac12\bigl(\lambda_{1,i}\lambda_{0,j}-\lambda_{j,i}\lambda_{0,1} \bigr),\quad  i<j,
\end{align*}
and the skein relation
$$
\lambda_{0,i}\lambda_{1,j}=\lambda_{1,i}\lambda_{0,j}+\lambda_{j,i}\lambda_{0,1},
$$
using which we can express the last bracket in the form
$$
 \{\lambda_{0,i},\lambda_{1,j}\}=\lambda_{1,i}\lambda_{0,j}-\frac12\lambda_{0,i}\lambda_{1,j},\quad  i<j,
$$
We introduce (still not completely decoration-independent) combinations
\be
x_j:=\lambda_{1,j}/\lambda_{0,j},\quad j=2,\dots,n-1.
\ee
The brackets between $x_j$ read
\be
\{x_i,x_j\}=x_ix_j-x_i^2,\quad 2\le i<j\le n-1,
\ee
that is, they have exactly the same form as brackets in the subalgebra of  ``local''  variables for a selected boundary component.
Then, taking
\be
r_i=\{x_2,\ i=2;\ x_i-x_{i-1}, \ 3\le i\le n-1\},
\ee
we obtain homogeneous Poisson relations $\{r_i,r_j\}=r_ir_j$ for $2\le i<j\le n-1$ and, finally, we introduce completely decoration-independent variables
\be\label{r-disc}
\widehat r_i:=r_i/r_2,\quad 3\le i\le n-1,
\ee
for which we have the proposition.
\begin{proposition}
in the disc with $n\ge 4$ bordered cusps we have $n-3$  Fenchel--Nielsen coordinates $\widehat r_i$ (\ref{r-disc}) that are
decoration independent and have homogeneous Poisson relations
\be\label{r-r1}
\{r_i,r_j\}=r_ir_j,\quad 3\le i<j\le n-1;
\ee
this algebra is nondegenerate for odd $n$ and it has a unique Casimir element
\be\label{r-C1}
C=\frac{r_3\cdots r_{n-3} r_{n-1}}{r_4\cdots r_{n-2}}
\ee
for even $n$.
\end{proposition}

\subsection{Mirzakhani's volumes of ${\mathcal M}_{g,s,n}$}
We can now define the Fenchel--Nielsen sympelctic form as an inverse to the Poisson brackets: we have $3g-3+s+s_h$ pairs of canonical length-twist coordinates, and for every hole with $n_i>0$ bordered cusps on its boundary we have $(n_i-1)$-dimensional Poisson algebra of $r$-variables associated with this hole. These algebras of $r$-variables are nontrivial provided $n_i\ge 3$. The corresponding Fenchel--Nielsen (or Weil--Petersson) 2-form then reads
$$
\omega_{\text{WP}}=\sum_{k=1}^{3g-3+s+s_h} d\ell_A\wedge d\tau_B + \sum_{i=1}^{s_h} \omega_j,
$$
where we let $w_j$ denote ``local'' 2-forms obtained by inverting Poisson structures (\ref{r-r}) (in the orthogonal complements to the Casimir elements (\ref{r-C}) if $n_i$ is even). The volume element is the relevant power of $\omega_{\text{WP}}$.

The seminal result of Mirzakhani \cite{Mir06} enables one to evaluate volumes of {\em moduli spaces} ${\mathcal M}_{g,s}$ (for $s>0$) obtained by factoring Teichm\"uller spaces ${\mathfrak T}_{g,s}$ by the action of the mapping-class group. Mirzakhani's theorem states that, provided we fixed perimeters $p_i$ of all holes, these volumes are finite and are polynomials in $p_i$. This assertion was based on McShane's (\cite{McS91},\cite{McS98}) and Mirzakhani's (\cite{Mir06}) identities for lengths of simple (i.e., without self-intersections) closed geodesics in $\Sigma_{g,s}$, recurrence relations derived in \cite{Mir07} was shown by Mulase, Safnuk \cite{MS06} and Do and Norbury \cite{DoN06} to satisfy Virasoro algebra relations, which eventually resulted in the construction of topological recursion model by Eynard and Orantin \cite{EO07} governing the generating function for these volumes; this generating function was immediately identified with a KdV hierarchy $\tau$-function.

Note that the Fenchel--Nielsen coordinates were instrumental in Mirzakhani's  derivation of volume formulas;  an advantage of this set of coordinates is that the integration over the twist coordinate goes from $0$ to $\ell_A$ in one copy of the corresponding moduli space, so the net effect of this integration is multiplication by $\ell_A$. For example, for a torus with a hole of perimeter $p$, we have the McShane--Mirzakhani identity \cite{McS91}, \cite{Mir06}, which, after differentiating w.r.t. $p$ takes a convenient form $\sum_\gamma \Bigl[\frac{1}{e^{\ell_\gamma+p/2}+1}+\frac{1}{e^{\ell_\gamma-p/2}+1}\Bigr]=1$, where the sum ranges all simple closed curves, each of which can be identified with an $A$-cycle in the corresponding copy of the moduli space, and the volume of the corresponding moduli space, on the one hand, is the integral of just the constant function $1$ and, on the other hand, is the integral from zero to infinity (the range of $\ell_\gamma$) of the function in the left-hand side multiplied by $\ell_\gamma$ due to the integration over the twist coordinate; the answer then reads
$$
\hbox{Vol}_{1,1}\simeq\int_0^\infty \Bigl[\frac{1}{e^{\ell_\gamma+p/2}+1}+\frac{1}{e^{\ell_\gamma-p/2}+1}\Bigr] \ell_\gamma d\ell_\gamma=\frac{\pi^2}{6}+\frac{p^2}{8},
$$
and amending for volumes of discrete automorphism groups for two terms in the right-hand side ($2$ and $6$ in this case), we finally obtain that $\hbox{Vol}_{1,1}=\frac{\pi^2}{12}+\frac{p^2}{48}$.

Since the mapping-class group is not affected by the presence of bordered cusps, and all closed curves that are perimeters of holes are preserved by the action of this group, we can fix all these perimeters $p_i$, $i=1,\dots,s$, exactly as in Mirzakhani's accounting. We now have twist coordinates $\tau_{B_i}$ for $s_h$ holes that contain bordered cusps. Note that normalizations of $\tau_{B_i}$ is not mapping-class group invariant depending on the choice of internal $A$-cycles, so only the differentials $d\tau_{B_i}$ are mapping-class group invariant objects and adding integrations over these twist coordinates results in multiplying the corresponding Mirzakhani's volume by $\prod_{j=1}^{s_h} p_j$. Adding one more cusp to a boundary component results in adding one more Casimir and does not affect the volume, but if we have more than two cusps on a component, we add an integration over the whole plane $\mathbb R^2$ making the corresponding volumes infinite.

\section{Concluding remarks and perspectives}
\setcounter{equation}{0}
We have shown that the Poisson bracket introduced by V.V.Fock on the set of shear coordinates induces, via the Goldman bracket, the Fenchel--Nielsen Poisson bracket on the set of length--twist coordinates, and we identified the canonical twist coordinates with geometric structures on pair-of-pant decompositions of Riemann surfaces $\Sigma_{g,s}$ with $s>0$ holes. Of course, the proof based on the Goldman bracket remains valid in the case of smooth surfaces $(g\ge2,s=0)$, but in this case we are lacking the first ingredient---the Fock bracket and the shear coordinate description, so we are lacking the cluster algebra description of the corresponding Teichm\"uller spaces. We generalized Fenchel--Nielsen coordinates to $\Sigma_{g,s,n}$---Riemann surfaces with $n>0$ bordered cusps on boundary components. Note that particular cases of such surfaces played a pivotal role in the description of monodromies of Painlev\'e equations in \cite{CMR2}. It is therefore natural to expect that results of this paper will find an application in the description of symplectic structures of the corresponding manifolds. In particular, the example in Fig.~\ref{fi:n-cusps} in the case of two cusps on the boundary component describes the monodromy manifold of the Painlev\'e~V equation .

A perspective direction of development of the method of this paper is related to dynamics of continued fractions based on Markov triples and their higher-dimensional generalizations. This addresses a hard problem of finding a way to ``stabilize'' geometric objects underlying a continued-fraction decompositions of projectivized real numbers in order to produce a variant of Poisson or quantum Thurston theory. The first attempt in this direction was done in \cite{ChP1}, and quite recently other methods emerged (see \cite{Ovs},\cite{BVes}) based on representations of $SL(2,\mathbb C)$-monodromies of Fuchsian equations.

\section*{Acknowledgements}
I am grateful to Ezra Getzler discussions with whom initiated this project. I am also grateful to my coauthors, Marta Mazzocco, Volodya Rubtsov and Misha Shapiro, for encouragement and the useful discussion.

\setcounter{section}{0}
\appendix{Wolpert's form and the Goldman bracket}
\setcounter{equation}{0}
Consider two oriented closed geodesics $\gamma_1$ and $\gamma_2$ intersecting at a point $P\in \gamma_1\# \gamma_2$. They may have other intersection points. We lift this pattern to the whole Poincar\'e upper half-plane: to every closed geodesic $\gamma$ we there set into correspondence a hyperbolic element $S_\gamma \in SL(2,\mathbb R)$ whose {\em invariant axis}, i.e., a unique geodesic that is invariant under the action of $S_\gamma$, becomes the closed geodesic upon identification by the action of $S_\gamma$, which acts as a shift by $\ell_\gamma$ along this invariant axis. Endpoints on the absolute of this axis are two distinct stable points of the hyperbolic element $S_\gamma$. Consider two elements $S_{\gamma_1}$ and $S_{\gamma_2}$. Let the invariant axis of $S_{\gamma_1}$ be the vertical half-line starting at the origin. Then the action of $S_{\gamma_1}$ is a dilatation, $z\to Rz$, $R=e^{\ell_{\gamma_1}}>1$. Let  stable points of the second element $S_{\gamma_2}$ be $-x_1$ and $x_2$ (we assume $x_i>0$), see the figure below. Then the (Euclidean) height $h$ of the axes intersection point is related to $x_i$ as $h^2=x_1x_2$, and the $SL(2,\mathbb R)$ element $S_{\gamma_2}$ has the form $z\to \frac{\alpha z-h^2}{z-\beta}$ with arbitrary $\alpha$, $\beta$. Consider now the composition $S_{\gamma_1}\circ S_{\gamma_2}:\ z\to R \frac{\alpha z-h^2}{z-\beta}= \frac{R\alpha z- Rh^2}{z-\beta}$. By the same consideration, the intersection point of the invariant axis of this element with that of $S_{\gamma_1}$ is now at height $h'=\sqrt{R}h$, which means that the geodesic distance between these two points (along the geodesic $\gamma_1$) is exactly $\ell_{\gamma_1}/2$.
 $$
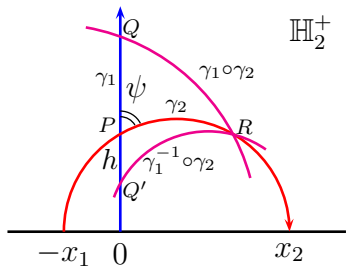
\begin{pspicture}(-5,-2)(5,1.5){\psset{unit=1.5}
\rput(0,0){
%
%\psclip{\pscircle[linewidth=1.5pt, linestyle=dashed](0,0){1}}
%
\rput(0,0){\psline[linewidth=1pt](-1,-1)(2,-1)}
\rput(0,0){\psline[linewidth=1pt, linecolor=blue]{->}(0,-1)(0,1)}
\rput(0,0){\psarc[linewidth=1pt,linecolor=red]{<-}(0.5,-1){1}{0}{180}}
\rput(0,0){\psarc[linewidth=0.5pt](0,-0.13){0.15}{30}{90}}
\rput(0,0){\psarc[linewidth=0.5pt](0,-0.13){0.2}{27}{90}}
%\endpsclip
\rput(0,0){\psarc[linewidth=1pt,linecolor=magenta](-0.625,-1){1.84}{15}{80}}
\rput(0,0){\psarc[linewidth=1pt,linecolor=magenta](0.78,-1){0.89}{55}{160}}
\rput(0.05,0.14){\makebox(0,0)[lb]{$\psi$}}
\rput(-0.5,-1.1){\makebox(0,0)[ct]{$-x_1$}}
\rput(0,-1.1){\makebox(0,0)[ct]{$0$}}
\rput(1.5,-1.1){\makebox(0,0)[ct]{$x_2$}}
\rput(-0.02,-0.35){\makebox(0,0)[rc]{\small$h$}}
\rput(-0.02,0.35){\makebox(0,0)[rc]{\tiny$\gamma_1$}}
\rput(0.5,0.05){\makebox(0,0)[cb]{\tiny$\gamma_2$}}
\rput(0.7,0.35){\makebox(0,0)[lb]{\tiny$\gamma_1{\circ}\gamma_2$}}
\rput(0.2,-0.25){\makebox(0,0)[lt]{\tiny$\gamma_1^{-1}{\circ}\gamma_2$}}
\rput(-0.04,-0.12){\makebox(0,0)[rb]{\tiny$P$}}
\rput(0.02,0.72){\makebox(0,0)[lb]{\tiny$Q$}}
\rput(1.5,0.62){\makebox(0,0)[lb]{$\mathbb H_2^+$}}
\rput(0.02,-0.55){\makebox(0,0)[lt]{\tiny$Q'$}}
\rput(1.02,-0.13){\makebox(0,0)[lb]{\tiny$R$}}
}
}
\end{pspicture}
$$
We therefore have the following general fact: given two hyperbolic elements $S_{\gamma_1}$ and $S_{\gamma_2}$ of a Fuchsian group $\Delta_{g,s}\subset PSL(2,\mathbb R)$ such that the corresponding closed geodesics $\gamma_1$ and $\gamma_2$ intersect at a point $P\in \gamma_1\# \gamma_2$, the closed geodesic $\gamma_1\circ_P \gamma_2$ (and $\gamma_1^{-1}\circ_P \gamma_2$) intersects the geodesics $\gamma_1$ and $\gamma_2$ at the respective points $Q$ and $R$ (at the respective points $Q'$ and $R'$) situated at geodesic distances that are exactly halves of the lengths of the corresponding closed geodesics, $|PQ|=|PQ'|=\ell_{\gamma_1}/2$ and $|PR|=|PR'|=\ell_{\gamma_2}/2$. Note that we then have that $|RQ|=\ell_{\gamma_1\circ_P \gamma_2}/2$ and $|RQ'|=\ell_{\gamma_1^{-1}\circ_P \gamma_2}/2$.

Let $\psi$ be the angle between $\gamma_1$ and $\gamma_2$ at the intersection point $P$. For a triangle $\{P,Q,R\}$ we have the hyperbolic cosine formula
$$
\cosh(|QR|)=-\sinh(|PQ|)\sinh(|PR|)\cos\psi +\cosh(|PQ|)\cosh(|PR|)
$$
applying which to two triangles $\{P,Q,R\}$  and $\{P,Q',R\}$ we obtain (replacing $\psi$ by $\pi-\psi$  in the second triangle)
\begin{align*}
\cosh\bigl(\ell_{\gamma_1\circ_P \gamma_2}\bigr)&=-\sinh(\ell_{\gamma_1}/2)\sinh(\ell_{\gamma_2}/2)\cos\psi +\cosh(\ell_{\gamma_1}/2)\cosh(\ell_{\gamma_2}/2)\\
\cosh\bigl(\ell_{\gamma_1^{-1}\circ_P \gamma_2}\bigr)&=\sinh(\ell_{\gamma_1}/2)\sinh(\ell_{\gamma_2}/2)\cos\psi +\cosh(\ell_{\gamma_1}/2)\cosh(\ell_{\gamma_2}/2).
\end{align*}
Adding the above two equalities we obtain the classical skein relation
$$
G_{\gamma_1\circ_P \gamma_2}+G_{\gamma_1^{-1}\circ_P \gamma_2}=G_{\gamma_1}G_{\gamma_2},
$$
and subtracting these relations, we obtain
$$
2\sinh(\ell_{\gamma_1}/2)\sinh(\ell_{\gamma_2}/2)\cos\psi=\cosh\bigl(\ell_{\gamma_1^{-1}\circ_P \gamma_2}\bigr)-\cosh\bigl(\ell_{\gamma_1\circ_P \gamma_2}\bigr),
$$
which generates the Wolpert formula out of the Goldman bracket: we can write the right-hand side as the Poisson bracket $\{G_{\gamma_1},G_{\gamma_2}\}=\sinh(\ell_{\gamma_1}/2)\sinh(\ell_{\gamma_2}/2)\{\ell_{\gamma_1},\ell_{\gamma_2}\}$, which gives that $\{\ell_{\gamma_1},\ell_{\gamma_2}\}=2\cos\psi$.

\end{document}